\documentclass[12pt, leqno]{amsart}

\usepackage{graphicx}
\usepackage{amssymb,amsmath,amsthm,amscd}
\usepackage{mathrsfs}
\usepackage{enumerate}
\usepackage{upgreek}
\usepackage{lscape}
\usepackage[usenames,dvipsnames]{color}
\usepackage[colorlinks=true, pdfstartview=FitV,
 linkcolor=blue,citecolor=blue,urlcolor=blue]{hyperref}
\usepackage{tikz}
\usepackage{xspace}
\usepackage{marginnote}
\usepackage{verbatim} 
 
\usepackage{framed}

\usepackage[all]{xy}

\allowdisplaybreaks[3]

\newcommand{\nc}{\newcommand}
\numberwithin{equation}{section}

\nc{\hs}{\hspace*}
\nc{\vs}{\vspace}
\nc{\ms}{\mspace}
\nc{\st}[1]{\{{#1}\}}

\nc{\qR}[1]{\ttq_{\mspace{-2mu}\raisebox{-.8ex}{${\scriptstyle{#1}}$}}}

\theoremstyle{plain}
\newtheorem{lemma}{Lemma}[section]
\newtheorem{proposition}[lemma]{Proposition}
\newtheorem{theorem}[lemma]{Theorem}

\newtheorem{corollary}[lemma]{Corollary}
\newtheorem{conjecture}{Conjecture}
\newtheorem*{convention}{Convention}

\newtheorem{thmx}{Theorem}

\theoremstyle{definition}
\newtheorem{remark}[lemma]{Remark}
\newtheorem{example}[lemma]{Example}

\newtheorem{definition}[lemma]{Definition}

\newcommand{\srt}[1]{{\langle {#1} \rangle}}
\newcommand{\ssrt}[1]{{( #1 )}}

\newcommand{\sseq}[1]{\boldsymbol{(}#1\boldsymbol{)}}

\newcommand{\tauQ}
{\tau_{\mspace{-2mu}\raisebox{-.4ex}{\scalebox{.6}{$Q$}}}}

\newcommand{\cox}[1]{ \tau_{\mspace{-2mu}\raisebox{-.4ex}{\scalebox{.6}{$#1$}}}}

\newcommand{\tauQp}
{\tau_{\mspace{-2mu}\raisebox{-.4ex}{\scalebox{.6}{$Q'$}}}}

\newcommand{\phiQ}
{\phi_{\mspace{-2mu}\raisebox{-.4ex}{\scalebox{.6}{$Q$}}}}

\newcommand{\phib}[1]{ \phi_{\mspace{-2mu}\raisebox{-.4ex}{\scalebox{.6}{$#1$}}}}

\renewcommand{\le}{\leqslant}
\renewcommand{\ge}{\geqslant}
\renewcommand{\preceq}{\preccurlyeq}

\newcommand{\ssqcup}{\mathop{\mbox{\normalsize$\bigsqcup$}}\limits}

\newcommand{\qm}{  \Updelta}
\newcommand{\seteq}{\mathbin{:=}}

\newcommand{\soplus}{\mathop{\mbox{\normalsize$\bigoplus$}}\limits}

\newcommand{\Rep}{{\rm Rep}}

\newcommand{\g}{\mathfrak{g}}
\newcommand{\h}{\mathfrak{h}}
\newcommand{\n}{\mathfrak{n}}

\newcommand{\R}{\mathbb{R}}
\newcommand{\C}{\mathbb{C}}
\newcommand{\Q}{\mathbb{Q}}
\newcommand{\Z}{\mathbb{Z}\ms{1mu}}

\newcommand{\al}{{\ms{1mu}\alpha}}
\newcommand{\ep}{\epsilon}
\newcommand{\la}{\lambda}
\newcommand{\be}{{\ms{1mu}\beta}}
\newcommand{\ga}{\gamma}
\newcommand{\ve}{\varepsilon}
\newcommand{\La}{\Lambda}

\newcommand{\wt}{{\rm wt}}

\newcommand{\de}{\mathfrak{d}}
\newcommand{\Dynkin}{\triangle}
\newcommand{\bDynkin}{\mathop{\mathbin{\mbox{\large $\blacktriangle$}}}}
\newcommand{\lf}{[\hspace{-0.3ex}[}
\newcommand{\rf}{]\hspace{-0.3ex}]}

\newcommand{\lt}{{\mathrm{l}}}
\newcommand{\rt}{{\mathrm{r}}}
\newcommand{\lwt}{\wt_\lt}
\newcommand{\rwt}{\wt_\rt}

\newcommand*\circled[1]{ \fontsize{6}{6}\selectfont \tikz[baseline=(char.base)]{
  \node[shape=circle,draw,inner sep=0.4pt] (char) {#1};} \fontsize{12}{12}\selectfont }

\newcommand*\bcircled[1]{\fontsize{6}{6}\selectfont \tikz[baseline=(char.base)]{
    \node[shape=circle, fill= black, draw=black, text=white,  inner sep=0.4pt] (char) {#1};} \fontsize{12}{12}\selectfont}

\newcommand*\rcircled[1]{\fontsize{6}{6}\selectfont \tikz[baseline=(char.base)]{
    \node[shape=circle, fill= gray, draw=black, text=white, inner sep=0.4pt] (char) {#1};} \fontsize{12}{12}\selectfont}

\newcommand{\Hom}{\operatorname{Hom}}

\newcommand{\res}{\mathrm{res}}

\newcommand{\tB}{\widetilde{B}}
\newcommand{\tsfB}{\ms{1mu}\widetilde{\sfB}}
\newcommand{\tX}{\widetilde{X}}

\newcommand{\tw}{{\widetilde{w}}}
\newcommand{\tm}{\widetilde{m}}
\newcommand{\tsfC}{{\widetilde{\sfC}}}

\newcommand{\tuC}{{\widetilde{\usfC}}}
\newcommand{\tuB}{\ms{1mu}\widetilde{\usfB}}
\newcommand{\tfb}{\widetilde{\sfb}}

\newcommand{\tde}{\widetilde{\de}}

\newcommand{\hPhi}{\widehat{\Phi}}

\newcommand{\hg}{\widehat{\g}}

\newcommand{\osfB}{\overline{\mathsf{B}}}
\newcommand{\usfB}{\underline{\mathsf{B}}}

\newcommand{\usfC}{\underline{\mathsf{C}}}
\newcommand{\uw}{{\underline{w}}}
\newcommand{\ucalN}{\ms{1mu}\underline{\calN}}

\newcommand{\utw}{\underline{\tw}}

\newcommand{\frakK}{\mathfrak{K}}

\newcommand{\sfC}{\mathsf{C}}
\newcommand{\sfS}{\mathsf{S}}
\newcommand{\sfc}{\mathsf{c}}
\newcommand{\sfB}{\mathsf{B}}

\newcommand{\sfD}{\mathsf{D}}
\newcommand{\sfP}{\mathsf{P}}
\newcommand{\sfQ}{\mathsf{Q}}
\newcommand{\sfW}{\mathsf{W}}

\newcommand{\sfh}{\mathsf{h}}

\newcommand{\sfn}{\mathsf{n}}

\newcommand{\sfb}{\mathsf{b}}
\newcommand{\sfg}{\mathsf{g}}

\newcommand{\bbA}{\mathbb{A}}

\newcommand{\bfB}{\mathbf{B}}

\newcommand{\calD}{\mathcal{D}}

\newcommand{\calA}{\mathcal{A}}
\newcommand{\calK}{\mathcal{K}}
\newcommand{\calI}{\mathcal{I}}

\newcommand{\calX}{\mathcal{X}}
\newcommand{\calY}{\mathcal{Y}}

\newcommand{\calN}{\mathcal{N}}
\newcommand{\calO}{\mathcal{O}}
\newcommand{\calT}{\mathcal{T}}
\newcommand{\calW}{\mathcal{W}}

\newcommand{\scrQ}{\mathscr{Q}}
\newcommand{\scrC}{\mathscr{C}}
\newcommand{\scrA}{\mathscr{A}}
\newcommand{\scrS}{\mathscr{S}}

\newcommand{\ttD}{\mathtt{D}}

\newcommand{\ttt}{\mathtt{t}}
\newcommand{\ttq}{{\ms{1mu}\mathtt{q}\ms{1mu}}}
\newcommand{\To}[1][{\hspace{2ex}}]{\xrightarrow{\,#1\,}}

\newlength{\mylength}
\newcommand{\hDynkin}{{\widehat{\Dynkin}}}
\newcommand{\hbDynkin}{\widehat{\bDynkin}}

\newcommand{\Oint}{\calO_{{\rm int}}}
\newcommand{\diag}{\mathrm{diag}}
\newcommand{\rmQ}{\mathrm{Q}}

\newcommand{\Qd}{\scrQ}
\newcommand{\cm}{\sfC}
\newcommand{\wl}{\sfP}
\newcommand{\rl}{\sfQ}
\newcommand{\sr}{\Pi}
\newcommand{\cwl}{\wl^\vee}
\newcommand{\scr}{\Pi^\vee}
\newcommand{\weyl}{\sfW}
\newcommand{\lan}{\langle}
\newcommand{\ran}{\rangle}
\newcommand{\Ug}{U_\nu(\g)}
\newcommand{\Ag}{A_\nu(\g)}
\newcommand{\An}{A_\nu(\n)}
\newcommand{\Upg}{U^+_\nu(\g)}

\newcommand{\sprt}[1]{\scalebox{0.95}{$(#1)$}}
\newcommand{\isoto}[1][]{\mathop{\xrightarrow%
[{\raisebox{.3ex}[0ex][.3ex]{$\scriptstyle{
#1}$}}]%
{{\raisebox{-.6ex}[0ex][-.6ex]{$\mspace{2mu}\sim\mspace{2mu}$}}}}}

\newcommand{\ee}{\end{enumerate}}
\newcommand{\ben}{\begin{enumerate}[{\rm (1)}]}
\newcommand{\bnum}{\begin{enumerate}[{\rm (i)}]}
\newcommand{\bna}{\begin{enumerate}[{\rm (a)}]}
\newcommand{\bc}{\begin{cases}}
\newcommand{\ec}{\end{cases}}

\newenvironment{myequation}
{\relax\setlength{\arraycolsep}{1pt}\begin{eqnarray}}
{\end{eqnarray}}
\newenvironment{myequationn}
{\relax\setlength{\arraycolsep}{1pt}\begin{eqnarray*}}
{\end{eqnarray*}}

\nc{\eq}{\begin{myequation}}
\nc{\eneq}{\end{myequation}}
\nc{\eqn}{\begin{myequationn}}
\nc{\eneqn}{\end{myequationn}}
\nc{\cl}{\colon}
\nc{\ake}[1][1ex]{\rule[-#1]{0ex}{1ex}}
\nc{\akew}[1][1ex]{\rule[-1ex]{#1}{0ex}}
\nc{\akeu}[1][1ex]{\rule[#1]{0ex}{1ex}}
\nc{\id}{\mathrm{id}}
\nc{\bl}{\bigl(}
\nc{\br}{\bigr)}
\nc{\qt}[1]{\quad\text{#1}}
\nc{\qtq}[1][{and}]{\quad\text{{#1}}\quad}
\nc{\eqs}[1][{*}]{\mathrel{\hs{.4ex}\raisebox{-.8ex}[0ex][.3ex]{$\scriptstyle{#1}$}\hs{-1.9ex}=\hs{.3ex}}}
\nc{\snoi}{\smallskip \noindent}
\nc{\mnoi}{\medskip \noindent}
\nc{\Mat}{\mathrm{Mat}}
\nc{\ol}{\overline}
\nc{\ang}[1]{\langle{#1}\rangle}
\nc{\ba}{\begin{array}}
\nc{\ea}{\end{array}}
\nc{\noi}{\noindent}
\nc{\Cor}{\begin{corollary}}
\nc{\encor}{\end{corollary}}
\nc{\up}{\mathrm{up}}
\nc{\Lemma}{\begin{lemma}}
\nc{\enlemma}{\end{lemma}}
\nc{\Prop}{\begin{proposition}}
\nc{\enprop}{\end{proposition}}
\nc{\ro}{{\rm(}}
\nc{\rfm}{{\rm)}}
\nc{\Aut}{\operatorname{Aut}}
\nc{\Bl}{\Bigl(}
\nc{\Br}{\Bigr)}
\nc{\Proof}{\begin{proof}}
\nc{\QED}{\end{proof}}
\nc{\wh}{\widehat}
\nc{\ws}{\wh{s}}
\nc{\wcox}[1][Q]{\wh{\tau}_{\mspace{-2mu}\raisebox{-.4ex}{\scalebox{.6}{$#1$}}}}
\nc{\rtl}{\mathsf{Q}}
\nc{\Ga}{\Gamma}
\nc{\teta}{\widetilde{\eta}}
\nc{\nn}{\nonumber}
\nc{\lec}[1][2]{\le_{#1}\ms{1mu}}
\nc{\LHS}{the left hand side\xspace}
\nc{\RHS}{the right hand side\xspace}

\newcounter{myce}
\newcounter{mycs}

\usepackage[colorinlistoftodos]{todonotes}

\title[The $(q,t)$-Cartan matrix specialized at $q=1$]{The $(q,t)$-Cartan matrix specialized at $q=1$ \\ and its applications\protect}

\date{\today}

\author[M.\ Kashiwara]{Masaki Kashiwara}
\thanks{The research of M.\ Kashiwara
was supported by Grant-in-Aid for Scientific Research (B)
20H01795, Japan Society for the Promotion of Science.}
\address[M.\ Kashiwara]{
Kyoto University Institute for Advanced Study,
Research Institute for Mathematical Sciences, Kyoto University,
Kyoto 606-8502, Japan \& Korea Institute for Advanced Study, Seoul 02455, Korea }
\email{masaki@kurims.kyoto-u.ac.jp}

\author[S.-j.~Oh]{Se-jin Oh}
\thanks{ The research of S.-j.\ Oh was supported by the Ministry of Education of the Republic of Korea and the National Research Foundation of Korea (NRF-2022R1A2C1004045).}
\address[S.-j.~Oh]{Ewha Womans University Seoul, 52 Ewhayeodae-gil, Daehyeon-dong, Seodaemun-gu, Seoul, South Korea}
\email{sejin092@gmail.com}


\begin{document}

\begin{abstract}
The $(q,t)$-Cartan matrix specialized at $t=1$, usually called the \emph{quantum Cartan matrix}, has deep connections with (i) the representation theory
of its untwisted quantum affine algebra, and (ii) quantum unipotent coordinate algebra, root system and quantum cluster algebra of \emph{skew-symmetric type}.
In this paper, we study
the $(q,t)$-Cartan matrix specialized at $q=1$, called the \emph{$t$-quantized Cartan matrix}, and investigate the relations with (ii$'$) its corresponding unipotent quantum coordinate algebra,
root system and quantum cluster algebra of \emph{skew-symmetrizable type}.
\end{abstract}

\setcounter{tocdepth}{1}

\maketitle \tableofcontents

\section{Introduction}

\subsection{$(q,t)$-Cartan matrix $\cm(q,t)$ }
Let $\g$ be a  finite-dimensional simple Lie algebra, $\cm=(\sfc_{i,j})_{i,j\in I}$ its Cartan matrix,
$\Dynkin$ its Dynkin Diagram, $\Phi^+$ its set of positive roots,
 and $U_\nu'(\widehat{\g})$ the untwisted quantum affine algebra associated to $\g$. When $\Dynkin$ is non simply-laced, we sometimes use $\bDynkin$ to
distinguish with $\Dynkin$  of an arbitrary type.
In \cite{FR99}, Frenkel-Reshetikhin introduced  a two-parameter
deformation $\cm(q,t)$ of $\cm$ to define
the deformation of $\calW$-algebra   $\calW_{q,t}(\g)$  of $\g$. We
call the $\cm(q,t)$ the \emph{$(q,t)$-Cartan matrix of $\g$}. Then
  {\rm (i)} it is  proved in  \cite{FR99B} that
the limit $t \to  1$ of $\calW_{q,t}(\g)$ recovers
 the commutative
Grothendieck ring of the category of finite-dimensional modules
$\scrC_{\widehat{\g}}$ over $U'_\nu (\widehat{\g})$, and  {\rm (ii)} it is expected in \cite{FH11} that 
 the limit $q \to  {\rm
exp}(\pi i/r )$ of $\calW_{q,t}(\g)$ recovers  the one of
$\scrC_{{}^L(\widehat{\g})}$  over the Langlands dual $U'_\nu
({}^L(\widehat{\g}) )$ of $U'_\nu (\widehat{\g})$  , where $r$ denotes the order of non-trivial Dynkin diagram automorphism $\sigma$ of simply-laced type Dynkin diagram, yielding the Dynkin diagram $\Dynkin$ of $\g$:
\begin{align} \label{eq: recovering rep theory}
\xymatrix@R=5ex@C=10ex{   K\left(\scrC_{\widehat{\g}} \right) \ar@/^1.5pc/@{<-->}[rr]|{\text{interpolation}} & \mathcal{W}_{q,t}(\g) \ar@{->}[l]^{ 1 \gets t }
\ar@{->}[r]_{ q \to {\rm exp}(\pi i /r)  } & K\left(\scrC_{{}^L\widehat{\g}}\right)}
\end{align}

In this paper, we study the specialization of $\cm(q,t)$
at $q=1$, which is not well-investigated, while the other specialization of $\cm(q,t)$ at $t=1$ is intensively investigated since late 1990's, in various view points.

\subsection{Quantum Cartan matrix $\cm(q)$ } The specialization of $\cm(q,t)$ at $t=1$, denoted by $\cm(q)$ and called the \emph{quantum Cartan matrix}, and its inverse $\widetilde{\cm}(q) \seteq \cm(q)^{-1}$ are known as the    \emph{controller} of the representation theory of $U'_\nu (\widehat{\g})$,
and has been extensively studied.
To name a few,  (i) we can read the denominator formulas of the normalized $R$-matrices between \emph{fundamental modules}
from the Laurent expansions of $\widetilde{\cm}(q)_{i,j}$'s at $q=0$, which \emph{determine} whether a given tensor product of fundamental modules is simple or not
(\cite{DO94, KKKII, Oh15, Fuj19, OS19B, OS19}),
(ii) ${\cm}(q)$ and $\widetilde{\cm}(q)$ are used as key ingredients for constructing  the \emph{quantum Grothendieck ring}  $K_\ttt(\scrC_{\widehat{\g}})$ of  $\scrC_{\widehat{\g}}$  (\cite{Nak04,VV02,H04,HL15}) and   the   $q$-character theory for $\scrC_{\widehat{\g}}$
(\cite{FR99,FR99B,FM01}),
(iii) the quantum Grothendieck ring $K_\ttt(\scrC^0_{\widehat{\g}})$ is isomorphic to   a   certain $\Z[\ttt^{\pm \frac{1}{2}}]$-algebra inside the quantum torus $\calY_{\ttt,\g}$, which is constructed from $\widetilde{\cm}(q)$ (see \cite{HL16,KKOP21}).
Here $\scrC^0_{\widehat{\g}}$
is the \emph{skeleton} subcategory of $\scrC_{\widehat{\g}}$,  also referred to as the Hernandez-Leclerc category.

Interestingly enough, the computation of $\widetilde{\cm}(q)$ of type $\g$ is related to the simple Lie algebra $\sfg$ of simply-laced type as below:
$$ (\g,\sfg) =  (A_n,A_{n}), \ (D_n,D_{n}), \ (E_{6,7,8},E_{6,7,8}), \  (B_n,A_{2n-1}),
\ (C_n,D_{n+1}),\  (F_4,E_6),\  (G_2,D_4).$$

It is  proved in \cite{HL15} that, when $\g=\sfg$ is of simply-laced type,
Laurent expansions of $\widetilde{\cm}(q)_{i,j}$'s at $q=0$ of type $\g$ can be computed by using any
Dynkin quiver $Q$ of the same type and the symmetric bilinear form $( \ , \ )$ on the root lattice $\sfQ_\sfg$.

As a generalization of \cite{HL15}, it is proved in
\cite{FO21} that Laurent expansions of $\widetilde{\cm}(q)_{i,j}$'s at $q=0$ of any finite type
$\g$  can be  computed by using any \emph{$\rmQ$-datum associated to $\g$} and
$( \ , \ )$ on $\sfQ_\sfg$.
We remark here that Laurent expansions of $\widetilde{\cm}(q)_{i,j}$'s at $q=0$ exhibit the remarkable periodicity and the positivity
related to the \emph{dual} Coxeter number $\sfh^\vee$ of $\g$

The notion of $\rmQ$-datum is introduced in~\cite{FO21}
 and can be understood as a generalization of Dynkin quiver.
The constituents and properties of $\rmQ$-datum associated to $\g$ are briefly summarized
as follows: \bna
\item A $\rmQ$-datum $\scrQ=(\Dynkin,\sigma,\xi)$ consists of (1) the Dynkin diagram $\Dynkin_\sfg$ of type $\sfg$,  (2)  the Dynkin diagram automorphism $\sigma$
yielding $\Dynkin_\g$ as an orbit of $\Dynkin_\sfg$ via $\sigma$, and (3) a height function $\xi$ on $\Dynkin$  satisfying certain axioms.
\item For each $\rmQ$-datum $\scrQ=(\Dynkin,\sigma,\xi)$, there exists a unique (generalized)-Coxeter element $\cox\scrQ$
in $\weyl_\sfg \rtimes \lan \sigma \ran$
satisfying certain properties.
\item
When $\g=\sfg$ and hence $\sigma={\rm id}$, the notion of $\rmQ$-datum $\scrQ=(\Dynkin,\sigma,\xi)$ is equivalent to the notion of Dynkin quiver $Q=(\Dynkin,\xi)$  of  type $ADE$.
\ee

\smallskip

When $\g=\sfg$ is of simply-laced type, we can consider the
path algebra $\C Q$ of a Dynkin quiver $Q=(\Dynkin_\sfg,\xi)$. Then it is well-known that the Auslander-Reiten(AR) quiver
$\Gamma_Q$ of $\C Q$  realizes the convex partial order
$\prec_Q$ on $\Phi^+_\sfg$ coming from the \emph{$Q$-adapted}
commutation class $[Q]$ of the longest element $w_0$ of
$\weyl_\sfg$.  On the other hand, $\Gamma_Q$ can be understood as a
\emph{heart} of the AR-quiver $\hDynkin=(\hDynkin_0,\hDynkin_1)$ of the
derived category $D^b({\rm Rep}(\C Q))$, in which $\hDynkin$ is
referred to as the \emph{repetition quiver}. In aspect of
combinatorics, (i) the set of vertices $\hDynkin_0$  is in one to
one correspondence with $\Phi^+_\sfg \times \Z$,
and  (ii) $\hDynkin$ satisfies
\emph{$\hDynkin$-additive property} as equations in the root lattice
$\sfQ_\sfg$ (see~\eqref{eq: classic additive}). In \cite{HL15}, Hernandez-Leclerc defined the \emph{heart}
subcategory $\scrC_Q$ of
$\scrC_{\widehat{\g}}=\scrC_{\widehat{\sfg}}$ by using the \emph{coordinate system} via $I \times \Z$ for vertices in  $\hDynkin$ and $\Gamma_Q$ (see~\eqref{eq: rep quiver} and~\eqref{eq: well-known}).
Then they proved
that the quantum  Grothendieck ring $K_\ttt(\scrC_Q)$ of $\scrC_Q$ is
isomorphic to the unipotent quantum coordinate algebra
$A_\nu(\mathsf{n})$ of the quantum group $U_\nu(\sfg)$. Here
$K_\ttt(\scrC_Q)$ is contained in a certain sub-torus
$\calY_{\ttt,Q}$ of $\calY_{\ttt,\sfg}$ determined by the coordinate
system of $\Gamma_Q$.

\smallskip

By understanding $\Gamma_Q$ as the Hasse quiver of $\prec_{Q}$, for each
$\rmQ$-datum $\scrQ=(\Dynkin_\sfg ,\sigma \ne {\rm id},\xi)$ associated to non simply-laced $\g$ ($\ne \sfg$), (i) the
(combinatorial) AR-quiver $\Gamma_\scrQ$,  the repetition quiver $\hDynkin^\sigma$ and the heart subcategory
$\scrC_\scrQ$ of $\scrC_{\widehat{\g}}$ are defined in~\cite{FO21,OS19B}
by developing $\scrQ$-adapted commutation class $[\scrQ]$ of $w_0 \in
\weyl_\sfg$,
(ii) it is proved that $\hDynkin^\sigma$ satisfies the \emph{$\hDynkin^\sigma$-additive property} and the set of vertices of $\hDynkin^\sigma$
is also in one to one correspondence with $\Phi_\sfg^+ \times \Z$,
 and (iii) it is proved in~\cite{KO18,HO19,FHOO} that the quantum Grothendieck ring $K_\ttt(\scrC_\scrQ) \subsetneq \calY_{\ttt,\scrQ}$ is also isomorphic to
$A_\nu(\mathsf{n})$ of $U_\nu(\sfg)$.   Furthermore, it is proved
in~\cite{FHOO} that the quantum Grothendieck ring
$K_\ttt(\scrC_{\widehat{\sfg}}^0)$ of $\scrC_{\widehat{\sfg}}^0$ and the one
$K_\ttt(\scrC_{\widehat{\g}}^0)$ of  $\scrC_{\widehat{\g}}^0$
  are isomorphic as an algebra,  preserving \emph{simple $(q,t)$-characters}. Along the proofs in
\cite{FHOO}, the isomorphism between the quantum torus
$\calY_{\ttt,\scrQ}$ and $\calT_{\nu,\scrQ}$ containing
$A_\nu(\mathsf{n})$ played important roles.

\smallskip

A \emph{cluster algebra} $\calA$, introduced by Fomin-Zelevinsky in \cite{FZ02} is   a  commutative $\Z$-algebra contained in the torus $\Z[X_{k}^{\pm 1} \mid k \in \mathsf{K}]$.
 By Berenstein and Zelevinsky \cite{BZ05},  the notion of cluster algebras
is generalized to the non-commutative one, \emph{quantum cluster
algebra}  $\calA_\upnu$, which is contained in the quantum torus
$\Z[\upnu^{1/2}][\widetilde{X}_{k}^{\pm 1} \mid k \in \mathsf{K}]$. From their
introductions, numerous connections and applications have been
discovered in various fields of mathematics including the
representation theory of quantum affine algebras.

Note that  the Grothendieck ring $K(\scrC^0_{\widehat{\g}})$ of $\scrC^0_{\widehat{\g}}$ has a
$\Uplambda$-cluster algebra structure (see \cite{KKOP20} for its definition).
Moreover, for any $\rmQ$-datum $\Qd$, it is proved in \cite{KKOP21,FHOO} (see
also \cite{HL16}) that the heart subcategory $\scrC_{\Qd}$ of
$\scrC_\g$ provides a monoidal categorification of the quantum cluster algebra $K_\ttt(\scrC_{\Qd})\simeq A_\nu(\sfn)$.

On the other hand, it is proved in
\cite{GLS13,GY17} that the unipotent quantum coordinate algebra
$A_\nu(\n)$ of the quantum group $U_\nu(\g)$ $(\g$ need not be the same as $\sfg)$ has a quantum cluster
algebra structure of skew-symmetrizable type, whose initial quantum seed  arises from the
combinatorics of $\weyl_\g$ and  the symmetric bilinear form on the
root lattice $\sfQ_\g$.

\subsection{Brief summary and comments on $\sfC(q)$}
To sum up, we can conclude that each $\rmQ$-datum $\scrQ$ associated to $\g$ having $\Dynkin$
as a simply-laced type $\sfg$ and its related objects, such as $\widetilde{\cm}(q)$,
$\cox\scrQ$, $\Gamma_\scrQ$ and $\hDynkin^\sigma$, encode key
information of the representation theories of  $U_\nu'(\widehat{\g})$
and  $A_\nu(\mathsf{n})$ of
the quantum group $U_\nu(\sfg)$.
However, for $\g$ with $\g \ne\sfg$,  $A_\nu(\n)$ of
 $U_\nu(\g)$ looks not related to the representation theory of quantum affine algebras   nor   known   $\rmQ$-data and  their   related objects, as far as the authors understand at this moment.

\subsection{$t$-quantized Cartan matrix $\usfC(t)$}
In the definition of $\rmQ$-datum,  $\bDynkin$ of type BCFG does \emph{not} appear as a constituent, and
$\bDynkin$ do have only the trivial Dynkin diagram automorphism ${\rm id}$.
In this paper, we mainly consider a Dynkin quiver $Q=(\bDynkin,\xi)$  as a generalization of the notion of $\rmQ$-datum  (see Section~\ref{subsec: D quiver}).

From this paper, we start to explore the representation theories which are controlled by  the specialization  $\usfC(t) \seteq \cm(1,t)$ and its inverse $\tuC(t)$.
We call  $\usfC(t)$  \emph{$t$-quantized Cartan matrix}.
As a starting point, we investigate  the relations among Dynkin quivers $Q$, $A_\nu(\n)$ and  $\tuC(t)$ in this paper.
As $\tsfC(q)$ and objects induced from $\tsfC(q)$ were used for representation theories of  $U_q'(\widehat{\g})$ and  $A_{\nu}(\mathsf{n})$,
we construct various mathematical objects from $\tuC(t)$ and Dynkin quivers $Q$,
and study the newly introduced objects and their applications.

The main achievements of this paper can be summarized as follows:
\smallskip

\fbox{ \ \parbox{81ex}{
\ben
\item[{\rm (\textbf{A})}] We introduce quivers $\Gamma_Q$ and $\hDynkin$ for Dynkin quivers $Q$ of any finite type $\g$, study their properties and   prove that the Laurent expansions of $\tuC(t)_{i,j}$'s at $t=0$ can be computable by any of them.
\item[{\rm (\textbf{B})}] Using $\tuC(t)_{i,j}$'s, we construct the quantum torus $\calX_\ttq$ and prove that its $\ttq$-commutation relation is controlled
by the root system of $\g$.
\item[{\rm (\textbf{C})}] For each Dynkin quiver $Q$, we define the quantum sub-torus $\calX_{\ttq,Q}$ of $\calX_\ttq$  and show that
it is isomorphic to the quantum torus $\calT_{\nu,[Q]}$ containing $A_{\nu}(\mathfrak{n})$.
\item[{\rm (\textbf{D})}]  For any $Q$-adapted sequence $\tw$, we prove that the pair  $(\La^\tw,B^\tw)$ is compatible, which recovers the compatible pair in \cite{HL15} as a particular case, by using the $\hDynkin$-additive property and $\tuC(t)_{i,j}$'s.
\ee
} \ }

\smallskip

\noi
{\rm (\textbf{A})} With the  notion of Dynkin quiver $Q =(\Dynkin,\xi)$, including
BCFG-types,  we can define the notion of $Q$-adaptedness for reduced
expressions of elements in $\weyl_\Dynkin$. Then we can see that the
set of all $Q$-adapted reduced expressions of $w_0 \in
\weyl_\Dynkin$  forms a commutation class $[Q]$. In particular, we
can see that $[Q]$ of type $BCFG$ enjoys the almost same properties
of $[Q]$ in $\weyl$ of type $ADE$ (Theorem~\ref{thm: unified 1}).

Similarly to the combinatorial feature of $\hDynkin$ for ADE-type, we define the repetition quiver $\hbDynkin=(\hbDynkin_0,\hbDynkin_1)$
for every Dynkin quiver
$Q=(\bDynkin,\xi)$ and prove that
\bnum
\item $\hbDynkin_0$  is in one  to one correspondence with $\Phi^+_{\blacktriangle} \times \Z$ via the map $\phiQ$
(Theorem~\ref{thm: bijection new}),
\item it satisfies the  $\hbDynkin$-additive property   (Theorem~\ref{thm: additive}).
\ee
Then we define the heart subquiver $\Gamma_Q$ of  $\hbDynkin$, also called the \emph{AR-quiver of $Q$}, and prove that it realizes the convex partial order $\prec_{[Q]}$ on $\Phi^+_{\blacktriangle}$ (Proposition~\ref{prop: quiver iso}).

\smallskip

Recall that  the   Laurent  expansion of $\widetilde{\cm}(q)_{i,j}$ at $q=0$ of type $\g$ can be   computed   by using  an arbitrary
$\scrQ=(\Dynkin_\sfg,\sigma,\xi)$.
Generalizing this phenomenon,  we prove  the following theorem:

\begin{thmx}  [{Theorem \ref{thm: inv}, Corollary~\ref{cor: computation setting}, Theorem~\ref{thm: ADE denom},~\ref{thm: BC denom}}] \label{them: main 1.1}
The  Laurent  expansion of $\tuC(t)_{i,j}$ of type $\g$ at $t=0$ can be
computed by using any Dynkin quiver $Q=(\Dynkin,\xi)$ or
$(\bDynkin,\xi)$ of the same type. Furthermore,
\bna
\item \label{it: pp}  the periodicity and the positivity for $\tuC(t)$ hold,
\item we obtain the closed formulae for $\tuC(t)_{i,j}$ for all Cartan types.
\ee
\end{thmx}

We remark here that in \cite{FM21}, the inverse of the
$\cm(q,t)$ itself is described in terms of bigraded modules over the
generalized preprojective algebras in the sense of
Gei\ss-Leclerc-Schr\"oer \cite{GLS17}.   In particular,  ~\eqref{it: pp} in Theorem~\ref{them: main 1.1}
can be also obtained from~\cite[Corollary 4.14]{FM21} by specialization at $q=1$ (see~\S~\ref{subsubsec: Quantum cartan and two parameter}). 

\noi
{\rm (\textbf{B})}
As we mentioned above, $\widetilde{\cm}(q)$ had used for the
constructions of the quantum torus, $(q,t)$-character theory and the
quantum Grothendieck ring $K_{\ttt}(\scrC_{\widehat{\g}}^0)$ of
$\scrC^0_{\widehat{\g}}$.
In \cite{FM01,FR99,Nak04,VV02,H03,H04}, the
algebraic constructions of $K_{\ttt}(\scrC_{\widehat{\g}}^0)$ and
its torus $\calY_{\ttt,\g}$ are described in terms of
$\widetilde{\cm}(q)$.
More precisely, the
quantum torus $\calY_{\ttt,\g}$ is constructed  from the  Laurent
expansion of $\widetilde{\cm}(q)_{i,j}$ at $q=0$ and proved that
$K_{\ttt}(\scrC_{\widehat{\g}}^0)$ is isomorphic to the intersection
of $\ttt$-screening operators $S_{i,\ttt}$ $(i \in I)$ on
$\calY_{\ttt,\g}$.

As $\widetilde{\usfC}(t)$-analogues, we
construct a new quantum torus $\calX_{\ttq}$ (Definition~\ref{def:
new q torus})  and \emph{quantum virtual Grothendieck ring} $\frakK_{\ttq}$
(Definition~\ref{def: virtual quantum G ring}) $ \subsetneq \calX_{\ttq}$, which \emph{recovers}
$\calY_{\ttt,\g}^0$  and $K_{\ttt}(\scrC_{\widehat{\g}}^0)$
respectively, when $\g=\sfg$.

\begin{thmx} [{Theorem~\ref{thm: calN},  Proposition~\ref{prop: YA com}}]  \label{thm: main 1.2}
The $\ttq$-commutation relation of $\calX_{\ttq}$ of type $\g$ is \emph{controlled} by {\rm (i)} the symmetric bilinear form on $\sfQ_\g$ and {\rm (ii)} the bijection $\phiQ$ between
$\hDynkin_0$ and $\Phi^+_\g \times \Z$
for any Dynkin quiver $Q$ of the same type.
\end{thmx}

\smallskip

\noi
{\rm (\textbf{C})}
Note that, for any $\g$ and a commutation $[\uw_0]$ of $w_0 \in
\weyl_\g$, $A_\nu(\n)$ of $U_\nu(\g)$ is contained in the quantum
torus $\calT_{\nu,[\uw_0]}$ and generated by the certain set of
unipotent quantum minors in $A_\nu(\n)$.   More precisely,  it is isomorphic to the quantum cluster
algebra $\calA_\nu(\Lambda^{[\uw_0]},B^{[\uw_0]})$ contained in the quantum torus
$\calT_{\nu,[\uw_0]}$  whose initial quantum seed
$\scrS = \{ (\Lambda^{[\uw_0]},B^{[\uw_0]}), \{x_i\}_{1 \le i \le \ell(w_0)} \}$
is given by the combinatorics of $\weyl_\g$ and $\sfQ_\g$ \cite{GLS13,GY17}. In particular, if
$\weyl$ is of simply-laced type $\sfg$ and $[\uw_0]=[\scrQ]$ for a $\rmQ$-datum $\scrQ=(\Dynkin_\sfg,\sigma,\xi)$, it is proved in \cite{HL15,FHOO} that $\calT_{\nu,[\scrQ]}$ is isomorphic to the quantum subtorus $\calY_{\ttt,\scrQ}$ of $\calY_\ttt$.
We prove the following theorem by considering the Dynkin quiver $Q$ of type BCFG also:

\begin{thmx} [Theorem~\ref{thm: tori iso}] \label{main thm: 1.3}
 For any $\g$ and a Dynkin quiver $Q$ of type $\g$, we have an isomorphism $\Uppsi_Q$ between $\calT_{\nu,[Q]}$ and
$\calX_{\ttq,Q}$, where $\calX_{\ttq,Q}$ denotes the quantum subtorus  of $\calX_\ttq$ determined by the coordinate system of $\Gamma_Q$ inside $\hDynkin$.
\end{thmx}

By restricting the map $\Uppsi_Q$ to the subalgebra  $\calA_\nu(\Lambda^{[Q]},B^{[Q]})$ of $\calT_{\nu,[Q]}$,  we can obtain the subalgebra $\frakK_{\ttq,Q}$ of $\calX_{\ttq,Q}$. We expect that
$\frakK_{\ttq,Q}$ is a subalgebra of $\frakK_{\ttq}$ for every Dynkin quiver $Q$ of type $\g$. Note that $\frakK_{\ttq,Q}$ is isomorphic to the quantum Grothendieck ring $K_\ttt(\scrC_Q)$  of  $\scrC_Q$ when $Q$ is a Dynkin quiver of type $ADE$.

\smallskip

\noi
{\rm (\textbf{D})}
For a Dynkin quiver $Q=(\Dynkin_\sfg,\xi)$ of simply-laced type and
a sequence $\tw$ of indices of $I_\sfg$ satisfying certain condition in~\eqref{cond:din},
we can construct a pair of matrices $(\Lambda^{\tw},B^{\tw})$  by
extending the combinatorics of $\weyl_\sfg$ and $\sfQ_\sfg$.
In particular, if $\tw$ is reduced, then the pair $(\Lambda^{\tw},B^{\tw})$ is  \emph{compatible with an integer $2$},
 $\calA_\nu(\Lambda^{\tw},B^{\tw})$ is isomorphic to $A_\nu(\sfn(w)) \subseteq A_\nu(\sfn)$ and hence
$K_\ttt(\scrC_{\widehat{\sfg}}^\tw)$ for a certain subcategory $\scrC^\tw_{\widehat{\sfg}}$ of
$\scrC^0_{\widehat{\sfg}}$ \cite{BZ05,HL15,KKKO18,KKOP20b}.
Motivated from this,  we conjecture for arbitrary $\g$ that the pair
$(\Lambda^{\tw},B^{\tw})$ satisfying ~\eqref{cond:din},   explicitly can be
calculated by using $\weyl_\g$ and
$\sfQ_\g$ $(\g$ need not to be the same as $\sfg)$, is compatible  (Conjecture~\ref{conj: compatible}). Under the $Q$-adapted condition in~\eqref{cond:din add}, we prove the conjecture:

\begin{thmx} [Theorem~\ref{thm: compatible conj}] \label{main thm: 1.4}
Let $\tw=\sseq{ i_k}_{1 \le k \le r}$ $(r \in \Z_{\ge1} \sqcup \{ \infty \})$ be a $Q$-adapted sequence of $I_\g$  for some Dynkin quiver $Q$
of type $\g$. Then the pair $(\La^\tw,B^\tw)$  is compatible; i.e., for $1 \le k,l \le r$, we have
$$  (\La^\tw B^\tw)_{k,l} =2   \delta(k=l)  d_k  \qquad \text{for some $d_k \in \Z_{>0}$.} $$

\end{thmx}

 We expect that  the quantum seed  associated to the compatible pair give a quantum
cluster algebra structure on certain
subalgebra of the quantum subtorus $\calX_{\ttq,\tw}$ inside
$\calX_\ttq$.

\subsection{Future works}
As we observe in this paper, $\tuC(t)$ is closely related to
the quantum coordinate ring
$A_\nu(\n)$ for any $\g$. Note that $A_\nu(\n)$ is
also isomorphic to the Grothendieck ring $K(\Rep(R_\g))$ of the
category $\Rep(R_\g)$ of finite  dimensional $\Z$-graded modules over the quiver
Hecke algebra $R_\g$ of type $\g$ (\cite{KL09,KL11} and \cite{R08}). In \cite{KO22}, we will prove
that the Laurent  expansions of $\tuC(t)_{i,j}$'s at $t=0$
\emph{determine} whether a given \emph{convolution product} of
\emph{$[Q]$-cuspidal modules} over  $R_\g$ is simple or not. It is known
for Dynkin quivers $Q$ of simply-laced type but not known for the non simply-laced type. This
phenomenon tells that $\tuC(t)$ controls the representation theory
of $R_\g$.

As far as the authors know,
there is \emph{no} Hopf algebra $\scrA$ whose Grothendieck ring is recovered by the
limit  $q \to 1$ of $\calW_{q,t}(\g)$ $(\g \ne \sfg)$ in the sense of~\eqref{eq: recovering rep theory}:
\vspace{-.5em}
\begin{align} \label{eq: recovering rep theory II}
 \raisebox{3.4em}{ \xymatrix@R=3ex@C=10ex{   &   K\left(\scrC_\scrA \right) \ar@{<-->}[ddr] \\
& \mathcal{W}_{q,t}(\g) \ar@{->}[dl]^{ 1 \gets t } \ar@{->}[u]_{ q \to 1  } \ar@{->}[dr]_{ q \to {\rm exp}(\pi i /r) \ \   }  \\
K\left(\scrC_{\widehat{\g}} \right)  \ar@{<-->}[uur]  \ar@{<-->}[rr] & & K\left(\scrC_{{}^L\widehat{\g}}\right)
}}
\end{align}

\vspace{-.5em}

The ultimate goal of this project is to construct new algebras $\scrA =\mathcal{U}'_\nu(\widehat{\g})$
whose representation theory is \emph{controlled} by $\tuC(t)$ and recovered by the limit $q \to 1$ of $\calW_{q,t}(\g)$ for non simply-laced type $\g$. We expect that
such an algebra might be constructed as a certain subquotient of the untwisted quantum affine algebra $U_\nu'(\widehat{\mathbf{g}})$ of simply-laced type $\mathbf{g}$ (different from $\sfg$ in general),
since $\Gamma_Q$ of type $\g$ Dynkin quiver $Q$ can be obtained from $\Gamma_{\overline{Q}}$ for the corresponding Dynkin quiver $\overline{Q}$ of type $\mathbf{g}$
via certain \emph{folding} (\cite{KO22}).
Here the correspondence $\g$ and $\mathbf{g}$ are given as below:
$$ (\g,\mathbf{g}) =  (B_n,D_{n+1}), \ \  (C_n,A_{2n-1}), \ \ (F_{4},E_{6}),  \ \ (G_2,D_4).$$

\medskip

\noindent \textbf{Remark.} When the authors almost finish the first
draft of this paper, Frenkel-Hernandez-Reshetikhin uploaded the
paper \cite{FHR21} at Arxiv, which looks deeply related to this
present paper.
They also consider the refined ring of
interpolating $(q,t)$-character $\overline{\calK}_{q,t}(\g)$ whose
specialization at $q=1$ and $\al = d$ ($d$ is  the lacing number of
$\g$) coincides with the specialization of $\frakK_\ttq$ at $\ttq
=1$ by observing  the ``folded'' phenomenon of $t$-characters  (see Remark~\ref{rmk: motivation} for more discussion).  They also raised the similar question related to \eqref{eq: recovering rep theory II}
in \cite[Remark 3.2]{FHR21}.

\smallskip

\noindent \textbf{Acknowledgments}\quad   The second author
   is grateful to
I.-S.\ Jang, Y.-H.\ Kim, K.-H.\ Lee and R.\ Fujita for helpful
discussions.

\begin{convention}  \hfill
\bnum
\item For a statement $\mathtt{P}$, $\delta(\mathtt{P})$ is $1$ or $0$ according to whether
$\mathtt{P}$ is true or not.  
\item
For a finite set $A$, we denote by $|A|$ the number of elements of $A$.
\item
We denote by $\lec$ the partial order on $\Z$ defined by:
$$m\lec n\Longleftrightarrow \text{$m\le n$ and $m\equiv n\bmod2$.}$$
\ee
\end{convention}

\section{Cartan datum and Dynkin quiver} \label{sec: Cartan Q datum}
In this section, we first fix the notation on Cartan data and Dynkin quivers. Then we investigate commutation classes of the longest element $w_0$ of Weyl group associated to Dynkin quivers   of an arbitrary finite type.

\subsection{Cartan datum}\label{sec:Cartan}
Let $I$ be an index set. A \emph{Cartan datum} is a quintuple $$(\cm,\wl,\sr,\cwl,\scr)$$ consisting of

\ben
\item a generalized symmetrizable Cartan matrix $\cm=(\sfc_{i,j})_{i,j\in I}$,
\item a free abelian group $\wl$, called the \emph{weight lattice},
\item $\sr=\{\al_i \mid  i \in I \} \subset \wl$, called the set of \emph{simple roots},
\item $\cwl \seteq \Hom_\Z(\wl,\Z)$, called the \emph{coweight lattice} and
\item $\scr=\{ h_i \in \cwl \mid  i \in I \}$, called the set of \emph{simple coroots},
\ee
 satisfying
\eq
&&\left\{\parbox{75ex}{\bna
\item $\lan h_i,\al_j\ran = \sfc_{i,j}$ for $i,j\in I$,
\item $\Pi$ is linearly independent over $\Q$,
\item for each $i \in I$, there exists $\varpi_i \in \wl$, called the \emph{fundamental weight}, such that  $\lan h_j,\varpi_i \ran=  \delta(i=j)$
for all $j \in I$ and
\item there exists a $\Q$-valued symmetric bilinear form $( \cdot  ,  \cdot  )$ on $\wl$ such that $\lan h_i,\la\ran = \dfrac{2(\al_i,\la)}{(\al_i,\al_i)}$
 and $(\al_i,\al_i)\in\Q_{>0}$ for any $i$ and $\la\in \wl$.  \label{it:4}
\ee}\label{cond:cartan}
\right.
\eneq
We set  $\h \seteq \Q \otimes_\Z \cwl$,  $\rl \seteq \bigoplus_{ i \in I} \Z\al_i$ and $\rl^+ \seteq \sum_{i \in I} \Z_{\ge 0} \al_i$,
and define $\wl^+ \seteq \{ \la \in P \mid \lan h_i,\la \ran \ge 0 \text{ for any } i \in I\}$, called the set of \emph{integral dominant weights}.
We denote by $\Phi$ the set of \emph{roots}, by $\Phi^+$ the set of \emph{positive roots} and by $\Phi^-$ the set of \emph{negative roots}.

\subsection{Finite Cartan datum} \label{subsec: FCD}

 For a Cartan datum
$(\cm,\wl,\sr,\cwl,\scr)$, we denote by $\Dynkin$
the corresponding Dynkin diagram. It is a graph with $I$ as the set of vertices and
the set of edges between $i, j\in I$ such that $c_{i,j}<0$.
We denote by $\Dynkin_0$ the set of vertices and $\Dynkin_1$ the set of edges.

For each finite Cartan datum, we take  a symmetric bilinear form $(  \cdot ,  \cdot )$ on $\h^*$
as in \eqref{cond:cartan} \eqref{it:4}.
We set $d_i \seteq (\al_i,\al_i)/2 = (\al_i,\varpi_i)$ for any $i\in I$.

For $i,j \in  I$, $d(i,j)$ denotes the number of edges between $i$ and $j$ in $\Dynkin$.
We have
\eqn
\ang{h_i,\al_j}&&=\bc
2&\text{if $i=j$,}\\
-\max\bl d_j/d_i,\;1\br&\text{if $d(i,j)=1$,}\\
0&\text{if $d(i,j)>1$,}
\ec \allowdisplaybreaks \\[2ex]
(\al_i,\al_j)&&=
\bc
2d_i&\text{if $i=j$,}\\
-\max(d_i,d_j)&\text{if $d(i,j)=1$,}\\
0&\text{if $d(i,j)>1$.}
\ec
\eneqn

In this paper, we   choose
the bilinear form $(  \cdot ,  \cdot )$ such that $(\al,\al)$=2 for short roots $\al$ in $\Phi^+$:
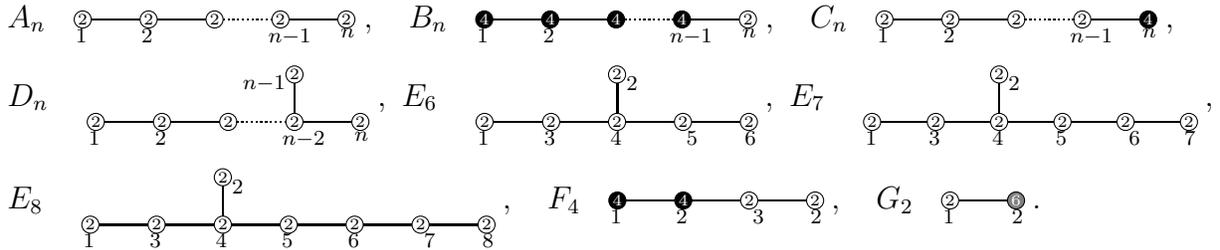
\begin{align*}
&A_n  \ \   \xymatrix@R=0.5ex@C=3.5ex{ *{\circled{2}}<3pt> \ar@{-}[r]_<{ 1 \ \  } & *{\circled{2}}<3pt> \ar@{-}[r]_<{ 2 \ \  } & *{\circled{2}}<3pt> \ar@{.}[r]
&*{\circled{2}}<3pt> \ar@{-}[r]_>{ \ \ n} &*{\circled{2}}<3pt> \ar@{}[l]^>{   n-1} }, \quad
 B_n  \ \   \xymatrix@R=0.5ex@C=3.5ex{ *{\bcircled{4}}<3pt> \ar@{-}[r]_<{ 1 \ \  } & *{\bcircled{4}}<3pt> \ar@{-}[r]_<{ 2 \ \  } & *{\bcircled{4}}<3pt> \ar@{.}[r]
&*{\bcircled{4}}<3pt> \ar@{-}[r]_>{ \ \ n} &*{\circled{2}}<3pt> \ar@{}[l]^>{   n-1} }, \quad
C_n  \ \   \xymatrix@R=0.5ex@C=3.5ex{ *{\circled{2}}<3pt> \ar@{-}[r]_<{ 1 \ \  } & *{\circled{2}}<3pt> \ar@{-}[r]_<{ 2 \ \  } & *{\circled{2}}<3pt> \ar@{.}[r]
&*{\circled{2}}<3pt> \ar@{-}[r]_>{ \ \ n} &*{\bcircled{4}}<3pt> \ar@{}[l]^>{   n-1} }, \  \allowdisplaybreaks\\
&D_n  \ \  \raisebox{1em}{ \xymatrix@R=2ex@C=3.5ex{  &&& *{\circled{2}}<3pt> \ar@{-}[d]_<{ n-1  }  \\
*{\circled{2}}<3pt> \ar@{-}[r]_<{ 1 \ \  } & *{\circled{2}}<3pt> \ar@{-}[r]_<{ 2 \ \  } & *{\circled{2}}<3pt> \ar@{.}[r]
&*{\circled{2}}<3pt> \ar@{-}[r]_>{ \ \ n} &*{\circled{2}}<3pt> \ar@{}[l]^>{   n-2} }},  \
E_{6}  \ \  \raisebox{1em}{   \xymatrix@R=2ex@C=3.5ex{  &&  *{\circled{2}}<3pt> \ar@{-}[d]^<{ 2\ \ }   \\
*{\circled{2}}<3pt>  \ar@{-}[r]_<{ 1 \ \  } & *{\circled{2}}<3pt> \ar@{-}[r]_<{ 3 \ \  } & *{\circled{2}}<3pt> \ar@{-}[r]_<{ 4 \ \  }
&*{\circled{2}}<3pt> \ar@{-}[r]_>{ \ \ 6} &*{\circled{2}}<3pt> \ar@{}[l]^>{   5} } },  \
E_{7}  \ \  \raisebox{1em}{    \xymatrix@R=2ex@C=3.3ex{  &&  *{\circled{2}}<3pt> \ar@{-}[d]^<{ 2\ \ }   \\
 *{\circled{2}}<3pt>  \ar@{-}[r]_<{ 1 \ \  } & *{\circled{2}}<3pt>  \ar@{-}[r]_<{ 3 \ \  } & *{\circled{2}}<3pt> \ar@{-}[r]_<{ 4 \ \  } & *{\circled{2}}<3pt> \ar@{-}[r]_<{ 5 \ \  }
&*{\circled{2}}<3pt> \ar@{-}[r]_>{ \ \ 7} &*{\circled{2}}<3pt> \ar@{}[l]^>{   6} } },  \allowdisplaybreaks \\
&  E_{8}  \ \  \raisebox{1em}{    \xymatrix@R=2ex@C=3.5ex{ &&  *{\circled{2}}<3pt> \ar@{-}[d]^<{ 2\ \ }   \\
 *{\circled{2}}<3pt>  \ar@{-}[r]_<{ 1 \ \  } & *{\circled{2}}<3pt>  \ar@{-}[r]_<{ 3 \ \  } & *{\circled{2}}<3pt> \ar@{-}[r]_<{ 4 \ \  } & *{\circled{2}}<3pt> \ar@{-}[r]_<{ 5 \ \  }   & *{\circled{2}}<3pt> \ar@{-}[r]_<{ 6 \ \  }
&*{\circled{2}}<3pt> \ar@{-}[r]_>{ \ \ 8} &*{\circled{2}}<3pt> \ar@{}[l]^>{   7} } }, \quad
F_{4}   \ \   \xymatrix@R=0.5ex@C=3.5ex{    *{\bcircled{4}}<3pt> \ar@{-}[r]_<{ 1 \ \  } & *{\bcircled{4}}<3pt> \ar@{-}[r]_<{ 2 \ \  }
&*{\circled{2}}<3pt> \ar@{-}[r]_>{ \ \ 2} &*{\circled{2}}<3pt> \ar@{}[l]^>{   3} }, \quad
 G_2  \ \   \xymatrix@R=0.5ex@C=3.5ex{  *{\circled{2}}<3pt> \ar@{-}[r]_<{ 1 \ \  }  & *{\rcircled{6}}<3pt> \ar@{-}[l]^<{ \ \ 2  } }.
\end{align*}
Here $ \xymatrix@R=0.5ex@C=3.5ex{    *{\circled{t}}<3pt>}_k$ means that $(\al_k,\al_k)=t$.

Note that we have $d_i\in\Z_{>0}$ for any $i\in I$ with our choice of the inner product.

\begin{remark}
In this paper, we use   an uncommon    convention for Dynkin diagram; i.e., we do \emph{not} use
doubly-laced (triply-laced) arrows but use circles with integers  instead of vertices which recover the
arrows in the usual convention.  We use this convention to describe the Dynkin quiver
for non simply-laced types.
\end{remark}

Note that the diagonal matrix $\sfD = {\rm diag}(d_i  \mid i \in I )$
symmetrizes $\cm$:
\begin{align*}
\text{The matrices $\osfB\seteq\sfD\cm=\bl(\al_i,\al_j)\br_{i,j\in I}$ and   $\usfB\seteq\cm\sfD^{-1}=\bl(\al^\vee_i,\al^\vee_j)\br_{i,j\in I}$ are symmetric, }
\end{align*}
where $\al^\vee_i=(d_i)^{-1}\al_i$. Note that  the entries in $\osfB$ are integers, while some entries in $\usfB$ are
not (see Example~\ref{Ex: usfB} below).

\begin{example} \label{Ex: usfB}
Note that, for the Cartan matrix $\cm$ of finite type ADE, $\usfB=\cm=\osfB$.
\bna
\item For the Cartan matrix $\cm$ of $B_n$ and $C_n$, $\usfB$ are
\begin{align*}
& \usfB_{B_n}=  \left(\begin{matrix}
1 & -\frac{1}{2} & 0& 0 &\cdots & 0 \\
-\frac{1}{2} & 1  & -\frac{1}{2} &   0 & \cdots & 0 \\
\vdots & \vdots & \ddots  &   \ddots & \cdots & 0  \\
0 & \cdots & \cdots  &   -\frac{1}{2} & 1 & -1  \\
0 & \cdots & \cdots &  0 & -1 & 2
\end{matrix}\right)  \quad \text{and}  \quad
 \usfB_{C_n}=\left(\begin{matrix}
2 & -1 & 0& 0 &\cdots & 0 \\
-1 & 2  & -1 &   0 & \cdots & 0 \\
\vdots & \vdots & \ddots  &   \ddots & \cdots & 0  \\
0 & \cdots & \cdots  &   -1 & 2 & -1  \\
0 & \cdots & \cdots &  0 & -1 & 1
\end{matrix}\right).
\end{align*}
\item For the Cartan matrix $\cm$ of $F_4$ and $G_2$, $\usfB$ are
\begin{align*}
& \usfB_{F_4} =
\left(\begin{array}{rrrr}
1 & -\frac{1}{2} & 0 & 0 \\
-\frac{1}{2} & 1 & -1 & 0 \\
0 & -1 & 2 & -1 \\
0 & 0 & -1 & 2
\end{array}\right)   \quad \text{and} \quad
\usfB_{G_2}=\left(\begin{array}{rr}
2 & -1 \\
-1 & \frac{2}{3}
\end{array}\right).
\end{align*}
\ee
\end{example}

\subsection{Weyl group and convex order}  \label{subsec: Weyl convex}
We denote by $\weyl$ the Weyl group associated to the finite Cartan datum.
It is   the subgroup of $\Aut(\wl)$   generated by simple reflections
$\{ s_i  \mid i \in I \}$:
$$   s_i \la = \la - \lan  h_i,\la \ran \al_i  \qquad (\la \in \wl).$$
Note that
\begin{eqnarray*}&&
\parbox{83ex}{
\bnum
\item there exists a unique element $w_0 \in \weyl$
with the biggest length,
\item $w_0$ induces the Dynkin diagram automorphism
$\;^*\cl I \to I$ sending $i \mapsto i^*$ , where $w_0(\al_i)=-\al_{i^*}$.
\ee}\label{eq: auto *}
\end{eqnarray*}

Let $\uw \seteq s_{i_1} \cdots s_{i_l}$ be a reduced expression of $w \in \weyl$,
 and define
\begin{align}\label{eq: beta_k uw_0}
\be^\uw_k \seteq s_{i_1} \cdots s_{i_{k-1}}(\al_{i_k})   \quad \text{ for } k=1,\ldots,l.
\end{align}
Then we have $\Phi^+ \cap w\, \Phi^-=\{ \be^\uw_1,\ldots,\be^\uw_l \} $ for any reduced expression of $w$  and
$|\Phi^+ \cap w\, \Phi^-|$ coincides with the length $l$ of $w$, denoted by $\ell(w)$.
In particular,  $\Phi^+ \cap w_0 \Phi^- = \Phi^+$ and $|\Phi^+|=\ell(w_0)$. It is well-known that the total order $<_{\uw}$  on $\Phi^+ \cap w\, \Phi^-$ defined by
$\be^\uw_a <_{\uw} \be^\uw_b $ for $a <b$
is \emph{convex} in the following sense:
 \eq \label{convex}
\text{if $\al,\be,\al+\be \in \Phi^+ \cap w\, \Phi^-$,}&&
\text{ we have either}\\
&& \al <_{\uw} \al+\be <_{\uw} \be \quad \text{ or } \quad   \be <_{\uw} \al+\be <_{\uw} \al.\nonumber
\eneq

Two reduced expressions $\uw$ and $\uw'$ of $w \in \weyl$ are said to be
\emph{commutation equivalent}, denoted by $\uw \sim \uw'$,  if
$\uw'$ can be obtained from $\uw$ by applying  the commutation relations $s_as_b=s_bs_a$ $(d(a,b)>1)$. Note that this relation $\sim$  is
 an equivalence relation and
an equivalence class  under $\sim$ is called a \emph{commutation class}.
we denote by $[\uw]$ the  commutation  class of $\uw$.

Now let us consider a reduced expression $\uw_0$ of the longest element $w_0 \in \weyl$.  Then,
each $\uw_0$ induces a convex total order $\le_{\uw_0}$ on $\Phi^+$. For a commutation class $[\uw_0]$, we define the \emph{convex partial order}
 $\preceq_{[\uw_0]}$  on $\Phi^+$  by:
\begin{align}\label{def: convex order}
\al \preccurlyeq_{[\uw_0]} \be  \quad \text{ if and only if }  \quad \al \le_{\uw_0'} \be  \quad \text{ for any } \uw_0' \in [\uw_0].
\end{align}
Note that \eqref{convex} still holds after replacing $<_{\uw_0}$
with $\prec_{[\uw_0]}$.

For a commutation class $[\uw_0]$ of $\uw_0$ and $\al \in \Phi^+$, we define the \emph{$[\uw_0]$-residue} of $\al$, denoted by $\res^{[\uw_0]}(\al)$, to be $i_k \in I$
if $\be^{\uw_0}_k=\al$ with $\uw_0=s_{i_1} \cdots s_{i_\ell}$. Note that this notion is well-defined; i.e, for any $\uw_0' = s_{j_1} \cdots s_{j_\ell} \in [\uw_0]$
with $\be^{\uw'_0}_t=\al$, we have $j_t=i_k$. Note that
\begin{align} \label{eq: residue}
(\al,\al)=(\al_i,\al_i) \qt{if $i=\res^{[\uw_0]}(\al)$.}
\end{align}

For a reduced expression $\uw_0=s_{i_1}s_{i_2} \cdots s_{i_\ell}$
of $w_0$, it is known that the expression $s_{i_2} \cdots s_{i_\ell}s_{i_1^*}$ is also a reduced expression
of $w_0$. This operation is sometimes referred to as a \emph{combinatorial reflection functor} and we write $r_{i_1}\uw_0 = \uw'_0$. Also it induces the operation on commutation classes of $w_0$ (i.e., $r_{i_1}[\uw_0] = [r_{i_1}\uw_0]$
is well-defined  if there exists a reduced expression $\uw'_0=s_{j_1}s_{j_2} \cdots s_{j_\ell}\in[\uw_0]$ such that $j_1=i_1$).
The relations $[\uw] \overset{r}{\sim} [r_i\uw]$ for $i \in I$ generate an equivalence relation, called the \emph{reflection equivalent relation} $\overset{r}{\sim}$, on the
set of commutation classes of $w_0$. For $\uw_0$ of $w_0$, the set of commutation classes $\lf \uw_0 \rf \seteq \{ [\uw_0']\; \mid \;[\uw_0'] \overset{r}{\sim}  [\uw_0] \}$
is called an \emph{$r$-cluster point}.

\subsection{Dynkin quiver $Q$} \label{subsec: D quiver}

A \emph{Dynkin quiver} $Q$ of $\Dynkin$ is an oriented graph whose
underlying graph is $\Dynkin$.   To each Dynkin quiver $Q$ of $\Dynkin$,
we can
associate a function $\xi\cl \Dynkin_0 \to \Z$, called  a
\emph{height function} of $Q$,  which satisfies the condition:
\begin{align*}
 \xi_i = \xi_j +1  \qquad \text{ if } d(i,j)=1 \text{ and } i \to j \text{ in } Q.
\end{align*}
Note that, since $\Dynkin$ is connected, height
functions of $Q$ differ by integers.
Conversely,  to  a Dynkin diagram $\Dynkin$ and a function $\xi:
\Dynkin_0 \to \Z$ satisfying $|\xi_i-\xi_j|=1$ for $i,j\in I$ with
$d(i,j)=1$, we can associate a Dynkin quiver $Q$ in a canonical way.

In this paper,  we abuse the terminology ``Dynkin quiver''
for a pair
$(\Dynkin,\xi)$ of a Dynkin diagram $\Dynkin$ and a height function
$\xi$ on $\Dynkin$.

For a Dynkin quiver $Q=(\Dynkin,\xi)$, we call $i \in \Dynkin_0$ a \emph{source} of $Q$ if $\xi_i > \xi_j$ for all $j \in \Dynkin_0$ such that $d(i,j)=1$.
We also call $i \in \Dynkin_0$ a \emph{sink} of $Q$ if $\xi_i  < \xi_j$ for all $j \in \Dynkin_0$ such that $d(i,j)=1$.

\begin{example} \label{Ex: BCFG Dynkin quiver} Here are examples of Dynkin quivers $Q^\circ$ of non simply-laced types.
\ben
\item \label{it: B Dynkin} $\xymatrix@R=0.5ex@C=6ex{ *{\bcircled{4}}<3pt> \ar@{->}[r]^<{ _{\underline{n}} \ \  }_<{1} & *{\bcircled{4}}<3pt> \ar@{->}[r]^<{ _{\underline{n-1}} \ \  }_<{2} & *{\bcircled{4}}<3pt> \ar@{.>}[r]
&*{\bcircled{4}}<3pt> \ar@{->}[r]^>{ \ \ _{\underline{1}}}_<{n-1} &*{\circled{2}}<3pt> \ar@{}[l]^<{\ \ n}_>{   _{\underline{2}} \ \ } }$ for $Q^\circ=(\bDynkin_{B_n},\xi)$,
\item \label{it: C Dynkin} $\xymatrix@R=0.5ex@C=6ex{ *{\circled{2}}<3pt> \ar@{->}[r]^<{ _{\underline{n}} \ \  }_<{1} & *{\circled{2}}<3pt> \ar@{->}[r]^<{ _{\underline{n-1}} \ \  }_<{2} & *{\circled{2}}<3pt> \ar@{.>}[r]
&*{\circled{2}}<3pt> \ar@{->}[r]^>{ \ \ _{\underline{1}}}_<{n-1} &*{\bcircled{4}}<3pt> \ar@{}[l]^<{\ \ n}_>{   _{\underline{2}}} }$  for $Q^\circ=(\bDynkin_{C_n}, \xi)$,
\item $ \  \xymatrix@R=0.5ex@C=6ex{    *{\bcircled{4}}<3pt> \ar@{->}[r]^<{ _{\underline{4}} \ \  }_<{1} & *{\bcircled{4}}<3pt> \ar@{->}[r]^<{  _{\underline{3}} \ \  }_<{2}
&*{\circled{2}}<3pt> \ar@{->}[r]^>{ \ \  _{\underline{1}}}_<{3 \ \ } &*{\circled{2}}<3pt> \ar@{}[l]^<{ \ \ 4}_>{    _{\underline{2}} \ \ } } \qquad \ \ \ $  for $Q^\circ=(\bDynkin_{F_4}, \xi)$,
\item $ \   \xymatrix@R=0.5ex@C=6ex{  *{\circled{2}}<3pt> \ar@{->}[r]_<{ 1 \ \  }^<{ _{\underline{2}} \ \  }  & *{\rcircled{6}}<3pt> \ar@{-}[l]^<{ \ \ 2  }_<{  \ \ _{\underline{1}}  } }\qquad\qquad\qquad\qquad\quad \ \ $
for $Q^\circ=(\bDynkin_{G_2},\xi)$.
\ee
Here
\bnum
\item an   underlined   integer $\underline{*}$ is the value $\xi_i$ at each vertex $i \in \Dynkin_0$,
\item an arrow $\xymatrix@R=0.5ex@C=4ex{  *{\circled{a}}<3pt> \ar@{->}[r]_<{ i \ \  }   & *{\circled{b}}<3pt> \ar@{-}[l]^<{ \ \ j  }}$ means that
$\xi_i = \xi_j+1$ and $d(i,j)=1$.
\ee
\end{example}

For a Dynkin quiver $Q=(\Dynkin,\xi)$ and its source $i$, we denote by $s_iQ$ the Dynkin quiver $(\Dynkin,s_i\xi)$ where
$s_i\xi$ is the height function defined as follows:
\begin{align} \label{eq: si height}
(s_i\xi)_j = \xi_j-  \delta(i=j)  \times 2.
\end{align}

Let $Q=(\Dynkin,\xi)$ be a Dynkin quiver and let $\weyl$ be the Weyl group associated to $\Dynkin$.
For a reduced expression $\uw=s_{i_1}\cdots s_{i_l}$ of $w \in \weyl$, $\uw$ is said to be \emph{adapted to} $Q$ (or \emph{$Q$-adapted})  if
$$ \text{ $i_k$ is a source of } s_{i_{k-1}}s_{i_{k-2}}\ldots s_{i_1}Q \text{ for all } 1 \le k \le l.$$

\begin{definition} For a  Weyl group $\weyl$, a \emph{Coxeter element} $\tau$ of $\weyl$ is a product of all simple reflections; i.e., there exists a reduced expression $s_{i_1} \cdots s_{i_{n}}$ of $\tau$
such that $\{ i_1,\ldots,i_{n} \}=\Dynkin_0$. Here $n=|\Dynkin_0|$.
\end{definition}

It is well-known that all of reduced expressions of   every   Coxeter element $\tau$ form a single commutation class and they are adapted to some Dynkin quiver $Q$.
Indeed, for a $Q$-adapted $\tau=s_{i_1} \cdots s_{i_{n}}$,
the height function $\xi$ of $Q$ satisfies $\xi_{i_k}=\xi_{i_l}+1$
for $1\le k<l\le n$ such that $d(i_k,i_l)=1$.
Conversely, for each Dynkin quiver $Q$,
there exists a unique Coxeter element $\tauQ$, all of whose reduced expressions are adapted to $Q$.
Furthermore, when $Q$ is of type $B_n$, $C_n$, $F_4$ or $G_2$,
we have
$$ \tauQ^{|\Phi^+|/|\Dynkin_0|} = w_0. $$

The following theorem is proved in \cite{HL15,FO21}
for a Dynkin quiver $Q=(\Dynkin,\xi)$ when $\Dynkin$ is simply-laced,
and we omit the proof in the general case since
it can be proved by similar arguments.

\begin{theorem} \label{thm: unified 1}
Let $Q=(\Dynkin,\xi)$ be a Dynkin quiver and let
$\sfh$ be the Coxeter number of $\Dynkin$.
\bnum
\item There exists a unique Coxeter element $\tauQ \in \weyl_{\Dynkin} $ such that it
has a reduced expression adapted to $Q$.
\item \label{it: tau re} We have $s_i(\tauQ) s_i = \cox{s_iQ}$ for a source $i$ of $Q$ and the order of $\tauQ$ is $\sfh  = \dfrac{ 2 |\Phi^+_{\Dynkin}|}{ | \Dynkin_0|}$.
\item  Any reduced expression $s_{i_1}\cdots s_{i_{n}}$ of $\tauQ$
is $Q$-adapted and
the height function $\xi'$ of the Dynkin quiver
$s_{i_{n}} \cdots s_{i_1}Q$ is given by
$$ \xi'_i = \xi_i - 2  \qt{for any $i \in \Dynkin_0$.}$$
\item There exists a $Q$-adapted reduced expression   of $w_0$ and  all $Q$-adapted reduced expressions of $w_0$ form a single commutation class, which is denoted by $[Q]$.
\item \label{item:iv}Let $\uw_0 = s_{i_1}s_{i_2} \cdots s_{i_\ell}$ be a $Q$-adapted reduced expression of $w_0$.  Then, the height function  $\xi'$ of  the Dynkin quiver
$s_{i_\ell}\cdots s_{i_1}Q$
is given by
$$ \xi'_{i} =\xi_{i^*} -  \sfh   \quad \text{ for any } i \in I_\Dynkin.$$
Moreover, $s_{i_2} \cdots s_{i_\ell}s_{i_\ell+1}$ is an $(s_{i_1}Q)$-adapted reduced expression of $w_0$ where $i_{\ell+1}={i_1}^*$.
\item \label{eq: rf} Let $Q'=(\Dynkin,\xi')$ be another Dynkin quiver.
Then we have $[Q]\overset{r}{\sim} [Q']$. Moreover if  $[Q] = [Q']$, then there exists $k \in \Z$ such that $\xi_i = \xi_i' +k$ for all $i \in \Dynkin_0$.
\item   For a Dynkin diagram $\Dynkin$,
the set $\{ [Q] \mid Q=(\Dynkin,\xi)  \}$ forms   an   $r$-cluster point $\lf \Dynkin  \rf$ and $\bigl\vert\ms{1mu} \lf \Dynkin  \rf\ms{1mu}\bigr\vert = 2^{ |\Dynkin_0|-1}$.
\ee
\end{theorem}

\section{Quivers} \label{sec: AR-quivers}
In this section, we first recall the notion of a
\emph{$($combinatorial\/$)$ Auslander-Reiten $($AR$)$ quiver} which realizes
the convex partial order $\preceq_{[\uw_0]}$  (see \eqref{convex})
for each commutation class $[\uw_0]$ of
$w_0$. Then we introduce quivers related to $[Q]$ and study their properties.

\subsection{Hasse quiver}
For a reduced expression $\uw_0 =s_{i_1} \cdots s_{i_\ell}$ of $w_0 \in \weyl$, we associates a quiver $\Upsilon_{\uw_0}$ to $\uw_0$ as follows \cite{OS19A}:
\begin{eqnarray} &&
\parbox{83ex}{
\bnum
\item The set of vertices is $\Phi^+ = \{ \be^{\uw_0}_k \mid 1 \le k \le \ell \}$.
\item We assign $(-\lan h_{i_k},\al_{i_l} \ran)$-many arrows from
$\be^{\uw_0}_k$ to $\be^{\uw_0}_l$ if and only if $1 \le l  < k \le \ell$, $d(i_k,i_l)=1$  and there is no index
$j$ such that $l < j < k$ and $i_j \in \{i_k,i_l\}$.
\ee
}\label{eq: com AR}
\end{eqnarray}
Hence the number of arrows from $\beta$ to $\gamma$ in $\Upsilon_{\uw_0}$
is either $0$ or
$\max\Bigl(\dfrac{(\ga,\ga)}{(\be,\be)},\;1\Bigr)$.

We say that a total order  $(\be_1< \cdots < \be_\ell)$ of $\Phi^+$ is a \emph{compatible reading} of $\Upsilon_{\uw_0}$ if
\begin{align}\label{eq: com reading}
\text{we have  $k \le l$
whenever there is a path  from $\beta_l$ to $\beta_k$
in $\Upsilon_{\uw_0}$.}
 \end{align}

\begin{remark}  In \cite{OS19A}, $ -( \al_{i_k},\al_{i_l} )$-many arrows was assigned
so that the assigning rule in~\eqref{eq: com AR} is slightly different from the one in \cite{OS19A}.
\end{remark}

\begin{theorem}[{\cite[Lemma 2.19, Proposition 2.20, Theorem 2.21 and 2.22]{OS19A}}] \label{thm: OS17}
The commutation class $[\uw_0]$ of
a reduced expression $\uw_0$ of $w_0$ satisfies the following properties. \
\bnum
\item  A reduced expression $\uw'_0$ of $w_0$ is commutation equivalent to $\uw_0$
if and only if
$\Upsilon_{\uw_0}=\Upsilon_{\uw'_0}$ as quivers.
 Hence $\Upsilon_{[\uw_0]}$ is well-defined.
\item For $\al,\be \in \Phi^+$ $\al  \preceq_{[\uw_0]}  \be$ if and only if there exists a path from $\be$ to $\al$ in $\Upsilon_{[\uw_0]}$. In other words, the quiver
$\Upsilon_{[\uw_0]}$ is the Hasse quiver of the partial ordering $\preceq_{[\uw_0]}$.
\item \label{it: noncom} For $\al,\be \in \Phi^+$,
if they are not comparable with respect to $\preceq_{[\uw_0]}$, then we have $(\al,\be)=0$.
\item \label{it: comp reading}
If
$\Phi^+=\{\be_1< \cdots < \be_\ell\}$ is  a compatible reading of $\Upsilon_{[\uw_0]}$, then
there is a unique reduced expression $\uw_0'=s_{j_1} \cdots s_{j_\ell}$
in $[\uw_0]$ such that
$\beta_k=\be^{\uw'_0}_k$ for any $k$.
\ee
\end{theorem}

We call $\Upsilon_{[\uw_0]}$ the \emph{combinatorial AR-quiver} of $[\uw_0]$.
By Theorem~\ref{thm: OS17},  $\Upsilon_{[\uw_0]}$ can be understood as   the  Hasse quiver of the ordered set $(\Phi^+,\preceq_{[\uw_0]})$
(if we forget the number of arrows).

\begin{example} \label{ex: B3}
Let us consider the following reduced expression $\uw_0$ of $w_0$ of type $B_3$:
\begin{align} \label{eq: B3 the reduced}
 \uw_0 = s_1s_2s_3s_1s_2s_3s_1s_2s_3.
\end{align}
Then one can easily check that $\uw_0$ is adapted to the following Dynkin quiver $Q$:
\begin{align} \label{eq: B3 Q}
Q: \xymatrix@R=0.5ex@C=4ex{ *{\bcircled{4}}<3pt> \ar@{->}[r]^<{\
\underline{3}}_<{1} &*{\bcircled{4}}<3pt>\ar@{->}[r]^>{\underline{2} \quad \   \qquad }_>{2 \quad \   \qquad }
&*{\circled{2}}<3pt> \ar@{-}^>{\ \ \underline{1}}_>{3} }.
\end{align}
The quiver $\Upsilon_{[\uw_0]}=\Upsilon_{[Q]}$
can be described as follows   (see Proposition~\ref{prop:GaUp} below):
\begin{align}\label{eq: B3 ex}
\Upsilon_{[Q]}= \raisebox{1.95em}{ \scalebox{0.7}{\xymatrix@!C=2ex@R=2ex{
(i\setminus p)  &-3& -2 &-1& 0 & 1& 2 & 3\\
1&&& \srt{1,3}\ar@{->}[dr]&& \srt{2,-3} \ar@{->}[dr] && \srt{1,-2}   \\
2&& \srt{2,3} \ar@{->}[dr]\ar@{->}[ur] && \srt{1,2} \ar@{->}[dr] \ar@{->}[ur]  && \srt{1,-3} \ar@{->}[ur] \\
3& \srt{3}  \ar@{=>}[ur] && \srt{2} \ar@{=>}[ur]  && \srt{1} \ar@{=>}[ur]  }}}.
\end{align}
Here,
we use the realization of root system of $B_n$ in $\R^n$
$$ \al_k =\ve_{k}-\ve_{k+1} \quad \text{ for $k<n$},  \quad  \al_n = \ve_n,
$$
where $\{ \ve_i \}_{1 \le i \le n}$ is an orthogonal basis with $(\ve_i,\ve_i)=2$,
and
\bna
\item $\lan a,\pm b \ran \seteq \ve_a \pm \ve_b $ and  $\lan c \ran \seteq
\ve_c $  for $1 \le a < b \le 3$ and $1 \le c \le 3$,
\item every positive root in the $i$-th layer has  $i$ as
its residue with respect to $[\uw_0]$.
\item every root with residue $i$ has the same squared length as $(\al_i,\al_i)$,
\item the indices $p$ will be defined by using the bijection~\eqref{eq: bijection} below.  
\ee

Using the \emph{compatible}   reading  on  $\Phi^+$
\begin{align}
\ang{1,-2}<\ang{1,-3}<\ang{1}<\ang{2,-3}<\ang{1,2}<\ang{2}<\ang{1,3}<\ang{2,3}<\ang{3},  \label{eq: com order B3}
\end{align}
we obtain  the reduced expression $\uw_0$ of $w_0$ in~\eqref{eq: B3 the reduced}.  
\end{example}

\subsection{Classical quivers} Throughout this subsection, we consider a Dynkin quiver
\eqn
&&\text{ $Q=(\Dynkin,\xi)$ of type $ADE$.}\label{ADE}
\eneqn
 We denote by
$\Rep(Q)$ the category of finite-dimensional modules of $Q$ over
$\C$. In this subsection, we recall a description of the
Auslander-Reiten(AR) quiver of $\Rep(Q)$ and that of the derived
category $\calD_Q \seteq D^b(\Rep(Q))$ in   the  aspect of
combinatorics.

\smallskip

  By definition,   the set of vertices of the AR-quiver of $\Rep(Q)$ (resp.\ $\calD_Q$)
is the set of the   isomorphism classes of indecomposable objects, denoted by ${\rm Ind} \;  \Rep(Q)$ (resp.\  ${\rm Ind} \; \calD_Q$).
Considering the bijection $\al \mapsto M_Q(\al)$ from $\Phi^+$ to ${\rm Ind}\; \Rep(Q)=\{ M_k(\al) \mid \al \in \Phi^+\}$  (resp.\ $(\al,k) \mapsto M_Q(\al)[k]$),
the set of vertices of AR-quiver of $\Rep(Q)$ (resp.\ ${\rm Ind} \;  \calD_Q$) can be labeled by $\Phi^+$ (resp.\ $\hPhi^+ \seteq \Phi^+ \times \Z$).
Here $[k]$ denotes the cohomological degree shift by $k$.

\smallskip

We define  the   \emph{repetition quiver $\hDynkin=(\hDynkin_0,\hDynkin_1)$}
 associated to $Q$   as follows:
\begin{equation}\label{eq: rep quiver}
\begin{aligned}
& \hDynkin_0 \seteq \{ (i,p) \in \Dynkin_0 \times \Z \mid  p -\xi_i \in 2\Z\}, \\
& \hDynkin_1 \seteq \{ (i,p) \to (j,p+1)   \mid (i,p) \in \hDynkin_0, \ d(i,j)=1  \}.
\end{aligned}
\end{equation}
Note that $\hDynkin$ depends only on the parity of the height
function of $Q$.

\begin{example} Here is an example of the repetition quiver $\hDynkin$ of simply-laced type.
When $\Dynkin$ is of type $D_{4}$,
the repetition quiver $\hDynkin$ is depicted as:
$$
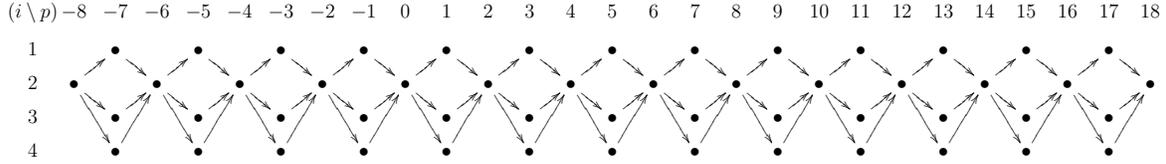
\raisebox{3mm}{
\scalebox{0.63}{\xymatrix@!C=0.3mm@R=2mm{
(i\setminus p) & -8 & -7 & -6 &-5&-4 &-3& -2 &-1& 0 & 1& 2 & 3& 4&  5
& 6 & 7 & 8 & 9 & 10 & 11 & 12 & 13 & 14 & 15 & 16 & 17& 18 \\
1&&\bullet \ar@{->}[dr]&& \bullet \ar@{->}[dr]&&\bullet \ar@{->}[dr]
&& \bullet \ar@{->}[dr]&& \bullet\ar@{->}[dr] && \bullet \ar@{->}[dr]&&\bullet \ar@{->}[dr]&
&\bullet \ar@{->}[dr]&&\bullet\ar@{->}[dr]&& \bullet \ar@{->}[dr]
&&\bullet \ar@{->}[dr]&& \bullet \ar@{->}[dr]&&\bullet \ar@{->}[dr] & \\
2&\bullet \ar@{->}[dr] \ar@{->}[ur]&& \bullet \ar@{->}[dr] \ar@{->}[ur] &&\bullet\ar@{->}[dr] \ar@{->}[ur]
&&\bullet \ar@{->}[dr] \ar@{->}[ur] && \bullet \ar@{->}[dr]\ar@{->}[ur] &&\bullet \ar@{->}[dr] \ar@{->}[ur]&&  \bullet \ar@{->}[dr] \ar@{->}[ur]
&&\bullet \ar@{->}[dr] \ar@{->}[ur] && \bullet \ar@{->}[dr] \ar@{->}[ur]&&\bullet \ar@{->}[dr] \ar@{->}[ur] && \bullet\ar@{->}[dr] \ar@{->}[ur]&&
\bullet\ar@{->}[dr] \ar@{->}[ur] && \bullet\ar@{->}[dr] \ar@{->}[ur]  && \bullet\\
3&& \bullet \ar@{->}[ur]&&\bullet \ar@{->}[ur]&&\bullet \ar@{->}[ur] &&\bullet \ar@{->}[ur]&& \bullet \ar@{->}[ur]
&&\bullet \ar@{->}[ur]&& \bullet \ar@{->}[ur] &&\bullet \ar@{->}[ur]&&\bullet \ar@{->}[ur]&&
\bullet \ar@{->}[ur]&&\bullet \ar@{->}[ur]&&\bullet \ar@{->}[ur]&&\bullet\ar@{->}[ur]&\\
4&& \bullet \ar@{<-}[uul]\ar@{->}[uur]&&\bullet \ar@{<-}[uul]\ar@{->}[uur]&&\bullet \ar@{<-}[uul]\ar@{->}[uur] &&\bullet \ar@{<-}[uul]\ar@{->}[uur]&& \bullet \ar@{<-}[uul]\ar@{->}[uur]
&&\bullet \ar@{<-}[uul]\ar@{->}[uur]&& \bullet \ar@{<-}[uul]\ar@{->}[uur] &&\bullet \ar@{<-}[uul]\ar@{->}[uur]&&\bullet \ar@{<-}[uul]\ar@{->}[uur]&&
\bullet \ar@{<-}[uul]\ar@{->}[uur]&&\bullet \ar@{<-}[uul]\ar@{->}[uur]&&\bullet \ar@{<-}[uul]\ar@{->}[uur]&&\bullet\ar@{<-}[uul]\ar@{->}[uur]&
}}}
$$
\end{example}
It is shown by Happel \cite{Ha87} that the AR-quiver of $\calD_Q$ is isomorphic to $\hDynkin$ as   a quiver  by the following bijection.
We set
\eq &&\ga_i^Q \seteq (1-\tauQ)\varpi_i \in\Phi^+ \qt{for $i \in \Dynkin_0$.}
\eneq
Then the bijection
between sets of vertices is given by
$$ \phiQ \cl \hDynkin_0 \to  {\rm Ind}\;  \calD_Q,  \qquad (i,p) \longmapsto
\uptau^{(\xi_i-p)/2}M_Q(\ga_i^Q),$$
where $\uptau$ denotes the AR-translation of $\calD_Q$.
  By
exchanging the labels of ${\rm Ind}\; \calD_Q$ to $\hPhi^+$, the
bijection $\phiQ\cl \hDynkin_0 \to \hPhi^+$ can be described in an
inductive way as follows (\cite[\S 2.2]{HL15}):
\eq &&
\left\{\parbox{75ex}{
\bnum
\item $\phiQ(i,\xi_i)=(\ga_i^Q,0)$,
\item if $\phiQ(i,p)=(\be,u)$, then we define

\vs{.5ex}
\quad$\phiQ(i,p\pm 2)=\bc
(\tauQ^{\mp1}(\be),u)&\text{if $\tauQ^{\mp1}(\be) \in \Phi^+$,}\\[1ex]
(-\tauQ^{\mp1}(\be),u \pm1)&\text{if $\tauQ^{\mp1}(\be) \in \Phi^-$.}
\ec
$
\ee
}\right.\label{eq: bijection}
\eneq

We say   that $(i,p)$ with $\phiQ(i,p)=(\be,u)$
is  the \emph{coordinate of $(\be,u) \in \hPhi^+$}.

The repetition quiver $\hDynkin$ satisfies the \emph{$\hDynkin$-additive property}
\footnote{In \cite{FO21}, ~\eqref{eq: classic additive} was called $\g$-additive property.
}:
For $i \in I$, $l \in \Z$ and any Dynkin quiver $Q=(\Dynkin,{\rm id},\xi)$, we have
\begin{equation} \label{eq: classic additive}
\tauQ^l(\gamma_i^Q)+\tauQ^{l+1}(\gamma_i^Q)  = \sum_{j; \ d(i,j)=1 } \tauQ^{l+(\xi_j-\xi_i+1)}(\gamma_j^Q) = \sum_{j \in I \setminus \{ i \} }  -\lan h_j, \al_i \ran  \tauQ^{l+(\xi_j-\xi_i+1)}(\gamma_j^Q)
\end{equation}
(see \cite[Proposition 2.1]{HL15}).  

The following statements  are well-known (\cite{B99}):
\eq &\akew[5ex]&
\parbox{80ex}{
\bna
\item The AR-quiver $\Gamma_Q$ of $\Rep(Q)$ is isomorphic to the full subquiver of $\hDynkin$, whose set of vertices $(\Gamma_Q)_0$ is given as follows:
$$
(\Gamma_Q)_0 = \phiQ^{-1}(\Phi^+ \times \{ 0\}) = \{ (i,p) \in \hDynkin_0 \mid \xi_i \ge p > \xi_{i^*}-\sfh  \}.
$$
\item $\Gamma_Q \simeq \Upsilon_{[Q]}$ as quivers.
\ee
}\label{eq: well-known}
\eneq

 \subsection{Repetition quivers in general case }

In this subsection, we consider a Dynkin quiver
$$\text{$Q=(\Dynkin,\xi)$ of an arbitrary finite type.}$$

In this case,
although there is no representation-theoretic interpretation
as for the simply-laced case,
we can generalize  (i) the coordinate system of $\Upsilon_{[Q]}$, (ii) the repetition quiver $\hDynkin$ containing $\Upsilon_{[Q]}$ and (iii) the additive property of $\hDynkin$.

\smallskip

For a Dynkin quiver $Q=(\Dynkin,\xi)$ with $|\Dynkin_0|=n$, recall that
there exists a unique Coxeter element $\tauQ=s_{i_1}\cdots s_{i_n}$with a  $Q$-adapted reduced expression.

For $1 \le k \le n$, we set
\begin{align}\label{eq: gamma k}
 \ga^Q_{i_k} \seteq (1 -\tauQ)\varpi_{i_{k}}=
(s_{i_1}\cdots s_{i_{k-1}}- s_{i_1}\cdots s_{i_k})\varpi_{i_k}=s_{i_1}\cdots s_{i_{k-1}}(\al_{i_k}).
\end{align}

In particular, we have
\eq\label{eq:gai}
(\ga_i^Q,\ga_i^Q)=(\al_i,\al_i)\qt{for any $i\in \Dynkin_0$.}
\eneq
We define the repetition quiver $\hDynkin=(\hDynkin_0,\hDynkin_1)$
associated to $Q$ similarly to the definition \eqref{eq: rep quiver}
in the simply-laced case:
\begin{align*}
& \hDynkin_0 \seteq \{ (i,p) \in I \times \Z \mid p -\xi_i \in 2\Z\}, \\
& \hDynkin_1 \seteq \{ (i,p)\To[{\;-\lan h_{i},\al_j \ran\;}] (j,p+1)   \mid  (i,p) \in \hDynkin_0, \ d(i,j)=1  \}.
\end{align*}
Here $(i,p) \To[{\;-\lan h_{i},\al_j \ran\;}](j,p+1)$ denotes
$(-\lan h_{i},\al_j \ran)$-many arrows from $(i,p)$ to $(j,p+1)$.

\begin{example}
Here is an example of the repetition quiver $\hDynkin$ of non simply-laced type.
When $\Dynkin$ is of type $B_{3}$,
the repetition quiver $\hDynkin$ is depicted as:
$$
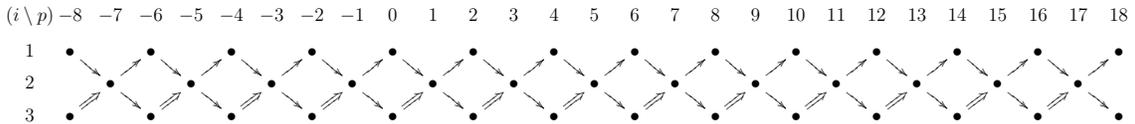
\raisebox{3mm}{
\scalebox{0.6}{\xymatrix@!C=0.5mm@R=2mm{
(i\setminus p) & -8 & -7 & -6 &-5&-4 &-3& -2 &-1& 0 & 1& 2 & 3& 4&  5
& 6 & 7 & 8 & 9 & 10 & 11 & 12 & 13 & 14 & 15 & 16 & 17& 18 \\
1&\bullet \ar@{->}[dr]&& \bullet \ar@{->}[dr] &&\bullet\ar@{->}[dr]
&&\bullet \ar@{->}[dr] && \bullet \ar@{->}[dr] &&\bullet \ar@{->}[dr] &&  \bullet \ar@{->}[dr]
&&\bullet \ar@{->}[dr] && \bullet \ar@{->}[dr] &&\bullet \ar@{->}[dr]  && \bullet\ar@{->}[dr] &&
\bullet\ar@{->}[dr] && \bullet\ar@{->}[dr]  && \bullet\\
2&&\bullet \ar@{->}[dr]\ar@{->}[ur]&& \bullet \ar@{->}[dr]\ar@{->}[ur] &&\bullet \ar@{->}[dr]\ar@{->}[ur]
&& \bullet \ar@{->}[dr]\ar@{->}[ur]&& \bullet\ar@{->}[dr] \ar@{->}[ur]&& \bullet \ar@{->}[dr]\ar@{->}[ur]&&\bullet \ar@{->}[dr]\ar@{->}[ur]&
&\bullet \ar@{->}[dr]\ar@{->}[ur]&&\bullet\ar@{->}[dr] \ar@{->}[ur]&& \bullet \ar@{->}[dr]\ar@{->}[ur]
&&\bullet \ar@{->}[dr]\ar@{->}[ur]&& \bullet \ar@{->}[dr] \ar@{->}[ur]&&\bullet \ar@{->}[dr]\ar@{->}[ur] & \\
3&\bullet  \ar@{=>}[ur]&& \bullet  \ar@{=>}[ur] &&\bullet \ar@{=>}[ur]
&&\bullet  \ar@{=>}[ur] && \bullet \ar@{=>}[ur] &&\bullet  \ar@{=>}[ur]&&  \bullet  \ar@{=>}[ur]
&&\bullet  \ar@{=>}[ur] && \bullet  \ar@{=>}[ur]&&\bullet  \ar@{=>}[ur] && \bullet \ar@{=>}[ur]&&
\bullet \ar@{=>}[ur] && \bullet \ar@{=>}[ur]  && \bullet }}}
$$
\end{example}

\medskip
Let us define the map $\phiQ\cl \Dynkin_0 \to \hPhi^+$ inductively as in~\eqref{eq: bijection}.

By the definition and \eqref{eq:gai}, we have
\eq\label{eq:length}
&&\text{if $\phiQ(i,p)=(\beta,k)$, then
$(\beta,\beta)=(\al_i,\al_i)$.}
\eneq

For $w\in\weyl$, we denote by $\wh{w}$ the
automorphism of $\hPhi^+$ defined as follows (\cite[Lemma B
in \S\,10]{Hum}):
\eq
\wh{w}(\be,k) &&\seteq \bc
(w\beta,k)&\text{if $w\beta\in\Phi^+$,}\\
(-w\beta,k-1)&\text{if $w\beta\in\Phi^-$.}
\ec\label{eq:hatw}
\eneq

Hence we have
\eqn
(\wh{w})^{-1}(\be,k) &&\seteq \bc
(w^{-1}\beta,k)&\text{if $w^{-1}\beta\in\Phi^+$,}\\
(-w^{-1}\beta,k+1)&\text{if $w^{-1}\beta\in\Phi^-$.}
\ec\eneqn

\Lemma\label{lem:hatw}
\smallskip
Let $(\beta,k)\in\hPhi^+$ and
let $s_{i_1}\cdots s_{i_r}$ be a reduced expression of an element $w$ of $\weyl$.
Then we have
$$\wh{w}=\ws_{i_1}\cdots \ws_{i_r}.$$
\enlemma \Proof First let us prove the following statement:
\eq&&\raisebox{1.2em}{
\hs{3ex}\parbox[t]{\mylength}{for a reduced expression
$s_{i_1}\cdots s_{i_r}$ of an element of $\weyl$ and
$\beta\in\Phi^+$, there exists $a$ such that $0\le a\le r$ and
$s_{i_s}\cdots s_{i_r}\beta$ is a positive root or a negative root
according to that $a<s\le r$ or $1\le s\le a$.} }
\label{claim:pn}\eneq

If
$s_{i_k}\cdots s_{i_r}\beta$ is a positive root for all $k$, then it
is enough to take $a=0$. Otherwise, let $a$ be the largest $s$ such
that $s_{i_s}\cdots s_{i_r}\beta$ is a negative root. Then
$s_{i_a}\cdots s_{i_r}\beta=-\al_{i_a}$ and $s_{i_s}\cdots
s_{i_r}\beta=-s_{i_k}\cdots s_{i_{a-1}}\al_{i_a}\in\Phi^-$ for $1\le
s\le a$. Thus we obtain \eqref{claim:pn}.

\smallskip
Then we can easily see that
$$\ws_{i_s}\cdots \ws_{i_r}(\beta,k)=
\bc
(s_{i_s}\cdots s_{i_r}\beta,k)&\text{if $a<s\le r$,}\\
(-s_{i_s}\cdots s_{i_r}\beta,k-1)&\text{if $1\le s\le a$.}
\ec
$$
by the descending induction on $s$.
\QED

By this lemma,
$\st{\ws_i}_{i\in\Dynkin_0}$ satisfies the braid relations. 
Hence it induces the braid group action on $\hPhi^+$.

\smallskip
The following lemma immediately follows from
the above lemma and the definition \eqref{eq: bijection}.

\Lemma\label{eq:phiex} For  any Dynkin quiver $Q$, $(i,p)\in\hDynkin_0$ and $t\in\Z$, we have
\eqn
&&\phiQ(i,p+2t)=\wcox{}^{-t}\phiQ(i,p),\\
&&\phiQ(i,p)=\wcox{}^{(\xi_i-p)/2}(\gamma_i^Q,0).
\eneqn
\enlemma

\Lemma \label{lem:refl}
Let $Q$ be a Dynkin quiver and $(i,p) \in \hDynkin_0$. If $\phiQ(i,p)=(\be,k)$, then  we have $(i^*,p\pm\sfh)\in\hDynkin_0$  and
$$  \phiQ(i^*,p\pm \sfh) = (\be,k   \pm   1).$$
 In particular, $\sfh +\xi_{i^*}-\xi_i\in2\Z$.
\enlemma

\begin{proof}
The simply-laced case is already known (\cite[Corollary 3.40]{FO21}). Hence, we assume that $\Dynkin$ is of non simply-laced type.

Since $(\tauQ)^{\sfh/2}=w_0=-\id$ and
$\ell(\tauQ)\cdot\frac{\sfh}{2}=\ell(w_0)$,
we have
$\wcox{}^{\sfh/2}=\wh{w_0}$.
Hence
\[ \phiQ(i,p\pm \sfh)=\wh{w_0}^{\mp1}(\beta,k)=(\beta,k\pm1). \qedhere \]
\end{proof}

\Prop\label{prop: i-source Q}
Let $i \in \Dynkin_0$ be a source of $Q$.
\bnum
\item\label{sit:1}
$\cox{s_iQ}=s_i\tauQ s_i$ and $\wcox[s_iQ]=\ws_i^{-1}\wcox\ws_i$.
\item \label{sit:2}$\gamma_i^Q=\al_i$, and
\eqn
\gamma_j^{s_iQ}=\bc
s_i\tauQ\gamma_i^Q&\text{if $j=i$,}\\
s_i\gamma_j^Q &\text{if $j\in \Dynkin_0\setminus\st{i}$.}
\ec
\eneqn
\item \label{sit:3}
For any $(j,p) \in \hDynkin_0$ and $k\in\Z$, we have
\eqn
\cox{s_iQ}{}^{((s_i\xi)_j-p)/2}(\ga_j^{s_iQ})&&
=s_i \tauQ{}^{(\xi_j-p)/2}(\ga_j^Q),\\
\wcox[s_iQ]{}^{((s_i\xi)_j-p)/2}(\ga_j^{s_iQ},k)&&=\ws_i^{-1} \wcox{}^{(\xi_j-p)/2}(\ga_j^Q,k).
\eneqn
\item\label{sit:4}
$\phib{s_iQ}=\ws_i^{-1}\circ\phib{Q}$.
\ee

\enprop
\Proof
\eqref{sit:1}\
Let us take a reduced expression $s_{i_1}\cdots s_{i_n}$
of $\tauQ$ such that $i_1=i$.
Then we have $\cox{s_iQ}=s_{i_2}\cdots s_{i_n}\cdot s_{i_1}$.
Hence we obtain \eqref{sit:1}.

\snoi
\eqref{sit:2}\
Since a reduced expression of $s_i\tauQ$ does not contain $s_i$, we have $s_i\tauQ\varpi_i=\varpi_i$.
Hence $\gamma_i^Q=(1-\tauQ)\varpi_i=(1-s_i)\varpi_i=\al_i$, and
$$
\gamma_i^{s_iQ}=(1-\cox{s_iQ})\varpi_i
=(s_i\tauQ-s_i\tauQ s_i)\varpi_i=s_i\tauQ(1-s_i)\varpi_i=
s_i\tauQ\gamma_i^Q.$$

Now assume that $j\not=i$.
Then $s_i\varpi_j=\varpi_j$,  and we have
\eqn
\ga_j^{s_iQ}=(1-\cox{s_iQ})\varpi_j=(s_i-s_i\cox{Q})\varpi_j=s_i\ga_j^{Q}.
\eneqn

\snoi
\eqref{sit:3}\
Assume first $j\not=i$.
Then we have
\eqn
\wcox[{s_iQ}]^{((s_i\xi)_j-p)/2}(\ga_j^{s_iQ},k)=\ws_i^{-1}
\wcox^{(\xi_j-p)/2}\ws_i(s_i\ga_j^Q,k)
=\ws_i^{-1}
\wcox^{(\xi_j-p)/2}(\ga_j^Q,k), 
\eneqn
as we desired.

Now assume that $i = j$.
Note that $\wcox=\ws_i\circ\wh{(s_i\tauQ)}$.
Since $(s_i\xi)_j=\xi_j-2$ , we obtain
\eqn
\wcox[s_iQ]{}^{((s_i\xi)_i-p)/2}(\ga_i^{s_iQ},k)
&&=\ws_i^{-1}\wcox{}^{(\xi_i-2-p)/2}\ws_i(s_i\tauQ\ga_i^{Q},k)\\
&&=\ws_i^{-1} \wcox{}^{(\xi_i-p)/2}\wcox^{-1}\ws_i\wh{(s_i\tauQ)}(\ga_i^Q,k)
=\ws_i^{-1} \wcox{}^{(\xi_i-p)/2}(\ga_i^Q,k).
\eneqn
as we desired.

\snoi
\eqref{sit:4} follows immediately from \eqref{sit:3} and Lemma~\ref{eq:phiex}.
\end{proof}

\smallskip
The statements in   the next theorem
are proved in the simply-laced case (see \cite{HL15,FO21} for instance), but they were not known in the case of  type BCFG.

\begin{theorem} \label{thm: bijection new}
Let $Q=(\Dynkin,\xi)$ be a Dynkin quiver.
\bnum
\item \label{it: bij}The map $\phiQ$ is a bijection.
\item \label{it:range}
Let $\Gamma_Q$ be the full subquiver of $\hDynkin$ whose set of vertices is $\phiQ^{-1}( \Phi^+  \times \{0\})$.
Then we have
\begin{align}\label{eq: range}
\phiQ^{-1}(\Phi^+\times \{0\})=\{ (i,p) \in \hDynkin_0 \mid  \xi_i \ge p > \xi_{i^*} - \sfh \}.
\end{align}
In particular, we have
\eqn  \phiQ(i,p)=\bl\tauQ^{(\xi_i-p)/2}(\ga_i^Q),0\br  \quad \text{for any $(i,p) \in (\Gamma_Q)_0$.}
\eneqn
\item \label{it:si}
If $i$ is a source of $Q$, then we have
$$\phib{s_iQ}(i^*,\xi_i-\sfh)=\phiQ(i,\xi_i)=(\al_i,0),$$
and
$$\bl \Gamma_{s_iQ} \br_0=\Bl \bl \Gamma_Q \br_0 \setminus\st{(i,\xi_i)}\Br\cup\st{(i^*,\xi_i-\sfh)}.$$

\ee
\end{theorem}

Note that
\eq\sfh\ge |\Dynkin_0|+1\ge \xi_{i^*}-\xi_i+2.
\label{ineq:cox}
\eneq
Note also that
$\xi_{i^*}\equiv\xi_{i}+\sfh\bmod2$
as seen in Lemma~\ref{lem:refl}.

Hence we have
\eqn
(\Gamma_Q)_0&&=\st{ (i,p) \in \hDynkin_0 \mid \xi_{i^*} - \sfh+2 \lec p\lec \xi_i }\\
&&=\Bigl\{(i,\xi_i-2k)\bigm| i\in\Dynkin_0,\ k\in\Z,\;
0\le k\le \dfrac{\sfh+\xi_i-\xi_{i^*}}{2}-1\Bigr\}.
\eneqn
(For $\lec$, see {\bf Convention} at the end of Introduction.)

We also call the quiver $\Gamma_Q $ in  Theorem~\ref{thm: bijection new}~
  \eqref{it:range}   \emph{the \ro combinatorial\/\rfm\ AR-quiver} of $Q$.

\begin{proof}[Proof of Theorem~\ref{thm: bijection new}]

Since the theorem is already known in the simply-laced case,
we shall prove it in the non simply-laced case. Recall that
$\sfh$ is even, $\tauQ^{\sfh/2}=w_0=-1$ on $\Phi$, $\ell(w_0)=\ell(\tauQ)\cdot(\sfh/2)$ and $i^*=i$ in this case.

\medskip
Let us prove first \eqref{it: bij} and \eqref{it:range}.
Set $A=\st{ (i,p) \in \hDynkin_0 \mid  \xi_i \ge p > \xi_{i} - \sfh}$.

\bna
\item
We start from the proof of
\eq\phiQ(A)\subset\Phi^+\times\st{0}.
\eneq
It means that
\eq \tauQ^k\,(\ga_i^Q)\in \Phi^+
\qt{for any $k$ such that $0\le k<\sfh/2$.} \label{eq:taup}
\eneq

Recall that $\rtl^+\seteq\sum_{i\in \Dynkin_0}\Z_{\ge0}\al_i$.
Since $\ell(\tauQ^{k+1})=\ell(\tauQ^{k})+\ell(\tauQ)$,
we conclude that $\tauQ^k\,(\ga_i^Q)=\tauQ^k\varpi_i-\tauQ^{k+1}\varpi_i$ belongs to $\rtl^+$. Hence it belongs to $\rtl^+\cap\Phi=\Phi^+$.
\item
Next we shall show that
\eqn
\parbox{70ex}{if $i,j\in \Dynkin_0$ and $m\in\Z$
 satisfies $0\le m<\sfh/2$ and $\tauQ^m\ga_i^Q=\ga_j^Q$,
then $i=j$ and $m=0$.}\eneqn
We have $\tauQ^{\sfh/2-m}\ga^Q_j=\tauQ^{\sfh/2}\ga_i^Q=w_0\ga_i^Q
=-\ga_i^Q\in\Phi^-$.
Hence \eqref{eq:taup} implies $\sfh/2-m=\sfh/2$, which implies that $m=0$.
Thus, we have $\ga_i^Q=\ga_j^Q$.
Hence \eqref{eq: gamma k} implies $i=j$.
\item  Now (b) implies that
$\phiQ\vert_A\cl A\to\Phi^+\times \st{0}$  is injective.

Since $|A|=\dfrac{\sfh}{2}\,| \Dynkin_0|=|\Phi^+|$, the restriction
$\phiQ\vert_A\cl A\to\Phi^+\times \st{0}$ is bijective.
Then Lemma~\ref{lem:refl} implies the bijectivity of $\phiQ$.
\ee

It completes the proof of \eqref{it: bij} and \eqref{it:range}.

\mnoi
\eqref{it:si} is a consequence of \eqref{it: bij}, \eqref{it:range}
and Proposition~\ref{prop: i-source Q}.
\QED

For a subset $S$ of $\hDynkin_0$,
we say that a sequence $\sseq{(i_k,p_k)}_{1\le k\le r}$ in $S$
is
a \emph{compatible reading} (of S) if
\begin{align} \label{eq: comreading}
\text{ $p_a>p_b$ if
$1\le a<b\le r$ and $d(i_a,i_b)\le1$ }
\end{align}
(and $S=\st{(i_k,p_k)\mid 1\le k\le r}$).
If $S=(\Ga_Q)_0$,  this condition is equivalent to the condition:
\eqn&&\text{we have  $a<b$ if
there is an arrow from  $(i_b,p_b)$ to $(i_a,p_a)$. }
\eneqn

We say that a sequence $\sseq{i_k}_{1\le k\le r}$ in $\Dynkin_0$ is {\em $Q$-adapted}
if $i_k$ is a source of
$s_{i_{k-1}}\cdots s_{i_1}Q$ for any $k$ such that $1\le k\le r$.

\Lemma\label{lem:admissible}
Let $Q$ be a Dynkin quiver and let $\sseq{i_k}_{1\le k\le r}$
be a $Q$-adapted sequence in $\Dynkin_0$ with $r\in\Z_{\ge1}$. Set
\eq
p_k=\xi_{i_k}-2\times \big\vert\st{s\in\Z\mid 1\le s<k, i_s=i_k}\big\vert.
\label{eq:ip}
\eneq
Then, we have the followings:
\bnum
\item
$\sseq{(i_k,p_k)}_{1\le k\le r}$ is a compatible reading in $\hDynkin_0$.
\item
$\phiQ(i_k,p_k)=\ws_{i_1}\cdots\ws_{i_{k-1}}(\al_{i_k},0)$ for any $k\in[1,r]$.
\item $s_{i_1}\cdots s_{i_k}\varpi_{i_k}=\tauQ^{(\xi_{i_k}-p_k)/2+1}\varpi_{i_k}$
for any $k\in[1,r]$.\label{item:svarp}
\item Set $\xi'=s_{i_r}\cdots s_{i_1}\xi$.
Then we have
\eq
\xi'_j=\xi_j-2\times \big\vert\st{s\in\Z\mid 1\le s\le r, i_s=j}\big\vert
\qt{for any $j\in\Dynkin_0$,}
\label{eq:xi'}\eneq
and
\eq\st{(i_k,p_k)\mid 1\le k\le r}=\st{(i,p)\in\hDynkin_0\mid \xi'_i <p\le \xi_i}.\label{eq:ipxi} \eneq
\ee
\enlemma
\Proof
Let us first prove (i),(ii),(iii).
We shall show that
\begin{equation*}
\addtocounter{equation}{1}
\setcounter{myce}{\value{equation}}
\setcounter{mycs}{\value{section}}
\left\{\parbox{70ex}{
\begin{enumerate}[(a)]
\item $\phiQ(i_k,p_k)=\ws_{i_1}\cdots\ws_{i_{k-1}}(\al_{i_k},0)$,
\item
$p_s>p_k$ for any $s\in\Z$ such that $1\le s<k$ and $d(i_s,i_k)=1$,
\item $s_{i_1}\cdots s_{i_k}\varpi_{i_k}=\tauQ^{(\xi_{i_k}-p_k)/2+1}\varpi_{i_k}$
\ee}
\right.\tag*{(\theequation)${}_{k,Q}$}
\end{equation*}
by an induction on $k$.
If $k=1$, it is evident.

Assume that $1<k\le r$.

We shall show (\themycs.\themyce)${}_{k,Q}$
assuming that (\themycs.\themyce)${}_{k-1,Q'}$ for $Q'=s_{i_1}Q$.
Define a sequence $\sseq{i'_s}_{1\le s\le r-1}$ by
$i'_s=i_{s+1}$.
Then  $\sseq{i'_s}_{1\le s\le r-1}$
is a $Q'$-adapted sequence.
For an integer $a$ such that $1\le a<\ell$, we have
\eqn
p'_a&&\seteq(s_{i_1}\xi)_{i'_{a}}
-2\times \big\vert\st{s\in\Z\mid 1\le s<a,\; i'_s=i'_{a}}\big\vert\\
&&=(s_{i_1}\xi)_{i_{a+1}}
-2\times \big\vert\st{s\in\Z\mid 1\le s<a, i_{s+1}=i_{a+1}}\big\vert\\
&&=\xi_{i_{a+1}}-2\,\delta(i_{a+1}=i_1)
-2\times\big\vert\st{s\in\Z\mid 2\le s<a+1, i_s=i_{a+1}}\big\vert\\
&&=\xi_{i_{a+1}}-2\times\big\vert\st{s\in\Z\mid 1\le s<a+1, i_s=i_{a+1}}\big\vert\\
&&=p_{a+1}.
\eneqn
Then the induction hypothesis (\themycs.\themyce)${}_{k-1,Q'}$ (a) implies that
$$\phib{Q'}(i'_{k-1},p'_{k-1})=\ws_{i'_1}\cdots\ws_{i'_{k-2}}(\al_{i'_{k-1}},0).$$
Hence we have
$$\phib{Q}(i_k,p_{k})=\ws_{i_1}\phib{Q'}(i'_{k-1},p'_{k-1})=
\ws_{i_1}\Bigl(\ws_{i'_1}\cdots\ws_{i'_{k-2}}\bl \al_{i'_{k-1}},0 \br  \Bigr)
=\ws_{i_1}\ws_{i_2}\cdots\ws_{i_{k-1}}(\al_{i_{k}},0)$$
by Proposition~\ref{prop: i-source Q}\;\eqref{sit:4}.
Thus we have obtained (\themycs.\themyce)${}_{k,Q}$(a).

Let us show (\themycs.\themyce)${}_{k,Q}$\,(b).
Since $p_s=p'_{s-1}$ if $1<s<k$,
(\themycs.\themyce)${}_{k-1,Q'}$\,(b)
implies
(\themycs.\themyce)${}_{k,Q}$\,(b) for $s\not=1$.
When $s=1$ and $d(i_s,i_k)=1$, we have
$$p_1=\xi_{i_1}>\xi_{i_k}\ge p_k.$$
Here the first inequality follows from the fact that $i_1$ is a source of $Q$.

\smallskip
Finally let us show (\themycs.\themyce)${}_{k,Q}$\,(c).
By (\themycs.\themyce)${}_{k-1,Q'}$\,(c), we have
$$\tauQp^{z}\varpi_{i'_{k-1}}=s_{i'_1}\cdots_{i'_{k-1}}\varpi_{i'_{k-1}},$$
where $z=\bl(s_{i_1}\xi)_{i'_{k-1}}-p'_{k-1}\br/2+1=\bl \xi_{i_k}-p_k-2\delta(i_1=i_k)\br/2+1$.

Hence we obtain
\eqn
s_{i_1}\cdots_{i_{k}}\varpi_{i_{k}}
&&=s_{i_1}\tauQp^{z}\varpi_{i'_{k-1}}
=\tauQ^z s_{i_1}\varpi_{i_k}.
\eneqn
Here the last equality follows
from Proposition~\ref{prop: i-source Q}\,\eqref{sit:1}.

If $i_1\not=i_k$, then $z=(\xi_{i_k}-p_k)/2+1$ and $s_{i_1}\varpi_{i_k}=\varpi_{i_k}$, which implies (c).
If $i_1=i_k$ then, we have
$$s_{i_1}\cdots_{i_{k}}\varpi_{i_{k}}
=\tauQ^{(\xi_{i_k}-p_k)/2+1}\tauQ^{-1}s_{i_1}\varpi_{i_k}.
$$
Since $s_{i_1}\tauQ$ does not contain $s_{i_1}$, we have
$\tauQ^{-1}s_{i_1}\varpi_{i_k}=\varpi_{i_k}$.
Thus we obtain (c).
This completes the proof of (i), (ii), (iii).

\medskip
Finally let us prove (iv).
We shall argue by induction on $r$.
Set $\xi''=s_{i_{r-1}}\cdots s_{i_1}\xi$.
Then the induction hypothesis implies
$\xi''_j=\xi_j-2\times \big\vert\st{s\in\Z\mid 1\le s<r, i_s=j}\big\vert$.
Hence, if $j\not=i_r$, then $\xi'_j=\xi''_j$, and
\eqref{eq:xi'} holds.
If $j=i_r$, then $\xi'_j=\xi''_j-2$ and
$$\big\vert\st{s\in\Z\mid 1\le s\le r, i_s=j}=\big\vert\big\vert\st{s\in\Z\mid 1\le s<r, i_s=j}\big\vert+1.$$
Thus we obtain  \eqref{eq:xi'}.

\smallskip
The formula \eqref{eq:ipxi} is obvious by the definition of $p_k$.
\QED

The following proposition says that the set of compatible readings of
$(\Ga_Q)_0$ has a one-to-one correspondence with the set of $Q$-adapted reduced expressions of $w_0$ (see \cite{B99} for ADE-cases).

\Prop\label{prop:GaUp} Let $Q$ be a Dynkin quiver and
let $\ell=|\Phi^+|$.
\bnum
\item
Let $\sseq{(i_k,p_k)}_{1\le k\le \ell}$ be a compatible reading of
$(\Gamma_Q)_0$.
Then we have
\ben
\item
$\uw_0\seteq s_{i_1}\cdots s_{i_\ell}$ is a $Q$-adapted reduced expression
of $w_0$,
\item
$\phiQ(i_k,p_k)=(\beta_k^{\uw_0},0)$ for any $k$ such that $1\le k\le \ell$.
\ee
\vs{.5ex}
\item
Conversely, let $\uw_0=s_{i_1}\cdots s_{i_\ell}$ be
a $Q$-adapted reduced expression of $w_0$, and set
$$p_k=\xi_{i_k}-2\times \big\vert\st{s\in\Z\mid 1\le s<k, i_s=i_k}\big\vert.$$
Then,
$\sseq{(i_k,p_k)}_{1\le k\le \ell}$ is a compatible reading of
$(\Gamma_Q)_0$, and
$$\phiQ(i_k,p_k)=(\beta_k^{\uw_0},0).$$
\ee
\enprop
\Proof
Let us first prove (i).

By the definition of a compatible reading,
$p_1=\xi_{i_1}$ and $i_1$ is a source of $Q$, and $\phiQ(i_1,p_1)=(\al_{i_1},0)$.
Set $\beta_k= s_{i_1}\cdots s_{i_{k-1}}\al_{i_k}$.

We shall show that
\begin{equation*}
\addtocounter{equation}{1}
\setcounter{myce}{\value{equation}}
\setcounter{mycs}{\value{section}}
\left\{\parbox{70ex}{
\bna
\item $\beta_k\in\Phi^+$,
\item $\phiQ(i_k,p_k)=(\beta_k,0)$,
\item $i_a$ is a source of $Q_a\seteq s_{i_{a-1}}\cdots s_{i_1}Q$ for any $a\le k$,
\item $\ell(s_{i_1}\cdots s_{i_k})=k$.
\ee}\right.\tag*{(\theequation)${}_{k,Q}$}
\end{equation*}
by induction on $k$.
If $k=1$, it is evident.

Assume that $1<k\le \ell$.

We shall show (\themycs.\themyce)${}_{k,Q}$
assuming that (\themycs.\themyce)${}_{k-1,Q'}$ with $Q'=s_{i_1}Q$.

By Theorem~\ref{thm: bijection new}\;\eqref{it:si}, we have
$$(\Gamma_{Q'})_0=\bigl((\Gamma_{Q})_0\setminus\st{(i_1,p_1)}\bigr)\cup\st{(i_1^*,\xi_{i_1}-\sfh)}.$$
Set
\eqn
(i'_k,p'_k)&&=
\bc(i_{k+1},p_{k+1})&\text{if $1\le k<\ell$,}\\
(i_1^*,\xi_{i_1}-\sfh)&\text{if $k=\ell$.}
\ec
\eneqn
Then we can easily see that
$\sseq{(i'_k,p'_k)}_{1\le k\le \ell}$ is a compatible reading
of $\Ga_{Q'}$.
Hence, by the induction hypothesis (\themycs.\themyce)${}_{k-1,Q'}$,
$\beta'_{k-1}\seteq s_{i'_1}\cdots s_{i'_{k-2}}\al_{i'_{k-1}}$ belongs to $\Phi^+$ and
$\phib{Q'}(i'_{k-1},p'_{k-1})=(\beta'_{k-1},0)$.
Hence
$\phiQ(i_k,p_k)=\ws_{i_1}(\beta'_{k-1},0)$ by Proposition~\ref{prop: i-source Q}.
Since it belongs to $\Phi^+\times\st{0}$,
we have $\phiQ(i_k,p_k)=(s_{i_1}\beta'_{k-1},0)=(\beta_k,0)$ and
$\beta_k\in\Phi^+$.
Hence we obtain $\ell(s_{i_1}\cdots s_{i_k})
=1+\ell(s_{i_1}\cdots s_{i_{k-1}})=1+(k-1)=k$.
Thus, the induction proceeds and (\themycs.\themyce)${}_{k,Q}$
holds for any $k$ such that $1\le k\le \ell$.
In particular, $s_{i_1}\cdots s_{i_\ell}$ has length $\ell$ and
$\uw_0$ is a reduced expression of $w_0$.

\mnoi
(ii)\  By the preceding Lemma~\ref{lem:admissible}, it remains to show that $(i_k,p_k)\in(\Ga_Q)_0$.
Since $\phiQ(i_k,p_k)=\ws_{i_1}\cdots\ws_{i_{k-1}}(\al_{i_k},0)$ by the preceding lemma and
$\beta_k\seteq s_{i_1}\cdots s_{i_{k-1}}(\al_{i_k})\in\Phi^+$,
Lemma~\ref{lem:hatw} implies that  $\phiQ(i_k,p_k)=(\beta_k,0)$. Hence
 $(i_k,p_k)\in(\Ga_Q)_0$.
\QED

The following proposition follows immediately from Proposition~\ref{prop:GaUp}.

\Prop \label{prop: quiver iso}
The map $\phiQ$ induces a quiver isomorphism
$\Gamma_Q\isoto\Upsilon_{[Q]}$.
\enprop

\begin{theorem}[$\hDynkin$-additive property] \label{thm: additive}
Let $Q=(\Dynkin,\xi)$ be a Dynkin quiver. For any $i \in I$ and $l \in \Z$, we have
\begin{align} \label{eq: additive BCFG}
\tauQ^l(\ga^Q_i)+ \tauQ^{l+1}(\ga^Q_i) = \sum_{j:\;d(i,j)=1 }   -\lan h_j,\al_i\ran \tauQ^{l +(\xi_j-\xi_i+1)/2} (\ga_j^Q).
\end{align}
\end{theorem}

\begin{proof} By Proposition~\ref{prop: i-source Q} \eqref{sit:3}, it suffices to prove the assertion when $i$ is a source of $Q$. Then
\LHS of~\eqref{eq: additive BCFG} is computed as
\begin{align*}
\tauQ^l(\ga^Q_i)+ \tauQ^{l+1}(\ga^Q_i) & = \tauQ^l(1-\tauQ)\varpi_i+ \tauQ^{l+1}(1-\tauQ)\varpi_i \\
& = \tauQ^l(1-\tauQ)(1+\tauQ)\varpi_i \eqs \sum_{j; \ d(i,j)=1} -\lan h_j,\al_i\ran \tauQ^l(1-\tauQ)\varpi_j,
\end{align*}
which implies our assertion. Here $\eqs$ holds by the fact that
$$  (1+\tauQ)\varpi_i = 2\varpi_i -\al_i =   \sum_{j; \ d(i,j)=1} -\lan h_j,\al_i\ran\varpi_j$$
under the assumption that $i$ is a source of $Q$.
\end{proof}

\begin{remark}
When $\Dynkin$ is of simply-laced type, the $\hDynkin$-additive property is closely related to the Auslander-Reiten theory for the path algebra $\C Q$, where $Q$ is a Dynkin quiver
on $\Dynkin$ (see~\cite{Ha87,S14} for instances).
For the non simply-laced type $\bDynkin$, the path algebra $\C Q$ has infinitely many indecomposable modules and the Auslander-Reiten theory for $\C Q$
is not well-investigated. Thus the relation between the representation theory of $\C Q$
and $\hbDynkin$ is not clear.
On the other hand, $\hDynkin$-additive property for simply-laced type $\Dynkin$ is also closely related to $T$-system
in the representation theory of quantum affine algebras \cite{Nak03,H06} and unipotent coordinate algebra \cite{FZ99}. Even in the non  simply-laced type $\bDynkin$,
the~\eqref{eq: additive BCFG} is related to the $T$-system among unipotent quantum minors in~\eqref{eq: T-system}
 below by considering weights of the unipotent quantum minors.
\end{remark}

\section{The $(q,t)$-Cartan matrix specialized at $q=1$ and AR-quivers} \label{Sec: specialization}
In this section, we first recall the $(q,t)$-Cartan matrix $\cm(q,t)$  introduced by Frenkel-Reshetikhin in \cite{FR99}. Then we prove the inverse $\tuC(t)$ of the matrix $\usfC(t) \seteq \cm(1,t)$
can be obtained from $\Gamma_Q$ for any Dynkin quiver $Q$ whose type is the same as the one of $\cm(q,t)$. Finally, we give explicit closed formula for $\tuC(t)$
for an arbitrary Dynkin diagram.

\subsection{The $(q,t)$-Cartan matrix} For an indeterminate $x$ and $k \in \Z$, we set,
$$   [k]_x \seteq \dfrac{x^k-x^{-k}}{x-x^{-1}}.$$
For an indeterminate $q$ and $i$, we set $q_i \seteq q^{d_i}$. For instance, when  $\Dynkin$ is of finite type $G_2$, we have $q_2 = q^{3}$.

For a given finite Cartan datum, we set $\calI=(\calI_{i,j})_{1\le i,j \le n}$ the \emph{adjacent matrix} of $\cm$ as follows:
\begin{align*}
 \calI_{i,j} =\; -\delta(i\not=j)\sfc_{i,j}\;=2\delta(i=j) - \sfc_{i,j}\in\Z_{\ge0}
\qt{for } i,j \in I.
\end{align*}

In~\cite{FR99}, the $(q,t)$-deformation of Cartan matrix
$\cm(q,t) = \bl \sfc_{i,j} (q,t)\br_{i,j \in I}  $ is introduced:
$$
\sfc_{i,j}(q,t) \seteq (q_it^{-1}+q_i^{-1}t) \;  \delta(i=j)   -[\calI_{i,j}]_{q}.
$$
The specialization of  $\cm(q,t)$ at $t=1$, denoted by $\cm(q)\seteq \cm(q,1)$, is usually called the \emph{quantum Cartan matrix}.

\subsection{$t$-quantized  Cartan matrix}

In this paper, we study the
specialization of $\cm(q,t)$ at $q=1$ and its symmetrizations by
  $\sfD  \seteq{\rm diag}(d_i \mid i \in I  )$ and $\sfD^{-1}$.

\begin{definition}
For each finite Cartan datum, we set
$$ \usfC(t) \seteq \cm(1,t)$$
and call it the \emph{$t$-quantized Cartan matrix}. We also set
$$ \usfB(t) \seteq \usfC(t) \sfD^{-1} = (\usfB_{i,j}(t))_{i,j\in I} \quad  \text{ and }   \quad   \osfB(t) \seteq \sfD \usfC(t) = (\osfB_{i,j}(t))_{i,j\in I}.$$
\end{definition}

Hence we have
\eq
\bl\usfB(t)\br_{i,j}&&=\bc
d_i^{-1}(t+t^{-1})&\text{if $i=j$,}\\
-\max\bl (d_i)^{-1},(d_j)^{-1}\br&\text{if $d(i,j)=1$,}\\
0&\text{if $d(i,j)>1$.}\ec
\eneq

\begin{example} Note that, for simply-laced types, we have $\usfB(t)= \usfC(t)$. The followings are $\usfB(t)$ for non simply-laced types:

\begin{align*}
& \usfB(t)_{B_n}=\scriptstyle{ \left(\begin{matrix}
\frac{t+t^{-1}}{2} & -\frac{1}{2} & 0& 0 &\cdots & 0 \\
-\frac{1}{2}&\frac{t+t^{-1}}{2} & -\frac{1}{2} &   0 & \cdots & 0 \\
\vdots & \vdots & \ddots  &   \ddots & \cdots & 0  \\
0 & \cdots & \cdots  &   -\frac{1}{2} & \frac{t+t^{-1}}{2}  & -1  \\
0 & \cdots & \cdots &  0 & -1 &  t+t^{-1}
\end{matrix}\right)}, \allowdisplaybreaks \\
& \usfB(t)_{C_n}=\scriptstyle{\left(\begin{matrix}
t+t^{-1} & -1 & 0& 0 &\cdots & 0 \\
-1 & t+t^{-1}  & -1 &   0 & \cdots & 0 \\
\vdots & \vdots & \ddots  &   \ddots & \cdots & 0  \\
0 & \cdots & \cdots  &   -1 & t+t^{-1} & -1  \\
0 & \cdots & \cdots &  0 & -1 & \frac{t+t^{-1}}{2}
\end{matrix}\right) }, \allowdisplaybreaks \\
& \usfB(t)_{F_4} =
\scriptstyle{\left(\begin{array}{rrrr}
 \frac{(t+t^{-1})}{2}& -\frac{1}{2} & 0 & 0 \\
-\frac{1}{2} &  \frac{t+t^{-1}}{2} & -1 & 0 \\
0 & -1 &t+t^{-1} & -1 \\
0 & 0 & -1 & t+t^{-1}
\end{array}\right)}, \allowdisplaybreaks\\
&  \usfB(t)_{G_2}=\scriptstyle{\left(\begin{array}{rr}
t+t^{-1} & -1 \\
-1 & \frac{t+t^{-1}}{3}
\end{array}\right)}.
\end{align*}

\end{example}

Note that $\usfB(t)\vert_{t=1}=\usfB \in {\rm GL}_{I}(\Q)$. We regard $\usfB(t)$ as an element of ${\rm GL}_{I}(\Q(t))$ and denotes its inverse
by $\tuB(t)=(\tuB_{i,j}(t))_{i,j\in I}$. Let
$$ \tuB_{i,j}(t) =\sum_{u\in\Z} \tfb_{i,j}(u)t^u$$
be the Laurent expansion of $\tuB(t)$ at $t=0$.
\smallskip

Since $t\usfB(t)$ has no pole at $t=0$ and
$\bl t\usfB(t)\br\vert_{t=0}=\sfD^{-1}$,   we obtain the following lemma:

\begin{lemma} \label{lem: basic tfb}
For any $i,j\in I$ and $u \in\Z$, we have
\ben
\item  $\tfb_{i,j}(u)=0$ if $u<1$,
\item  $\tfb_{i,j}(1)=  \delta(i=j)  d_i$.
\ee
\end{lemma}

Since $\usfB(t)$ is symmetric, $\tuB(t)$ is also symmetric: $\tfb_{i,j}(u) = \tfb_{j,i}(u)$ for all $i,j \in I$ and $u \in \Z$.

\subsection{Computations via AR-quivers} In this subsection, we show that
$\tuB(t)$ can be computed by reading   the  AR-quiver $\Gamma_Q$ for an arbitrary Dynkin quiver $Q$.  This
result is already proved  for Dynkin quivers
of type $ADE$ in~\cite{HL15} (see also \cite{FM21}).

\smallskip
We fix a finite Dynkin diagram  $\Dynkin$.
Throughout this subsection, $Q$ denotes a Dynkin quiver on $\Dynkin$.

\begin{definition} For a Dynkin quiver $Q=(\Dynkin,\xi)$ and $i,j \in
\Dynkin_0$, we define a function $\eta_{i,j}^Q\cl\Z \to \Z$ by
$$
\eta_{i,j}^Q(u) \seteq \bc
\bl\varpi_i,\tauQ^{(u+\xi_j-\xi_i-1)/2}(\ga_j^Q)\br & \text{ if } u+\xi_j-\xi_i-1\in 2\Z, \\
0& \text{ otherwise}.
\ec
$$
Note that $\xi_i-\xi_j\equiv d(i,j)\bmod 2$.
\end{definition}

\begin{lemma} \label{lem: eta ij well-defined}
Let $Q'$ be another Dynkin quiver   on $\Dynkin$.
Then we have $\eta_{i,j}^Q = \eta_{i,j}^{Q'}$.
\end{lemma}

\begin{proof}
For our assertion, it suffices to prove that
$$   \eta_{i,j}^Q (u)= \eta_{i,j}^{s_kQ}(u) , \quad \text{
if $k \in I$ is a source of $Q$ and $u+\xi_j-\xi_i-1\in 2\Z$.}  $$

When $k \ne i$, we have $s_k \varpi_i =  \varpi_i$ and $(s_k\xi)_i =\xi_i$,
and hence Proposition~\ref{prop: i-source Q} tells that
$$
\eta_{i,j}^Q(u) = (\varpi_i,s_k \tauQ^{(u+\xi_j-\xi_i-1)/2}(\ga_j^Q))
= (\varpi_i, \cox{  s_k  Q}^{(u+(s_k\xi)_j-(s_k\xi)_i-1)/2}(\ga_j^{s_kQ})) =\eta_{i,j}^{s_kQ}(u).
$$
When  $k =i$, we have $s_i \tauQ \varpi_i=\varpi_i$ since $s_i \tauQ$ has a reduced expression without $s_i$. Using Proposition~\ref{prop: i-source Q} once again, we have
$$
\eta_{i,j}^Q(u) = \bl s_k \tauQ \varpi_i,s_k \tauQ^{(u+\xi_j-\xi_i+1)/2}(\ga_j^Q)\br = \bl\varpi_i, \cox{  s_k  Q}^{(u+(s_k\xi)_j-(s_k\xi)_i-1)/2}(\ga_j^{s_kQ})\br =\eta_{i,j}^{s_kQ}(u),$$
which completes our assertion.
Note that $(s_k\xi)_i=\xi_i-2\delta(k=i)$.
\end{proof}

By the above lemma, the following notation is well-defined:
$$ \eta_{i,j} \seteq  \eta_{i,j}^Q,$$
where $Q$ is an arbitrary Dynkin quiver on $\Dynkin$.

\begin{lemma} \label{lem: key}
The function $\eta_{i,j}$ satisfies the following properties: for any $u \in \Z$, we have
\bnum
\item \label{it: first} $\eta_{i,j}(u+\sfh)=-\eta_{i,j^*}(u)$,
\item \label{it: add} $\eta_{i,j}(u-1)+ \eta_{i,j}(u+1) = \displaystyle
\sum_{k\,;\,d(k,j)=1}  -\lan h_k,\al_j \ran \eta_{i,k}(u)$,
\item \label{it:ske} $-\eta_{i,j}(-u)=\eta_{j,i}(u)$,
\item \label{it: vanish}   $\eta_{i,j}(0)=0$,
\item \label{it: d_i}   $\eta_{i,j}(1)=d_i\;   \delta(i=j) $.
\ee
\end{lemma}

\begin{proof}
\eqref{it: first} For a Dynkin quiver $Q$ of simply-laced type, it is proved in \cite[Lemma 3.7]{Fuj19}. For a Dynkin quiver $Q$ of non simply-laced type, it is a consequence of the fact that
$\tauQ^{\sfh/2}=-1$.

\snoi
\eqref{it: add}   follows from the $\hDynkin$-additive property in Theorem~\ref{thm: additive}.

\snoi
\eqref{it:ske} Take a Dynkin quiver $Q$ and $u \in \Z$ such that $u+\xi_i-\xi_j-1 \in 2\Z$. Then we have
\begin{align*}
-\eta_{i,j}(-u) & = (\varpi_i,\tauQ^{(-u+\xi_j-\xi_i-1)/2}(\tauQ-1)\varpi_j) =  (\tauQ^{(u+\xi_i-\xi_j+1)/2}(\tauQ^{-1}-1)\varpi_i,\varpi_j) \\
& =  (\varpi_j,\tauQ^{(u+\xi_i-\xi_j-1)/2}(1-\tauQ)\varpi_i) = \eta_{j,i}(u).
\end{align*}

\snoi
\eqref{it: vanish},\;\eqref{it: d_i}\hs{1.3ex}By Lemma~\ref{lem: eta ij well-defined},
we may assume that $Q$ satisfies that $\xi_j=1$ and
$\xi_k \in \{0,1\}$ for any $k\in\Dynkin_0$. Then for any $k\in I$ with $\xi_k=1$, $k$ is a source and hence $\ga_k^Q=\al_k$ by~\eqref{eq: gamma k}.
For $i ,  j\in I$ with $\xi_i=\xi_j=1$, we have
$\eta_{i,j}(1) = (\varpi_i,\al_j)= \delta(i=j)$.
When $i ,  j\in I$ with $\xi_i=0 \ne \xi_j=1$,
we have  $\eta_{i,j}(0) = (\varpi_i,\al_j)=0$.
\end{proof}

\begin{theorem} \label{thm: inv}
 We have $\tfb_{i,j}(u)=\eta_{i,j}(u)$ for any $i,j \in I$ and $u \in \Z_{\ge 0}$.
In other words, we have
$$
\tfb_{i,j}(u) = \bc \bl\varpi_i,\tauQ^{(u+\xi_j-\xi_i-1)/2}(\ga_j^Q) \br &
\qt{if $u\in\Z_{\ge0}$ and $u+\xi_j-\xi_i-1\in2\Z$,} \\[.5ex]
0 & \qt{otherwise}
\ec
$$
for any Dynkin quiver $Q$, $i,j \in I$ any $u\in\Z$.
In particular, $\tfb_{i,j}(u) \in \Z$.
\end{theorem}

\begin{proof}
Define
$$  H_{i,j}(t) \seteq \sum_{u\ge0} \eta_{i,j}(u)t^u
\eqs\sum_{u\ge1} \eta_{i,j}(u)t^u \in   t  \Z\lf t \rf \quad \text{ for } i,j \in I.$$
Here, $\eqs$ follows from Lemma~\ref{lem: key}~\eqref{it: vanish}.

It is enough to show
\begin{align}\label{eq: to show inv}
\sum_{k \in I } H_{i,k}(t)\usfB_{k,j}(t) =  \delta(i=j) \qt{for any $i,j\in I$.}
\end{align}
Let $x_{i,j}(u)$ be the coefficient of $t^u$ in \LHS of~\eqref{eq: to show inv}.
Since $ \usfB_{k,j}(t)\in t^{-1}\Q[[t]]$, it is enough to show
$x_{i,j}(u)=\delta(u=0)  \delta(i=j)$ for $u\ge0$.
Now we have
$$   d_j x_{i,j}(u)= \delta(u> -1) \eta_{i,j}(u+1)+ \delta(u\ge 1)\,\eta_{i,j}(u-1) +\hs{-1ex}\sum_{k\,; \; d(k,j)=1}\hs{-1ex} \lan h_k,\al_j \ran \delta(u> 0)\eta_{i,k}(u)$$
for any $u\in\Z$.
By Lemma~\ref{lem: key}~\eqref{it: add},  $x_{i,j}(u)$ vanishes for $u \ge 1$.
Finally, Lemma~\ref{lem: key}  \eqref{it: d_i} implies
\[ d_jx_{i,j}(0)= \eta_{i,j}(1)=d_i  \delta(i=j). \qedhere \]
\end{proof}

\begin{corollary} \label{cor: computation setting}
The coefficients $\{ \tfb_{i,j}(u) \mid i,j \in I, \ u \in \Z \}$ enjoy the following properties:
\bnum
\item \label{it: parityvan}
  $\tfb_{i,j}(u)=0$ unless $u\ge1$ and $u\equiv d(i,j)+1\bmod2$.
\item \label{it: h nega} $\tfb_{i,j}(u+\sfh)=-\tfb_{i,j^*}(u)$ for $u \ge 0$.  In particular, $\tfb_{i,j}(\sfh)=0$.  
\item \label{it: 2h pd}  $\tfb_{i,j}(u+2\sfh)=\tfb_{i,j}(u)$ for $u \ge 0$.
\item \label{it: h pd}$\tfb_{i,j}(\sfh-u)=\tfb_{i,j^*}(u)$ for $0 \le u \le \sfh$.
\item \label{it: 4th}  $\tfb_{i,j}(2\sfh-u)=-\tfb_{i,j}(u)$ for $0 \le u \le 2\sfh$.
\item\label{it: tfb positive} $\tfb_{i,j}(u) \ge 0$ for $0 \le u \le \sfh$.
\item \label{it:6th} $\tfb_{i,j}(u)\le 0$ for $\sfh \le u \le 2\sfh$.
\item \label{it:7th} For any Dynkin quiver
$Q$ with a height function $\xi$, $i,j\in I$ and $u \in \Z$, we have
\begin{align}\label{eq: tfb general u}
\tfb_{i,j}(u) -\tfb_{i,j}(-u) = \bc (\varpi_i,\tauQ^{(u+\xi_j-\xi_i-1)/2}(\ga_j^Q)  ) & \text{ if $u+\xi_j-\xi_i-1\in2\Z$}, \\
0 & \text{ otherwise.}
\ec
\end{align}
\ee
\end{corollary}

\begin{proof}
\eqref{it: parityvan} follows from Theorem~\ref{thm: inv} and the fact that
$\xi_j-\xi_i\equiv d(i,j)\bmod2$.

\snoi
(\ref{it: h nega},\;\ref{it: 2h pd}) follow from Lemma~\ref{lem: key}~\eqref{it: first}, Theorem~\ref{thm: inv} and ~\eqref{it: parityvan}.

\snoi
\eqref{it: h pd} \
 As we have already seen in Lemma~\ref{lem:refl},
we have $\xi_{j^*}-\xi_{j}\equiv \sfh\bmod 2$.
Hence we may assume that
$\sfh-u+\xi_j-\xi_i-1\equiv u+ \xi_{j^*}-\xi_i -1\equiv0\bmod 2$.
Then, we have
\begin{align*}
\tfb_{i,j}(\sfh-u) = \eta_{i,j}(\sfh-u)=-\eta_{i,j^*}(-u)=\eta_{j^*,i}(u) = \tfb_{j^*,i}(u) = \tfb_{i,j^*}(u).
\end{align*}

\snoi
\eqref{it: 4th}
 is a consequence of \eqref{it: h nega} and \eqref{it: h pd}.

\snoi
\eqref{it: tfb positive}\
Let us take a Dynkin quiver $Q$ whose height function $\xi$ satisfies
$\xi_i=0$ and $\xi_k \in \{0,1\}$ for all $k \in \Dynkin_0$.
it is enough to show
that $\eta_{i,j}(u)\ge0$ if $0<u<\sfh$ and $u+\xi_j-\xi_i-1\in2\Z$.
Then $\beta\seteq\tauQ^{(u+\xi_j-\xi_i-1)/2}(\ga_j^Q)$ is a positive root by
Theorem~\ref{thm: bijection new}\;\eqref{it:range} and
$$0\le \dfrac{u+\xi_j-\xi_i-1}{2}<\dfrac{\sfh+\xi_j-\xi_{j^*}}{2}$$
Hence
$\eta_{i,j}(u)=(\varpi_i,\be) \ge 0$.

\snoi
\eqref{it:6th} follows from \eqref{it: tfb positive} and \eqref{it: h nega}.

\snoi
\eqref{it:7th} \ For $u \ge 0$,
\RHS of~\eqref{eq: tfb general u} is equal to $\tfb_{i,j}(u)$ and
hence it is nothing but Theorem~\ref{thm: inv}.
As a function   in $u$,
\RHS of~\eqref{eq: tfb general u} is $2\sfh$-periodic since
$\tauQ^\sfh=1$.
Hence it is enough to show that \LHS
of~\eqref{eq: tfb general u} is also $2\sfh$-periodic:
\eq
\tfb_{i,j}(u+2\sfh) -\tfb_{i,j}(-u-2\sfh)= \tfb_{i,j}(u) -\tfb_{i,j}(-u)
\label{eq:per}
\eneq
for all $u\in\Z$.
If $u\ge0$, then we have
$\tfb_{i,j}(u+2\sfh) -\tfb_{i,j}(-u-2\sfh)=\tfb_{i,j}(u)
= \tfb_{i,j}(u) -\tfb_{i,j}(-u) $
by \eqref{it: 2h pd}.
If $u\le-2\sfh$, then
we have
$$\tfb_{i,j}(u+2\sfh) -\tfb_{i,j}(-u-2\sfh)=-\tfb_{i,j}(-u-2\sfh)
=\tfb_{i,j}(-u)= \tfb_{i,j}(u) -\tfb_{i,j}(-u) $$ also by \eqref{it: 2h
pd}. Finally if $-2\sfh<u<0$, then \eqref{eq:per} holds since
$\tfb_{i,j}(u+2\sfh)= -\tfb_{i,j}(-u)$ and $
-\tfb_{i,j}(-u-2\sfh)=\tfb_{i,j}(u)$ by \eqref{it: 4th}.
\end{proof}

\begin{proposition}\label{prop:Xi}
For any $i,j \in I$ and $l \in \Z_{\ge 0}$, we have
\begin{align}   \label{eq: Xi ij}
\tfb_{i,j}(u)=0\qt{for $u\le d(i,j)$.}
\end{align}
\end{proposition}

\begin{proof}
Since $\tfb_{i,j}(u)=0$ for $u \le 0$, it suffices to
prove~\eqref{eq: Xi ij} for $u=d(i,j)-2l-1$ with an integer $l$ such that
 $0 \le 2l < d(i,j)-1$. Let us take a height function $\xi$ such that
$\xi_j-\xi_i=d(i,j)$ and $\xi_j \ge \xi_{j^*}$.
By~\eqref{eq: tfb general u}, we have
\begin{align}\label{eq: Xi step}
\tfb_{i,j}(-u) -\tfb_{i,j}(u) = (\varpi_i,\tauQ^{l}(\ga_j^Q) ).
\end{align}
Since $-u\le0$, we have $\tfb_{i,j}(-u)=0$ by Lemma~\ref{lem: basic tfb}. Hence \LHS of~\eqref{eq: Xi step} is equal to $ -\tfb_{i,j}(u)$, which is non-positive by Corollary~\ref{cor: computation setting}~\eqref{it: tfb positive}
since $u\le d(i,j)<\sfh$ (see \eqref{ineq:cox}).
On the other hand, \RHS of~\eqref{eq: Xi step} is   non-negative  since
$\tauQ^l(\gamma_j^Q) \in \Phi^+$ by
Theorem~\ref{thm: bijection new}\;\eqref{it:range} and
$$0\le l<\dfrac{d(i,j)-1}{2}<\dfrac{\sfh}{2}\le\dfrac{\sfh+\xi_j-\xi_{j^*}}{2}.$$
Thus, we obtain  $\tfb_{i,j}(u)=0$, as desired.
\end{proof}

Together with Corollary~\ref{cor: computation setting},
we have
\eq\text{ $\tfb_{i,j}(u)=0$ unless $d(i,j)+1\lec u$.} \eneq
Here $u\lec v$ means that $u\le v$ and $u\equiv v\bmod 2$.

\begin{corollary} \label{cor: additive application}
For $i,j \in \Dynkin_0$, let us define an  even function $\teta_{i,j} : \Z \to \Z $ as follows:
\begin{align} \label{eq: teta}
 \teta_{i,j}(u) = \tfb_{i,j}(u)+\tfb_{i,j}(-u) \qquad \text{ for } u \in \Z.
 \end{align}
Then we have
$$
\teta_{i,j}(u-1) + \teta_{i,j}(u+1) + \sum_{k:\;d(k,j)=1} \lan h_k,\al_j \ran \teta_{i,k}(u) =
2d_i\;\delta(u=0)\, \delta(i=j).
$$
\end{corollary}

\begin{proof}
It is a direct consequence of Lemma~\ref{lem: key} and Theorem~\ref{thm: inv}.
\end{proof}

\subsection{Explicit computation}
With the help of Theorem~\ref{thm: inv}, we can obtain
$\tfb_{i,j}(u)$ from any AR-quiver $\Gamma_Q$. In this subsection,
we will explicitly compute $\tfb_{i,j}(u)$  for non simply-laced
type.

\begin{remark} By Corollary~\ref{cor: computation setting}\;\eqref{it: h nega}, it is enough to compute
$$   \tde_{i,j}(t) \seteq (1+t^{\sfh})\tuB_{i,j}(z) = \sum_{u=0}^{\sfh-1} \tfb_{i,j}(u)t^u \quad \text{ for each $i,j\in I$}.$$
\end{remark}

\subsubsection{Simply-laced type} For each simply-laced type $\Dynkin$, $(\tde_{i,j}(t))_{i,j \in \Dynkin_0}$ was explicitly calculated as follows:

\begin{theorem}  [\cite{DO94,HL15,KKKII,Fuj19,OS19}] \label{thm: ADE denom}  Note that $\tde_{i,j}(t) = \tde_{j,i}(t)$  for $i,j \in I$ .
\ben
\item For $\Dynkin$ of type $A_{n}$, and $i,j\in I =\{1,\ldots,n\}$, $\tde_{i,j}(t)$ is given as follows:
\begin{align} \label{eq: A formula}
\tde_{i,j}(t)&   = \sum_{s=1}^{\min(i,j,n+1-i,n+1-j)}t^{|i-j|+2s-1}.
\end{align}
\item For $\Dynkin$ of type $D_{n+1}$, and $i,j\in I =\{1,\ldots,n,n+1\}$, $\tde_{i,j}(t)$ is given as follows:
\begin{align}\label{eq: D formula}
\tde_{i,j}(t) &  = \bc
\displaystyle\sum_{s=1}^{\min(i,j)} \bl t^{|i-j|+2s-1}+\delta(\max(i,j)<n)\,t^{2n-i-j+2s-1})
 & \text{ if } \min(i,j)<n,\\[3ex]
\ \displaystyle \sum_{s=1}^{\lfloor (n+ \delta_{i,j}) /2 \rfloor}  t^{4s-1 -2\ms{1mu}\delta(i,j)} & \text{ otherwise.}
\ec
\end{align}
\item   For $\Dynkin$ of type $E_6$  and $i \le j\in I =\{1,\ldots,6\}$, $\tde_{i,j}(t)$ is given as follows:
\begin{align*} 
&\tde_{1,1}(t) =  t + t^{7},  && \tde_{1,2}(t) =  t^{4} + t^{8},   \allowdisplaybreaks\\
& \tde_{1,3}(t)  = t^{2} +t^{6}+t^{8}, && \tde_{1,4}(t)  = t^{3}+t^{5}+t^{7}+t^{9}, \allowdisplaybreaks\\
& \tde_{1,5}(t)  = t^{4}+t^{6}+t^{10}, &&  \tde_{1,6}(t)  = t^{5}+t^{11},  \allowdisplaybreaks\\
&  \tde_{2,2}(t)  = t^{1}+t^{5}+t^{7}+t^{11}, && \tde_{2,3}(t)  = t^{3}+t^{5}+t^{7}+t^{9},  \allowdisplaybreaks\\
& \tde_{2,4}(t)  = t^{2}+t^{4}+2t^{6}+t^{8}+t^{10},  && \tde_{3,3}(t) = t^{1}+t^{3}+t^{5}+2t^{7}+t^{9}, \allowdisplaybreaks\\
&\tde_{3,4}(t) = t^{2}+2t^{4}+2t^{6}+2t^{8}+t^{10},   && \tde_{3,5}(t) = t^{3}+2 t^{5}+ t^{7}+ t^{9}+ t^{11},\allowdisplaybreaks \\
& \tde_{4,4}(t) = t^{1}+ 2t^{3}+ 3t^{5}+ 3t^{7}+ 2t^{9}+ t^{11},  && \tde_{i,j}(t) =
t^\sfh\tde_{i,j^*}(t^{-1}) = \tde_{j,i}(t) =t^\sfh \tde_{j,i^*}(t^{-1}).
\end{align*}
\item For $E_7$ and $E_8$, see {\rm Appendix~\ref{appeA: tde}.}
\ee
\end{theorem}

\subsubsection{$B_n$ and $C_n$ cases} [$B_n$ case]
By using an orthogonal basis $\{ \ve_i \mid 1 \le i \le n\}$ of $\R^n$
with $(\ve_i,\ve_i)=2$,
the simple roots
 and $\Phi^+_{B_n}$ can be realized as follows:
\begin{align*}
 \al_i  &=\ve_i- \ve_{i+1} \quad \text{ for $i<n$}, \quad \quad \al_n =\ve_n, \\
\Phi^+_{B_n} & = \Bigl\{  \ve_i = \sum_{k=i}^n \al_k  \bigm|  1 \le i \le n \Bigr\} \ssqcup \Bigl\{  \ve_i-\ve_{j} = \sum_{k=i}^{j-1} \al_k \bigm| 1 \le i <j \le n
\Bigr\} \allowdisplaybreaks\\
& \quad \ssqcup \Bigl\{  \ve_i+\ve_{j} = \sum_{k=i}^{j-1} \al_k + 2\sum_{s=j}^n \al_s
\bigm | 1 \le i <j \le n \Bigr\}.
\end{align*}

 Recall the notations in Example~\ref{ex: B3}:
 $\ang{a,b}\seteq\ve_a-\ve_b$ and $\ang{c}\seteq\ve_c$ for $1\le a<b\le n$ and
 $1 \le c\le n$.

For the Dynkin quiver $Q^\circ$ in Example~\ref{Ex: BCFG Dynkin quiver}~\eqref{it: B Dynkin},
we have
\begin{align}\label{eq: tau circ}
 \cox{Q^\circ} =s_{i_1} s_{i_2} \cdots s_{i_n} =s_{1}s_2 \cdots s_n
\end{align}
and the lemma below, by direct calculation:
\begin{lemma} \label{lem: the Q B} \hfill
\bnum
\item $\ga^{Q^\circ}_{i_k}=\ga^{{Q^\circ}}_k = \lan 1, -(k+1) \ran $ for $k <n$ and $\ga_n =\lan  1  \ran $.
\item For $k <n$, we have $\cox{Q^\circ}^{s}(\ga^{{Q^\circ}}_k) = \bc \lan 1+s, -(k+1+s) \ran & \text{ for } 0 \le s < n-k,\\ \lan k+s+1-n, 1+s \ran & \text{ for }n-k \le s \le n-1.\ec$
\item  $\cox{Q^\circ}^{s}(\ga^{{Q^\circ}}_n) = \lan s+1 \ran$ for $0 \le s \le n-1$.
\ee
\end{lemma}

[$C_n$ case]  Note that $\weyl$ of type $C_n$ is isomorphic to the Weyl group of type $B_n$.  By using the orthonormal basis $\{ \ep_i \ | \  1 \le i \le n \}$
 of $\R^n$, the simple roots and $\Phi^+_{C_n}$ can be realized as follows:
\begin{align*}
 \al_i &=\ep_i- \ep_{i+1} \quad \text{ for $i<n$}  \quad \text{ and } \quad   \al_n =2\ep_n, \\
\Phi^+_{C_n}  &= \Bigl\{  2\ep_i = \sum_{k=i}^n \al_k
\bigm| 1 \le i \le n \Bigr\} \ssqcup \Bigl\{  \ep_i-\ep_{j} = \sum_{k=i}^{j-1} \al_k \bigm| 1\le i <j \le n \Bigr\}
\\
& \quad\qquad \ssqcup \Bigl\{  \ep_i+\ep_{j} = \sum_{k=i}^{j-1} \al_k + 2\sum_{s=j}^{n-1} \al_s +\al_n \bigm| 1 \le i <j \le n \Bigr\}.
\end{align*}

For the Dynkin quiver $Q^\circ$ in Example~\ref{Ex: BCFG Dynkin quiver}~\eqref{it: B Dynkin},
$\cox{Q^\circ}$ is the same as~\eqref{eq: tau circ} and we have the following in a similar way as $B_n$-case:

\begin{lemma} \label{lem: the Q C} Setting $\lan a, \pm b \ran \seteq \ep_a \pm \ep_b$ for $1 \le a <b \le n$ and $\lan c, c \ran \seteq 2 \ep_c$ for $1 \le c \le n$,
we have
\bnum
\item $\ga^{Q^\circ}_{i_k}=\ga^{{Q^\circ}}_k = \lan 1, -(k+1) \ran $ for $k <n$ and $\ga_n =\lan  1,1  \ran $,
\item $\cox{Q^\circ}^{s}(\ga^{{Q^\circ}}_k) = \bc \lan 1+s, -(k+1+s) \ran & \text{ for } 0 \le s < n-k,\\ \lan k+s+1-n, 1+s \ran & \text{ for }n-k \le s \le n-1.\ec$
\ee
\end{lemma}

For $Q$ of type $B_n$ or $C_n$, recall
$$(d_i)_{i \in I} = \bc
(2,2,\ldots,2,1) & \text{ if $Q$ is of type $B_n$} , \\
(1,1,\ldots,1,2) & \text{ if   $Q$ is of type $C_n$} ,
\ec \quad\text{ and } \quad \sfh=2n. $$

\begin{theorem} \label{thm: BC denom}
For $\bDynkin$ of type $B_n$ or $C_n$, and $i,j\in I =\{1,\ldots,n\}$, the closed formula of $\tuB_{i,j}(t)$ is given as follows:
for any $i,j\in I$ such that $i\le j$, we have $\tuB_{j,i}(t) = \tuB_{i,j}(t)$  and
\eq\label{eq: BC formula}
&&
\tde_{i,j}(t) =  
\bc
\max(d_i,d_j)\displaystyle\sum_{s=1}^{i} t^{n-i-1+2s}\\[-.3ex]
\hs{8ex}=\max(d_i,d_j)\,t^{n}\,\dfrac{t^{i}-t^{-i}}{t-t^{-1}}   & \text{if $i \le j=n$},
\\[1ex]
\max(d_i,d_j)\displaystyle\sum_{s=1}^{i}( t^{j-i+2s-1} +  t^{2n -j-i+2s-1} ) \\
\hs{8ex}=\max(d_i,d_j)\,t^{n}\,\dfrac{(t^{n-j}+t^{j-n})(t^i-t^{-i})}{t-t^{-1}}
  &\text{if $i\le j<n$.}
\ec
\eneq
\end{theorem}

\begin{proof}
By Theorem~\ref{thm: inv}  and Corollary~\ref{cor: computation setting},
\eqref{eq: BC formula} comes from Lemma~\ref{lem: the Q B} and  Lemma~\ref{lem: the Q C} which  gives  explicit computations for $\tau^k_{Q^\circ}(\ga_j^{Q^\circ})$ $(0 \le k <n, \ j \in I)$ in terms of $\lan a, \pm b\ran$ ($1 \le a \le b \le n$)   and  $\lan c \ran$ ($1 \le c \le n$).
\end{proof}

\begin{example}
For $\bDynkin$ of types $B_3$ and $C_3$, $\tuB^{B_3}(t)$ and $\tuB^{C_3}(t)$ are given as follows:
\begin{align*}
 \tuB^{B_3}(t) &=  \scriptstyle{ \dfrac{1}{1+t^6}
\left(\begin{array}{lll}
2( t^{5} +  t) & 2 (t^{4} + t^{2}) & 2 \, t^{3} \\
2 (t^{4} + t^{2}) & 2( t^{5} + 2  t^{3} +  t) & 2 (t^{4} + t^{2}) \\
2 t^{3} & 2 ( t^{4} +  t^{2}) & t^{5} + t^{3} + t
\end{array}\right)},  \allowdisplaybreaks \\
\tuB^{C_3}(t) & = \scriptstyle{  \dfrac{1}{1+t^6} \left(\begin{array}{lll}
t^{5} + t & t^{4} + t^{2} & 2 \, t^{3} \\
t^{4} + t^{2} & t^{5} + 2 \, t^{3} + t & 2 (t^{4} +  t^{2}) \\
2 \, t^{3} & 2 (t^{4} +  t^{2}) & 2 ( t^{5} + t^{3} +  t)
\end{array}\right)}.
\end{align*}
\end{example}

\subsubsection{$F_4$ and $G_2$ cases}The simple roots for these cases can be expressed as follows:
\ben
\item For $F_4$-case, we use the notation $(a,b,c,d)\seteq a\ve_1+b\ve_2+c\ve_3+d\ve_4$ where $\st{\ve_i}_{1 \le  i \le 4}$ is an orthogonal basis with the square length $2$.  Then
$$\akew\al_1 = (0,1,-1,0), \, \al_2 = (0,0,1,-1), \, \al_3 = (0,0,0,1) \; \text{and} \; \al_4 = (1/2,-1/2,-1/2,-1/2).$$
\item For $G_2$-case, we use an orthonormal basis $\{ \ep_i \mid 1 \le i \le 3 \}$ of $\R^3$ and the notation $\ssrt{a,b,c}\seteq a\ep_1+b\ep_2+c\ep_3$. Then
$$ \al_1 = \ssrt{0,1,-1} \quad \text{ and }  \quad \al_2 = \ssrt{1,-2,1}.$$
\ee

The AR-quivers
$\Gamma_{Q^\circ}$ for Dynkin quivers in Example~\ref{Ex: BCFG Dynkin quiver} are depicted as follows:
\begin{align}
& \Gamma^{F_4}_{Q^\circ}=  \hspace{-2ex} \raisebox{4em}{ \scalebox{0.61}{\xymatrix@!C=3.3ex@R=2ex{
(i\setminus p) & -9 & -8 & -7 & -6  & -5 & -4 & -3 & -2 & -1 & 0 & 1 & 2 & 3 & 4\\
1&&&& \sprt{1,0,-1,0}\ar[dr] && \sprt{0,0,1,1} \ar@{->}[dr]   &&   \sprt{1,0,0,-1} \ar@{->}[dr]  &&\sprt{0,1,0,1}\ar@{->}[dr]&& \sprt{0,0,1,-1}\ar@{->}[dr]&&  \sprt{0,1,-1,0} \\
2&&& \sprt{1,-1,0,0} \ar@{->}[dr]\ar@{->}[ur]   &&  \sprt{1,0,0,1} \ar@{->}[dr]\ar@{->}[ur]     &&  \sprt{1,0,1,0} \ar@{->}[dr]\ar@{->}[ur]&& \sprt{1,1,0,0}\ar@{->}[dr]\ar@{->}[ur]&& \sprt{0,1,1,0}\ar@{->}[dr]\ar@{->}[ur]&& \sprt{0,1,0,-1}  \ar@{->}[ul]\ar@{->}[ur] \\
3&&\sprt{\frac{1}{2},-\frac{1}{2},-\frac{1}{2},\frac{1}{2}} \ar@{->}[dr]\ar@{=>}[ur]   &&\sprt{\frac{1}{2},-\frac{1}{2},\frac{1}{2},\frac{1}{2}}\ar@{->}[dr]\ar@{=>}[ur]   &&  \sprt{1,0,0,0}\ar@{->}[dr]\ar@{=>}[ur] && \sprt{\frac{1}{2},\frac{1}{2},\frac{1}{2},\frac{1}{2}}\ar@{->}[dr]\ar@{=>}[ur] &&  \sprt{\frac{1}{2},\frac{1}{2},\frac{1}{2},-\frac{1}{2}}\ar@{->}[dr]\ar@{=>}[ur]   &&\sprt{0,1,0,0} \ar@{=>}[ur]\\
4& \sprt{\frac{1}{2},-\frac{1}{2},-\frac{1}{2},-\frac{1}{2}}   \ar@{->}[ur]  &&  \sprt{0,0,0,1}\ar@{->}[ur]      &&  \sprt{\frac{1}{2},-\frac{1}{2},\frac{1}{2},-\frac{1}{2}}\ar@{->}[ur]  && \sprt{\frac{1}{2},\frac{1}{2},-\frac{1}{2},\frac{1}{2}} \ar@{->}[ur]  &&\sprt{0,0,1,0}\ar@{->}[ur]    & & \sprt{\frac{1}{2},\frac{1}{2},-\frac{1}{2},-\frac{1}{2}} \ar@{->}[ur]
}}}, \label{eq: AR F}
\end{align}
and
\begin{align}
& \Gamma^{G_2}_{Q^\circ}=  \raisebox{2.3em}{ \scalebox{0.8}{\xymatrix@!C=4ex@R=2ex{
(i\setminus p) &   -3 &  -2 & -1 & 0 & 1 & 2\\
1&   \ssrt{0,1,-1}  \ar@{=>}[dr]  && \ssrt{1,0,-1}  \ar@{=>}[dr]     &&\ssrt{1,-1,0} \ar@{=>}[dr] \\
2&&  \ssrt{1,1,-2}  \ar@{-}[ul]\ar@{->}[ur]  && \ssrt{2,-1,-1} \ar@{-}[ul]\ar@{->}[ur]  &&\ssrt{1,-2,1} \ar@{-}[ul]
}}}.\label{eq: AR G}
\end{align}

Using AR-quivers in~\eqref{eq: AR F} and~\eqref{eq: AR G},  we can obtain $(1+t^{\sfh})\tuB_{i,j}(z)$ for these cases:
\fontsize{8}{8}
\eqn
&&(1+t^{12})\,\tuB^{F_4}(t)= \\
&& \hs{-1ex}\left(\begin{array}{llll}
2(t^{11} + t^{7} + t^{5} + t) & 2(t^{10} + t^{8} + 2 t^{6} + t^{4} +t^{2}) & 2( t^{9} + t^{7} +  t^{5} +  t^{3}) & 2(t^{8} + t^{4}) \\
2(t^{10} + t^{8} + 2 t^{6} + t^{4} +t^{2})\akew   & 2(t^{11} + 2t^{9} + 3t^{7} + 3t^{5} +2t^{3} + t)\akew & 2 (t^{10} + 2t^{8} + 2 t^{6} + 2 t^{4} +  t^{2})\akew & 2 (t^{9} +t^{7} + t^{5} + t^{3}) \\
2( t^{9} + t^{7} +  t^{5} +  t^{3})  &  2 (t^{10} + 2t^{8} + 2 t^{6} + 2 t^{4} +  t^{2})\;  \; & t^{11} + 2 t^{9} + 3 t^{7} + 3 t^{5} + 2 t^{3} + t \; \;  & t^{10} + t^{8} + 2 t^{6} + t^{4} + t^{2} \\
2(t^{8} + t^{4}) &  2 (t^{9} +t^{7} + t^{5} + t^{3}) & t^{10} + t^{8} + 2 t^{6} + t^{4} + t^{2} & t^{11} + t^{7} + t^{5} + t
\end{array}\right)
\eneqn  \fontsize{12}{12}
and
\begin{align*}
& (1+t^6)\tuB^{G_2}(t)= \left(\begin{array}{ll}
t^{5} + 2 \, t^{3} + t  & 3 (t^{4} +  t^{2}) \\
3 (t^{4} +  t^{2}) & 3 (t^{5} + 2 t^{3} + t)
\end{array}\right).
\end{align*}

\subsubsection{Remark on the inverses of $\sfC(q)$ and $\sfC(q,t)$} \label{subsubsec: Quantum cartan and two parameter}
Note that $\usfB(t)$ can be obtained from the quantum Cartan matrix $\sfC(q)$
by just replacing $q$ with $t$, when $\g$ is of simply-laced type.
The inverse $\tsfC(q)$ of quantum Cartan matrix $\sfC(q)$ also enjoys the similar properties of
$\tuB(t)$ even in the non simply-laced type. We refer \cite{FO21} for $\tsfC(q)$ of non simply-laced type.

Very recently, in \cite{FM21}, Fujita-Murakami investigated the behavior of the inverse $\tsfC(q,t)$ of $\sfC(q,t)$, as an matrix with entries in $\Z( \hspace{-.3ex}(q,t) \hspace{-.3ex})$, which implies the several properties of $\tuB(t)$ in Corollary~\ref{cor: computation setting}  and an implicit computation of $\tuB(t)$ via the specialization at $q=1$. However, as far as the authors understand,
they do (i) \emph{not} use the Coxeter elements and (ii) \emph{not}
give explicit formulas of entries of $\tsfC(q,t)$ and hence $\tuB(t)$.

\section{Quantum torus associated with a Dynkin diagram} \label{sec: Quantum torus}

In this section, we construct a quantum torus related to
the $t$-quantized Cartan matrix.
Then we investigate the structure of the quantum torus.

\begin{definition} \label{def: new q torus} Let $\ttq$ be an indeterminate,
and let $\calX_\ttq $ be the $\Z[\ttq^{\pm \frac{1}{2} }]$-algebra with the generators
$\{  \tX^{\pm1}_{i,p}  \mid (i,p) \in \hDynkin_0   \}$  and the following
defining relations:
\begin{itemize}
\item
$\tX_{i,p}\tX_{i,p}^{-1}=\tX_{i,p}^{-1}\tX_{i,p}=1$ for any $(i,p) \in \hDynkin_0 $,
\item
$\tX_{i, p}\tX_{j, s} = \ttq^{\ucalN(i,p;j,s)}\tX_{j,s}\tX_{i,p}$ for $(i,p)$, $(j,s) \in \hDynkin_0 $,
\end{itemize}
where we set (see \eqref{eq: teta}):
\eq \label{eq: uN}
\hs{5ex}\ucalN(i,p;\,j,s) &&\seteq \;
\teta_{i,j}(p-s-1)-\teta_{i,j}(p-s+1)\\
&&=\tfb_{i,j}(p-s-1)-\tfb_{i,j}(s-p-1)-\tfb_{i,j}(p-s+1)+\tfb_{i,j}(s-p+1).
\nn\eneq
We call $\calX_\ttq$  the \emph{quantum torus associated with $\usfC(t)$}.
\end{definition}

\begin{remark} \hfill
\bnum
\item For simply-laced type $\Dynkin$, the quantum torus  $\calX_\ttq$ is already defined in \cite{Nak04,VV02,H04} in the context
of the quantum Grothendieck ring of $\scrC^0_{\g}$ of  $U'_\nu(\widehat{\g})$.
More precisely, the $\ttq$-commutation relation in~\eqref{eq: uN} coincides with the \emph{$\ttt$-commutation relation} in~\cite{Nak04,VV02,H04} when $\Dynkin$ is of simply-laced type, which is defined by using quantum Cartan matrix $\cm(q)$ and its inverse $\tsfC(q)$.  
\item Note that
$$\ucalN(i,p;j,s) = \ucalN(j,p;i,s)=-\ucalN(i,s;j,p)=-\ucalN(j,s;i,p)$$
since $\tuB(t)$ is symmetric. Also we have
\begin{align}\label{eq: pos}
\ucalN(i,p;j,s) \seteq \tfb_{i,j}(p-s-1)-\tfb_{i,j}(p-s+1) \qquad  \text{ if } p >s
\end{align}
by Lemma~\ref{lem: basic tfb}.
\ee
\end{remark}

We say that $\tm \in \calX_\ttq$ is a {\em $\calX_\ttq$-monomial} if it is a product of the generators $\tX_{i,p}^{\pm1}, \ttq^{\pm1/2}$.
For a  $\calX_\ttq$-monomial $\tm$,
writing $\tm=\ttq^c\prod_{(i,p)\in\hDynkin_0}\tX_{i,p}^{\nu_{i,p}}$ with $c\in \frac{1}{2}\Z$,
$\nu_{i,p}\in\Z$
(product is taken in some order),
we set $u_{i,p}(\tm)=\nu_{i,p}$.

For $\calX_\ttq$-monomials $\tm$ and $\tm'$, we set
\eq
&&\ucalN(\tm,\,\tm') \seteq \sum_{(j,p),(i,s)\in \hDynkin_0}
u_{j,p}(\tm)u_{i,s}(\tm')\ucalN(j,p;i,s).\label{def:Nmm}\eneq
Hence we have
\eq\tm\;\tm'=\ttq^{\ucalN(\tm,\,\tm')}\tm'\;\tm.\eneq

Note that there exists a $\Z$-algebra anti-involution $\overline{( \cdot )}$ on $\calX_\ttq$ given by
\begin{align}\label{eq: involution}
 \ttq^{1/2 } \mapsto  \ttq^{-1/2}  \quad \text{ and } \quad  \tX_{i,p}   \mapsto  \ttq^{d_i} \tX_{i,p}.
\end{align}
Then, for any $\calX_\ttq$-monomial $\tm \in \calX_\ttq$, there exists a unique $r \in \frac{1}{2}\Z$ such that
$\ttq^{r}\tm$ is $\overline{( \cdot )}$-invariant.
An element of this form is called a \emph{bar-invariant} monomial.

\begin{definition}  [{cf.\ \cite[Definition 5.5]{FHOO}}] \label{def: quantum tori XqQ}
For a  subset $\sfS \subset  \hDynkin_0$, we   denote by  $\calX_{\ttq,\sfS}$
the   quantum subtorus  of $\calX_{\ttq}$ generated by $\tX^{\pm1}_{i,p}$ for $(i,p) \in \sfS \subset \hDynkin_0$.
In particular, for a Dynkin quiver $Q=(\Dynkin,\xi)$, we denote by $\calX_{\ttq,Q}$ the quantum subtorus of $\calX_{\ttq}$ generated by $\tX^{\pm1}_{i,p}$ for $(i,p) \in (\Gamma_Q)_0 \subset \hDynkin_0$.
\end{definition}

The following theorem is a generalization of~\cite[Proposition 3.1]{HL15} and ~\cite[Proposition 5.21]{FO21}:

\begin{theorem}\label{thm: calN} Let $(i,p)$, $(j,s)\in\hDynkin_0$
and let $Q$ be a Dynkin quiver  on $\Dynkin$. Set
$\phiQ(i,p)=(\al,k)$ and $\phiQ(j,s)=(\be,l)$.
Then, we have
\begin{align}\label{eq: calN}
\ucalN(i,p;j,s) =  (-1)^{k+l + \delta( p \ge s ) } \, \delta\bl (i,p) \ne (j,s) \br\;  (\al,\be)
\end{align}
\end{theorem}

\begin{proof}
First we assume that $p > s$. Then by Theorem~\ref{thm: inv} and~\eqref{eq: pos}, we have
\begin{align*}
\ucalN(i,p;j,s)
& =  -\bl\varpi_{i}, \cox{Q}^{(p-s+\xi_{j}-\xi_{i})/2 }(\gamma^{Q}_{j}) - \cox{Q}^{(p-s+\xi_{j}-\xi_{i})/2 -1}(\gamma^{Q}_{j})\br \allowdisplaybreaks\\
& =-  (\cox{Q}^{(\xi_{i} - p)/2} (1-\cox{Q})\varpi_{i}, \cox{Q}^{(\xi_{j}-s)/2}(\gamma^{Q}_{j}))=  - (\cox{Q}^{(\xi_{i} - p)/2}(\gamma^{Q}_{i}), \cox{Q}^{(\xi_{j}-s)/2}(\gamma^{Q}_{j})) \\& =  (-1)^{1+ k+l}(\alpha, \beta).
\end{align*}
By the skew-symmetricity, we obtain the assertion for $p<s$.
When $p=s$, the left hand side obviously vanished
and the right hand side vanishes also   since
$(\al,\beta)=0$ if $i\not=j$ by Theorem~\ref{thm: OS17}~\eqref{it: noncom}.
\end{proof}

For $a,b\in\Z$ such that $a\le b$ and $i \in I$, we define a $\calX_\ttq$-monomial $\tm^{(i)}[a,b]$  by
\eq  \tm^{(i)}[a,b] = \prod_{(i,p) \in \hDynkin_0, \; a\le p\le b }   \tX_{i,\,p}.\label{def:m}
\eneq

Note that (see \eqref{def:Nmm})
\eqn
\ucalN\bl\tm^{(i)}[p,p'],\; \tm^{(j)}[s, s']\br
\seteq\sum_{\substack{(i,x),\;(j,y)\in\hDynkin_0;\\
p\le x\le p',\ s\le y\le s'}}\ucalN(i,x;j,y).
\eneqn

The following proposition is a generalization of  \cite[Proposition 8.4]{FHOO}:

\begin{proposition} \label{prop: NnKR}
Let $Q=(\Dynkin, \xi)$ be a Dynkin quiver.
Let $(i, p), (i,p'), (j, s), (j, s') \in \hDynkin_0$
with $p \le p'$ and $s \le s'$. Assuming
$p - s \le d(i,j)$ and $s' - p' \le d(i,j)$,
we have
\eq\label{eq:KR coeff}
&&\ba{l}\ucalN\bl\tm^{(i)}[p,p'],\; \tm^{(j)}[s, s']\br\\[1ex]
\hs{13ex}= \Bigl(\cox{Q}^{(\xi_{i} -p)/2 + 1}\varpi_{i} + \cox{Q}^{(\xi_{i} - p')/2}\varpi_{i},\;
\cox{Q}^{(\xi_{j}-s)/2 + 1}\varpi_{j} - \cox{Q}^{(\xi_j -s')/2}\varpi_j\Bigr).
\ea\eneq
\end{proposition}

\begin{proof}
Recall that $u\lec v$ means that $u\le v$ and $u\equiv v\bmod 2$.

By the definition of $\ucalN\bl\tm^{(i)}[p,p'],\; \tm^{(j)}[s, s']\br$, we have
\eqn
\ucalN\left(\tm^{(i)}[p,p'], \tm^{(j)}[s, s'] \right)
&&= \sum_{x:\;p\lec\, x\,\lec p'}\hs{2ex}\sum_{y:\;s\lec\,y\,\lec s'}  \ucalN(i, x ; j, y) \allowdisplaybreaks\\
&& = \sum_{x:\;p\lec\, x\,\lec p'}\hs{2ex}\sum_{y:\;s\lec\,y\,\lec s'} \bl \teta_{i,j}\bl(x-y-1)-\teta(x-y+1)\br\\
&&= \sum_{y:\;s\lec\,y\,\lec s'}\bl \teta_{i,j}(p-y-1)-\teta_{i,j}(p'-y+1)\br \\
&&=  \sum_{y:\;s\lec\,y\,\lec s'}\bl\tfb_{i,j}(y- p+1)  - \tfb_{i,j}(p'-y+1)\br,
\eneqn
where the last equality follows from Proposition~\ref{prop:Xi}.
On the other hand, we compute
\begin{align*}
&\Bigl(\cox{Q}^{(\xi_{i} -p)/2 + 1}\varpi_{i} + \cox{Q}^{(\xi_{i} - p')/2}\varpi_{i},  \cox{Q}^{(\xi_{j}-s)/2 + 1}\varpi_{j} - \cox{Q}^{(\xi_j -s')/2}\varpi_j\Bigr) \allowdisplaybreaks \\
&\hs{5ex}=\Bigl(\cox{Q}^{(\xi_{i} -p)/2 + 1}\varpi_{i} + \cox{Q}^{(\xi_{i} - p')/2}\varpi_{i},  \sum_{y:\;s\lec\,y\lec\, s'}(\cox{Q}^{(\xi_{j}-y)/2+1}\varpi_j- \cox{Q}^{(\xi_j-y)/2}\varpi_j)\Bigr) \allowdisplaybreaks \\
&\hs{5ex}= -\sum_{y:\;s\lec\,y\lec\, s'}  \Bigl(\cox{Q}^{(\xi_{i} -p)/2 + 1}\varpi_{i} + \cox{Q}^{(\xi_{i} - p')/2}\varpi_{i},  \cox{Q}^{(\xi_{j}-y)/2}\gamma_{j}^{Q}\Bigr) \allowdisplaybreaks\\
&\hs{5ex}= -\sum_{y:\;s\lec\,y\lec\, s'}\left(\varpi_{i}, \cox{Q}^{((p-y-1) + \xi_{j} -\xi_{i}-1)/2}\gamma^{Q}_{j}  + \cox{Q}^{((p'-y+1) + \xi_{j} -\xi_{i}-1)/2}\gamma^{Q}_{j}\right)  \allowdisplaybreaks\\
&\hs{5ex}\underset{(1)}{=}-\hs{-1ex}\sum_{y:\;s\lec\,y\lec\, s'}\Bigl(\tfb_{i,j}(p-y-1) - \tfb_{i,j}(y-p+1) +\tfb_{i,j}(p'-y+1) - \tfb_{i,j}(y-p'-1) \Bigr)\allowdisplaybreaks\\
&\hs{5ex}\underset{(2)}{=} -\sum_{y:\;s\lec\,y\lec\, s'}\bigl(\tfb_{i,j}(p'-y+1) - \tfb_{i,j}(y-p+1)\bigr),
\end{align*}
where $\underset{(1)}{=}$ follows from~\eqref{eq: tfb general u}, and we again used Proposition~\ref{prop:Xi} for $\underset{(2)}{=}$.
From the above two computations, we obtain the conclusion.
\end{proof}

For a $\calX_\ttq$-monomial $\tm$, we define
\begin{align*}
\wt_Q(\tm) &  \seteq \sum_{(i,p) \in \hDynkin_0}  u_{i,p}(\tm) \pi(\phiQ(i,p)) \in  \rl.
\end{align*}
Here we set
$$\text{$\pi\bl(\beta,k)\br \seteq (-1)^k\be$ for  $(\be,k)\in\hPhi^+$.}$$
 We call
$\wt_Q(\tm)$ the \emph{$Q$-weight of $\tm$}.
With this definition, \eqref{eq: calN} reads as
\eq\label{eq:Nij}
\ucalN(i,p;j,s) =  (-1)^{ \delta( p \ge s ) }
 \delta\bl (i,p) \ne (j,s)\br\; \bl\wt_Q(\tX_{i,p}),\,\wt_Q(\tX_{j,s})\br.
\eneq

\medskip
Recall that $\cm=(\sfc_{i,j})_{i,j}$ is the Cartan matrix.
For $(i,p+1)\in\hDynkin_0$,  we set
\eq \label{eq: tB}
\tB_{i,p}\seteq\ttq^{n_{i,p}}\tX_{i,p-1}\tX_{i,p+1}\prod_{j:\;d(i,j)=1}
\tX_{j,p}^{\ \sfc_{j,i}},
\eneq
where we choose $n_{i,p}\in\frac{1}{2}\Z$ such that
$\tB_{i,p}$ is a bar-invariant $\calX_\ttq$-monomial.

The following corollary (cf. \cite[(5.3)]{FO21}) follows from
Theorem~\ref{thm: additive}.
\begin{corollary} \label{lem: A zero}
For $(i,p) \in I \times \Z$ with $(i,p+1) \in \hDynkin_0$ and any Dynkin quiver $Q$, we have
$$\wt_Q(\tB_{i,p})  =0.$$
\end{corollary}
\Proof
It is easy to see that
\eq\pi\bl\wh{w}^m(\beta,k)\br=w^m\bl\pi(\beta,k)\br
\qt{for any $w\in\weyl$, $m\in\Z$ and $(\beta,k)\in\hPhi^+$.}
\eneq
Hence, we have
\eqn
\wt_Q(\tB_{i,p})&&=
\pi\bl\phiQ(i,p-1)\br+\pi\bl\phiQ(i,p+1)\br+\sum_{j:\;d(i,j)=1}\sfc_{j,i}
\pi\bl\phiQ(j,p)\br\\
&&=\pi\bl\wcox\,{}^{(\xi_i-p+1)/2}(\ga_i^Q,0)\br
+\pi\bl\wcox\,{}^{(\xi_i-p-1)/2}(\ga_i^Q,0)\br
+\hs{-1.8ex}\sum_{j:\;d(i,j)=1}\hs{-1.6ex}\sfc_{j,i}\pi\bl\wcox\,{}^{(\xi_j-p)/2}(\ga_j^Q,0)\br\\
&&=\tauQ{}^{(\xi_i-p+1)/2}(\ga_i^Q)
+\tauQ\,{}^{(\xi_i-p-1)/2}(\ga_i^Q)
+\sum_{j:\;d(i,j)=1}\sfc_{j,i}\tauQ\,{}^{(\xi_j-p)/2}(\ga_j^Q),
\eneqn
which vanishes by Theorem~\ref{thm: additive}.
\QED

\begin{example} Let $Q$ be a Dynkin quiver of type $B_3$ in~\eqref{eq: B3 Q}.
Then we have
$$
\wt_Q(\tX_{1,1})= \srt{2,-3},  \   \wt_Q(\tX_{2,2})= \srt{1,-3}, \     \wt_Q(\tX_{2,0})= \srt{1,2},\     \wt_Q(\tX_{3,1})= \srt{1},
$$
by using~\eqref{eq: B3 ex}, and hence
$$
\wt_Q(\tX_{2,0} \tX_{2,2} \tX_{1,1}^{-1} \tX_{3,1}^{-2})=\wt_Q(\tB_{2,1})=0.
$$
 Also, one can check that all  the monomials in the following element in $\calX_q$ have the same $Q$-weight as $\ve_2-\ve_3$:
\begin{equation}
\begin{aligned} \label{eq: elt in frakK B3}
&q\tX_{1,1} +q \tX_{2,2}\tX_{1,3}^{-1} + q^2 \tX_{3,3}^2\tX_{2,4}^{-1} + (q^{-1} +q)\tX_{3,3}\tX_{3,5}^{-1} \\
&\hspace{30ex} + q^2 \tX_{2,4} \tX_{3,5}^{-2} +  q \tX_{1,5}\tX_{2,6}^{-1} + q^{-1} \tX_{1,7}^{-1}.
\end{aligned}
\end{equation}
\end{example}
 
The following proposition is a generalization of \cite[Proposition 3.12]{H04}:
\begin{proposition} \label{prop: YA com}
For $i,j\in I$ and $p,s,t,u \in \Z$ with $(i,p), (j,s+1), (i,t+1), (j,u+1) \in \hDynkin_0$, we have
$$
\tX_{i,p} \tB_{j,s}^{-1} = \ttq^{\,\be(i,p;j,s)}\,\tB_{j,s}^{-1}\tX_{i,p}
\quad \text{and} \quad
\tB^{-1}_{i,t}\tB_{j,u}^{-1} = \ttq^{\,\al(i,t;j,u)}\,\tB_{j,u}^{-1}\tB_{i,t}^{-1}.
$$
Here,
\eq
\be(i,p;j,s) &&=   \delta(i=j) (-\delta(p-s=1)+\delta(p-s=-1))(\al_i,\al_i),
\label{eq;beta} \allowdisplaybreaks \\ [1ex]
\al(i,t;j,u)  &&=
\bc
\pm(\al_i,\al_i)&\text{if $(i,t)=(j,u\pm2)$,}\\
\pm2(\al_i,\al_j)&\text{if $d(i,j)=1$ and $t=u\pm1$,}\\
0&\text{otherwise.}\ec\label{eq;alpha}
\eneq
\end{proposition}

\begin{proof}
Let us prove first \eqref{eq;beta}.

The indices in $ \hDynkin_0$ appearing in
$\tB_{j,s}$ can be described as follows:
$$
 \raisebox{2.5em}{ \xymatrix@!C=8ex@R=1ex{
&  (j'_1,s) \ar[dr]^{-\sfc_{j'_1,j}} \\
 (j,s-1)  \ar[ur]^{-\sfc_{j,j'_1}} \ar[dr]_{-\sfc_{j,j'_r}} & \vdots & (j,s+1) \\
& (j'_r,s)  \ar[ur]_{-\sfc_{j'_r,j}}
 }} \quad \text{ for } j'_k \text{ with } d(j'_k,j)=1.
$$

\noi
Let us first assume that (a) $|p - s| > 1$ or (b) $d(i,j)> 1$ and $|p - s| \le 1$. In these case,  we have
$$   \{ (i,p) \}  \cap \Bigl(\st{ (j,s+1),  (j,s-1)}\cup
\st{(j',s) \mid i'\in \Dynkin_0, \  d(j',j)=1}\Bigr)   = \emptyset.$$

\bna
\item In this case, we have $p > s+1$ or $p < s-1$. Then \eqref{eq:Nij}
 implies that
$$  \ucalN \bl \tX_{i,p}, \tB_{j,s}^{-1} \br  =  (-1)^{\delta(p>s+1)}   \bl \wt_Q(\tX_{i,p}) , \wt_Q(\tB_{j,s}^{-1} ) \br=0. $$
where the last equality follows from Corollary~\ref{lem: A zero}.

\item  In this case,  there are no path between $(i,p)$ and $(j,s \pm1)$,
as well as between $(i,p)$ and $(j',s)$ with $d(j',j)=1$.
Thus Theorem~\ref{thm: OS17}~\eqref{it: noncom} tells that
$$  \bl \wt_Q(\tX_{i,p}) ,  \wt_Q(\tX_{j,s \pm 1}) \br =   \bl \wt_Q(\tX_{i,p}) , \wt_Q(\tX_{j',s})  \br =0.$$
Thus $\ucalN(   \tX_{i,p}, \tB_{j,s}^{-1} ) =0$.
\ee

\noi

\bna
\item[({\rm c})]Assume that $d(i,j)=1$ and $p=s$. Then, the relevant terms are
$$
\xymatrix@!C=7ex@R=2ex{
&(j',s)\ar[dr]\\
 (j,s-1) \ar[dr]\ar[ur]&& (j,s+1)&&d(j',j)=1,\ j'\not=i. \\
& (i,s) \ar[ur]
 }
$$

In this case, Theorem~\ref{thm: OS17}~\eqref{it: noncom} tells that
$$  \bl \wt_Q( \tX_{i,s}), \wt_Q(\tX_{j',s} ) \br =0.$$
Also, \cite[Proposition 2.15]{OS19A} tells that
$$  \bl \wt_Q( \tX_{i,s}), \wt_Q(\tX_{j,s-1} ) \br  =\bl \wt_Q( \tX_{i,s}), \wt_Q(\tX_{j,s+1} ) \br.  $$
Then Theorem~\ref{thm: calN}  tells that
\eqn
&& \ucalN \bl \tX_{i,s}, \tB_{j,s}^{-1} \br   = -\sum_{ j' \ne i,j }\sfc_{j',j} (-1)^{\delta(s\ge s)}\bl \wt_Q( \tX_{i,s}), \wt_Q(\tX_{j',s} ) \br \\
&&\hs{9ex} - (-1)^{\delta(s\ge s+1)}\bl  \wt_Q(\tX_{i,s}), \wt_Q(\tX_{j,s+1})  \br -(-1)^{\delta(s\ge s-1)} \bl  \wt_Q(\tX_{i,s}), \wt_Q(\tX_{j,s-1})  \br\\
&&\hs{4ex}=- \bl  \wt_Q(\tX_{i,s}), \wt_Q(\tX_{j,s+1})  \br + \bl  \wt_Q(\tX_{i,s}), \wt_Q(\tX_{j,s-1})  \br =0.
\eneqn
\ee

\item[(d)]\ Now let us assume that $i=j$ and $|p -s|=1$, which is the only remained case.
Then, $(i,p)$ is one of $(j,s \pm 1)$ and
$$
({\rm 1}) \hs{.4ex} \raisebox{1.7em}{ \xymatrix@!C=7ex@R=2ex{
(i,p)= (j,s-1) \ar[dr]&& (j,s+1) \\
& (j',s) \ar[ur]
 }}\hs{1ex} \text{or}\hs{1ex}  ({\rm 2}) \hs{.4ex} \raisebox{1.7em}{ \xymatrix@!C=7ex@R=2ex{
 (j,s-1) \ar[dr]&& (i,p)=(j,s+1) \\
& (j',s) \ar[ur]
 }}\hs{-1.8ex} .
$$

\begin{enumerate}
\item[(1)]
In this case, Theorem~\ref{thm: calN}  tells that
\eqn
\ucalN \bl \tX_{i,s-1}, \tB_{j,s}^{-1} \br  && =- \sum_{ j'\not=j } \sfc_{j',j}
(-1)^{\delta(s-1\ge s)} \bl \wt_Q( \tX_{j,s-1}), \wt_Q(\tX_{j',s} ) \br \\
&&\hs{10ex}- (-1)^{\delta(s-1\ge s+1)}\bl  \wt_Q(\tX_{j,s-1}), \wt_Q(\tX_{j,s+1}) \br \\
&&=  \bl \wt_Q( \tX_{j,s-1}),  -\wt_Q(\tB_{j,s})+\wt_Q(\tX_{j,s-1})  \br \\
&&= \bl \wt_Q( \tX_{j,s-1}),  \wt_Q(\tX_{j,s-1}) \br.
\eneqn
Here the  third equality holds by Corollary~\ref{lem: A zero}.   Finally the assertion for this case is completed by~\eqref{eq:length}.

\item[(2)] In the case $p=s+1$, we have
\eqn
 \ucalN \bl \tX_{i,s+1}, \tB_{j,s}^{-1} \br  && = -\sum_{  j'\not=j }  \sfc_{j',j}
(-1)^{\delta(s+1\ge s)} \bl \wt_Q( \tX_{j,s+1}), \wt_Q(\tX_{j',s} )\br \allowdisplaybreaks\\
&&\hs{5ex} -(-1)^{\delta(s+1\ge s-1)}\bl  \wt_Q(\tX_{j,s+1}), \wt_Q(\tX_{j,s-1}) \br \allowdisplaybreaks\\
&&=  \bl \wt_Q( \tX_{j,s+1}), \wt_Q(\tB_{j,s}) - \wt_Q(\tX_{j,s-1})  \br\allowdisplaybreaks \\
&&=   - \bl \wt_Q( \tX_{i,s+1}),  \wt_Q(\tX_{i,s-1}) \br.
\eneqn
\ee

\snoi
Thus we complete  the proof of \eqref{eq;beta}.

\medskip
Now, let us prove \eqref{eq;alpha}.
We may assume that $(i,t)\not=(j,u)$.

Assume that
$$   \ucalN \bl \tX_{i',k}, \tB_{j,u}^{-1} \br\not=0    $$
for some factor $\tX_{i',k}$ of $\tB_{i,t}^{-1}$.
Then, \eqref{eq;beta} implies that we have
either
(1) $i=j$ and $|t - u| = 2$, or (2) $d(i,j)=1$ and  $|u-t|=1$.

Hence, we can assume (1) or (2).

\snoi
(1)\ Assume that $i = j$ and $|t - u| = 2$.
Then \eqref{eq;beta} implies that
\begin{align*}
\al(i,t;j,u) &=\bc
 -\beta(i,u-1;j,u ) = - (\al_i,\al_i)  & \text{if $(i,t)=(j,u-2)$,}\\[1ex]
 -\beta(i,u+1;j,u) = (\al_i,\al_i)  & \text{if $(i,t)=(j,u +2)$.}
\ec
\end{align*}

\snoi
(2)\ Assume that $d(i ,j) = 1$ and $|t - u| = 1$.
By~\eqref{eq;beta}, we have
\begin{align*}
\al(i,t;j,u)   & = \bc
-  \sfc_{j,i} \beta(j,u-1;j,u )=- \sfc_{j,i}(\al_j,\al_j)= -2 (\al_i,\al_j)
& \text{if $t-u =-1$,} \\[1ex]
-\sfc_{j,i}\beta(j,u+1;j,u)=\sfc_{j,i}(\al_j,\al_j)= 2 (\al_i,\al_j) & \text{if $t-u =1$, }
\ec
\end{align*}
as we desired.
\end{proof}

\begin{definition} \label{def: virtual quantum G ring}
For $i\in\Dynkin_0$, we denote by $\frakK_{i,\ttq}$
the $\Z[\ttq^{\pm \frac{1}{2}}]$-subalgebra of $\calX_\ttq$ generated by
\begin{align} \label{eq: generators}
\tX_{i,l}(1+\ttq^{-d_i}\tB_{i,l+1}^{-1}), \quad \tX_{j,l} \ \  (j\in\Dynkin_0\setminus\st{i}, \; l \in \Z).
\end{align}
We set
$$\frakK_\ttq \seteq  \bigcap_{i \in I} \frakK_{i,\ttq}$$
and call it the \emph{quantum virtual Grothendieck ring associated to $\usfC(t)$}.
\end{definition}

 \begin{remark} \label{rmk: motivation}
When $\sfg$ is associated with $\Dynkin$ of simply-laced type, $\frakK_\ttq$ is isomorphic to the quantum Grothendieck ring $\frakK_t(\sfg) \seteq K_\ttt(\scrC^0_{\widehat{\sfg}})$ invented in
\cite{Nak04,VV02,H04}. Following the construction of $\frakK_\ttt(\g)$ for all $\g$ in~\cite{H04},
$\frakK_{i,\ttt}$ is isomorphic to the kernel of $t$-deformed screening operator $S_{i,\ttt}$ on $\frakK_\ttt(\g)$.
Then it is proved that (i) $\frakK_{i,\ttt}$ is generated by the elements corresponding to~\eqref{eq: generators},
(ii) $ K(\scrC^0_{\widehat{\g}})  \simeq  \frakK_\ttt(\g)  |_{\ttt=1}$ (see also \cite{FM01}).
In the construction of $\frakK_\ttt(\g)$, the quantum Cartan matrix $\sfC(q)$ of type $\g$ and its inverse $\tsfC(q)$ are crucially used.
In Definition~\ref{def: virtual quantum G ring}, we define the $\usfC(t)$-analogue $\frakK_\ttq$ of $K_\ttt(\scrC^0_{\hg})$, which is new one for non simply-laced $\g$, to the best knowledge of the authors.  As we mentioned in the introduction, the specialization of $\frakK_\ttq$ at $\ttq=1$ coincides with the
the specialization $\overline{\calK}_{t}(\g)$ of $\overline{\calK}_{q,t}(\g)$ at $q=1$ and $\al=d$, which is commutative and defined in \cite{FHR21} (see~\cite{JLO22} also).
For non simply-laced $\g$, it is proved in \cite[Theorem 4.3]{FHR21} that $\overline{\calK}_{t}(\g)$ is a homomorphic image of $\overline{\calK}_{t}(\sfg)$, where $\sfg$
is of simply-laced type and contains $\g$ as its non-trivial subalgebra. Hence the specialization $B_{i,p}$ of the monomial  $\tB_{i,p}$ in~\eqref{eq: tB} at $\ttq=1$ can be understood
as an image of $A_{i,p}$ under the surjection, which is deformation of the simple root $\al_i$ in  the $q$-character theory (\cite{FR99,FR99B,FM01}).
In \cite{JLO22}, the ring $\frakK_\ttq$ will be investigated in more precise way.
\end{remark}

Note that one can check that the element in~\eqref{eq: elt in frakK B3} is contained in $\frakK_\ttq$ of type $B_3$.

\section{Unipotent quantum coordinate algebra and the quantum tori isomorphism} \label{Sec: U alg and iso}

In this section, we review the unipotent quantum coordinate algebra $A_\nu(\n)$  associated to a Dynkin diagram $\Dynkin$ and
the quantum torus $\calT_{\nu,[Q]}$ generated by
quasi-commuting unipotent  quantum minors associated to $[Q]$. Then we shall prove that
the quantum torus $\calT_{\nu,[Q]}$ and $\calX_{\ttq,Q}$ are isomorphic.

\subsection{Quantum group and its representations}  Let $(\cm,\wl,\sr,\cwl,\scr)$ be a finite Cartan datum. Let $\nu$ be an indeterminate.
 We denote by $U_\nu(\g)$ the quantum group associated to the Cartan datum, which is  an  algebra over $\Q(\nu)$
generated by $e_i,f_i$ $(i\in I)$ and $\nu^h$ $(h \in \cwl)$. For $\be \in \rl$, $U_\nu(\g)_\be$ denotes the weight space of $U_\nu(\g)$ with weight $\be$ and we set $\wt(x)=\be$ for $x \in U_\nu(\g)_\be$.
We denote by $U_\nu^+(\g)$ (resp.\ $U_\nu^-(\g)$) the subalgebra of $U_\nu(\g)$ generated by $e_i$ (resp.\ $f_i$) for $i \in I$. For $n \in \Z_{\ge0}$ and $i \in I$, we set $e_i^{(n)} \seteq e_i/[n]_i!$ and
$f_i^{(n)} \seteq f_i/[n]_i!$, where we set
$$ \nu_i \seteq \nu^{d_i},\  [n]_i\seteq\dfrac{\nu_i^n-\nu_i^{-n}}{\nu_i-\nu_i^{-1}} \quad \text{ and } [n]_i !\seteq \prod_{k=1}^n [k]_i. $$

We set $\bbA \seteq \Z[\nu^{\pm1}]$ and denote  by $U_\bbA^\pm(\g)$ the $\bbA$-subalgebra of $U_\nu(\g)^\pm$ generated by $e_i^{(n)}$ (resp.\ $f_i^{(n)}$) for $i \in I$ and $n \in \Z_{\ge 0}$.

Note that $U_\nu(\g)$ is a Hopf algebra with the comultiplication
$$\Delta  \cl U_\nu(\g) \to U_\nu(\g) \otimes U_\nu(\g)$$
given by
$$    \Delta(e_i) = e_i \otimes 1 + \nu_i^{h_i}  \otimes e_i , \quad
  \Delta(f_i) = f_i \otimes \nu_i^{-h_i} + 1  \otimes f_i, \quad
 \Delta(\nu^{h}) = \nu^{h} \otimes \nu^{h}.$$

Let $\varphi$ and $*$ be the $\Q(\nu)$ anti-automorphism of $U_\nu(\g)$ defined by
\begin{align*}
& \varphi(e_i) = f_i,  \quad \varphi(f_i) = e_i,\quad \varphi(\nu^h) = \nu^h \ \ \text{ and } \ \  e_i^* =e_i, \quad f_i^* = f_i, \quad (\nu^h)^* = \nu^{-h}.
\end{align*}

There is also a $\Q$-algebra homomorphism $\ol{\phantom{a} }$ of $\Ug$ given by
$$
\overline{\nu} =\nu^{-1},\quad
\overline{e_i}=e_i, \quad \overline{f_i}=f_i  \qtq \overline{\nu^h}=\nu^{-h}.
$$

For a left $U_\nu(\g)$-module $N$, we denote by $N^r$  the right $U_\nu(\g)$-module $\{ m^r \mid m \in M\}$ with the right action induced by $\varphi$:
$$ (m^r)x = (\varphi(x)m)^r  \qquad \text{ for $m \in M$ and $x \in \Ug$}.$$

Let $M$ be a left $\Ug$-module. For $\la \in \wl$, let $$M_\la\seteq\{ m \in M \mid \nu^h m = \nu^{\lan \la,h \ran}m \text{ for every } h \in \cwl\}$$ be the corresponding weight space.
We say a left $\Ug$-module $M$ \emph{integrable} if (i) $M = \soplus_{\eta \in \wl} M_\eta$ with $\dim M_\eta < \infty$ and (ii) the actions of $e_i$ and $f_i$ are locally nilpotent for all $i \in I$.
We denote by $\Oint(\g)$ the category of integrable left $\Ug$-modules $M$ satisfying that there exist finitely many $\la_1,\ldots,\la_m$ such that $\wt(M) \seteq \{ \eta \in \wl \ | \ M_\eta \ne 0\} \subset \bigcup_j (\la_j+\rl^-).$
The category $\Oint(\g)$ is semisimple with its simple objects being isomorphic to the highest weight modules $V(\la)$ with the highest weight vector $u_\la$ of the highest weight $\la \in \wl^+$.
We denote by $\Oint^r(\g)$ the category of integrable right $\Ug$-modules $M^r$ such that $M \in \Oint(\g)$. Then
$\Oint^r(\g)$ is also semisimple with its simple objects being isomorphic to the highest weight modules $V^r(\la)\seteq\bl V(\la)\br^r$ with the highest weight vector $u^r_\la$ for some $\la \in \wl^+$.

Note that the highest weight module $V(\la)$ $(\la \in \wl^+)$ has a unique non-degenerate symmetric bilinear form $( \ , \ )_\la$ such that
$$ (u_\la,u_\la)_\la=1 \quad \text{and}\quad (xu,v)_\la = (u,\varphi(x)v)_\la \quad \text{ for } u,v \in V(\la) \text{ and } x \in \Ug.$$

The tensor product of $\Q(\nu)$-spaces $V^r(\la) \otimes V(\la)$ $(\la \in \wl^+)$ has the natural structure of a $\Ug$-bimodule given by
$$x \cdot (u^r\otimes v) \cdot y = (u^r \cdot y) \otimes (x \cdot v)$$
and the bilinear form $( \ ,  \ )_\la$ on $V(\la)$ induces the non-degenerate bilinear form   $$\lan \cdot,\cdot\ran_\la
\cl V^r(\la) \times V(\la)\to\Q(\nu)$$
given by  $\lan u^r,v\ran_\la=(u,v)_\la$.

\subsection{Unipotent quantum coordinate algebra and quantum minors}
Let $\Ug^*$ be the space $$\Hom_{\Q(\nu)}(\Ug,\Q(\nu)).$$ Then the comultiplication $\Delta$ induces the multiplication on $\Ug^*$ as follows:
$$ \Ug^* \otimes \Ug^* \to (\Ug \otimes \Ug)^* \overset{\Delta^*}{\to} \Ug^*. $$
Namely, for $f,g \in \Ug^*$ and $x \in \Ug$, we have
$$  (fg) (x) = f(x_{(1)})g(x_{(2)}),$$
where $\Delta(x)=x_{(1)} \otimes x_{(2)}$ written by Sweedler's notation.

Note that $\Ug^*$ has also a $\Ug$-bimodule structure given by
$$  x  \cdot (fg) \cdot y = (x_{(1)} \cdot f \cdot y_{(1)})(x_{(2)} \cdot g \cdot y_{(2)})  \ \text{ for } f,g \in \Ug^* \text{ and } x,y \in \Ug,$$
where $(x\cdot f \cdot y)(z) \seteq f(yzx))$ for $z \in \Ug$.

\begin{definition}
We define the \emph{quantum coordinate algebra} $A_\nu(\g)$ as follows:
$$  A_\nu(\g) = \{   f \in \Ug^* \mid  \Ug f \text{ belongs to } \Oint(\g) \text{ and } f \Ug \text{ belongs to } \Oint^r(\g)  \}.$$
\end{definition}

We sometimes denote by $e_i^*$ and $f_i^*$ the operators on $A_\nu(\g)$ acting at the right.

\begin{proposition}[{\cite[Proposition 7.2.2]{Kas93}}]\label{prop: P-W-Kas} We have an isomorphism $\Upphi$ of $\Ug$-bimodules
$$  \soplus_{\la \in \wl^+}  V(\la) \otimes V^r(\la) \isoto A_\nu(\g)  $$
given by
$$
\Upphi(u \otimes v^r)(x) = \lan v^r,\,x\cdot u\ran_\la =\lan v^r \cdot x, u\ran_\la =(v,\,x\cdot u)_\la
$$
for $u \in V(\la)$, $v^r \in V^r(\la)$ and $x \in \Ug$.
\end{proposition}

By Proposition~\ref{prop: P-W-Kas}, $\Ag$ admits a biweight space decomposition $\Ag = \soplus_{\eta,\zeta\in \wl} \Ag_{\eta,\zeta}$
where
$$
\Ag_{\eta,\zeta} \seteq \{ \psi \in \Ag \mid  \nu^{h_\lt} \cdot \psi \cdot  \nu^{h_\rt} = \nu^{\lan h_\lt,\eta\ran+\lan h_\rt,\zeta\ran}\psi \}.
$$
For $\phi \in \Ag_{\eta,\zeta}$, we write $ \lwt(\phi)=\eta$ and $\rwt(\phi)=\zeta$.

\begin{definition}
For $\varpi \in \wl^+$ and $\mu,\zeta \in \sfW\varpi$, we defined the \emph{generalized quantum minor} $\qm(\mu,\zeta)$ as follows:
$$ \qm(\mu,\zeta) \seteq \Upphi(u_\mu \otimes u_\zeta^r) \in \Ag_{\mu,\zeta}.$$
\end{definition}

For elements $y,y' \in \Ag$, we write $y \equiv y'$ if there exists $r \in \Z$ such that $y = \nu^ry'$.

\begin{lemma} \cite[(9.13)]{BZ05} \label{lem: BZ}
For $u,v \in \weyl$ and $\la,\mu \in \wl^+$, we have
$$   \qm(u\la,v\la)\qm(u\mu,v\mu) \equiv \qm(u(\la+\mu),v(\la+\mu))   $$
\end{lemma}

The following proposition is $\nu$-analogue of \cite[Theorem 1.17]{FZ99} and can be proved by following the same arguments in \cite[\S 3.2, 3.3]{GLS13}:
\begin{proposition}
For $u,v \in \weyl$, assume that $\ell(us_i)=\ell(u)+1$ and $\ell(vs_i)=\ell(v)+1$. Then we have
\begin{equation} \label{eq: T-system among delta}
\begin{aligned}
& \qm(us_i(\varpi_i),vs_i(\varpi_i))\qm(u(\varpi_i),v(\varpi_i)) = \\
& \hspace{10ex} \nu^{-d_i}\qm(us_i(\varpi_i),v(\varpi_i))\qm(u(\varpi_i),vs_i(\varpi_i)) + \prod_{j\ne i} \qm(u(\varpi_j),v(\varpi_j))^{-\sfc_{j,i}}.
\end{aligned}
\end{equation}
Here $\qm(u\la,v\la) = \prod_{j\ne i} \qm(u(\varpi_j),v(\varpi_j))^{-\sfc_{j,i}}$ for $\la = s_i\varpi_i+\varpi_i$ by the {\rm Lemma~\ref{lem: BZ}}.
\end{proposition}

The tensor product $\Upg \otimes \Upg$ has the algebra structure defined by
$$ (x_1 \otimes x_2) \cdot (y_1 \otimes y_2) = \nu^{-(\wt(x_2),\wt(y_1))} (x_1y_1\otimes x_2y_2)  $$
for homogeneous elements $x_1,x_2,y_1,y_2 \in \Upg$. We define the $\Q(\nu)$-algebra homomorphism $\Delta_\n\cl\Upg \to \Upg \otimes \Upg$ by
$$  \Delta_\n(e_i) =e_i \otimes 1 +1 \otimes e_i  \quad \text{ for $i \in I$}. $$

We set
$$\An \seteq \soplus_{\be \in \rl^-} \An_\be  \quad \text{ where }  \An_\be \seteq \Hom_{\Q(\nu)}(\Upg_{-\be},\Q(\nu)).$$

Defining the bilinear form $\lan \cdot,\cdot \ran \cl  (\An \otimes \An) \times (\Upg \otimes \Upg) \to \Q(\nu) $ by
$$ \lan \psi\otimes \theta , x \otimes y \ran = \theta(x)\psi(y), $$
the multiplication on $A_\nu(\n)$ is given by
$$
(\phi \cdot \theta)(x) = \lan \phi \otimes \theta,\Delta_\n(x) \ran  \quad\text{ for } \phi,\theta \in A_\nu(\n), \  x \in \Upg.
$$
The algebra $\An$ is called the \emph{unipotent quantum coordinate algebra}. We set
$$ A_{\bbA}(\n) =\{ \psi \in \An \mid \psi(U_\bbA^+(\g)) \subset \bbA \}.$$

We define the bar-involution $\bar{\quad}$ of $\An$ by
$$ \overline{\psi}(x) = \overline{\psi(\bar{x})} \quad \text { for } \psi\in \An \text{ and } x \in \Upg.$$

Note that Lusztig \cite{Lus90,Lus91} and Kashiwara \cite{K91} have constructed a specific $\bbA$-basis   $\bfB^\up$  of $\An$, which is called the \emph{dual canonical/upper global basis}.

The homomorphism $p_\n\cl \Ag \to \An$ induced by $\Upg \to \Ug$ is given by
$$ \lan p_n(\psi),x\ran = \psi(x) \quad \text{ for any } x \in \Upg.$$
Then we have
$$ \wt(p_\n(\psi)) = \lwt(\psi)-\rwt(\psi).$$

\begin{lemma} [{\cite[Proposition 8.5.4]{KKKO18}}] \label{lemma: exp pn}
For $\psi,\theta \in \Ag$, we have
$$ p_\n (\psi\theta)= \nu^{(\zeta_1,\zeta_2-\eta_2)} p_\n (\psi)p_\n(\theta)$$
where $\psi \in\Ag_{\eta_1,\zeta_1}$ and $\zeta \in \Ag_{\eta_2,\zeta_2}$.
\end{lemma}

\begin{definition}
For $\varpi \in \wl^+$ and $\mu,\zeta \in \sfW\varpi$, we define the \emph{unipotent quantum minor} $D(\mu,\zeta)$ as follows:
$$ D(\mu,\zeta) \seteq p_\n(\qm(\mu,\zeta)) \in \An_{\zeta-\mu}.$$
\end{definition}
\noindent
By the definition, $D(\mu,\zeta)$ vanishes if $\zeta\not\preceq \mu$.

\smallskip

The behavior of $D(\mu,\zeta)$ is investigated intensively (see \cite{BZ05,GLS13,KKKO18}):

\begin{proposition} \hfill
\ben
\item For $\La \in \wl^+$ and $\mu,\zeta \in \sfW\La$, $D(\mu,\zeta)$ is a member of $\bfB^\up$ if $\mu \preceq \zeta$, and $D(\mu,\zeta)=0$ otherwise.
\item For $\la,\mu \in \wl^+$ and $s,t,s',t' \in \sfW$ such that
\bnum
\item $\ell(s's)=\ell(s')+\ell(s)$ and $\ell(t't)=\ell(t')+\ell(t)$,
\item $s's\la \preceq t'\la$ and $s'\mu \preceq t't\mu$,
\ee
we have
$$   D(s'\mu,t't\mu)D(s's\la,t'\la)= \nu^{(s's\la+t'\la,t't\mu-s'\mu)}D(s's\la,t'\la) D(s'\mu,t't\mu).$$
\item For $u,v \in \sfW$ and $i \in I$ satisfying $u<us_i$ and $v < vs_i$, we have
\begin{equation}\label{eq: T-system}
\begin{aligned}
& \nu^{(vs_i\varpi_i,v\varpi_i-u\varpi_i)}D(us_i\varpi_i,vs_i\varpi_i)D(u\varpi_i,v\varpi_i)  \\
& \hspace{5ex} = \nu^{(v\varpi_i,vs_i\varpi_i-u\varpi_i)-d_i}D(us_i\varpi_i,v\varpi_i)D(u\varpi_i,vs_i\varpi_i)+D(u\la,v\la),
\end{aligned}
\end{equation}
where $\la=s_i\varpi_i+\varpi_i=2\varpi_i -\al_i$.
\ee
\end{proposition}

The equation~\eqref{eq: T-system} is referred to as $T$-system among unipotent quantum minors.  

\medskip
 For a sequence of indices $\tw = \sseq{i_1,\ldots, i_r} \in I^r$, $1 \le  k \le r$ and $i \in I$, we use the following notations:
\begin{equation}\label{eq: notations +,-}
\begin{aligned}
w_{\le k}  &\seteq s_{i_1} \cdots s_{i_k}, \quad w_{\le 0}\seteq  1,  \\
 k^+ &\seteq \min\bl \{ k< j \le r\mid i_j =i_k \} \cup \{ r+1 \} \br, \\
 k^-  &\seteq \max\bl \{ 1 \le  j <k \mid i_j =i_k \} \cup \{  0 \} \br, \\
 k^-(i) & \seteq \max\bl \{ 1 \le  j <k \mid i_j =i \} \cup \{  0 \} \br, \\  k^+(i)  &\seteq \min\bl \{ k < j  \le r \mid i_j =i \} \cup \{  r+1 \} \br.
\end{aligned}
\end{equation}

For a reduced expression $\uw_0=s_{i_1}\cdots s_{i_\ell}$ of the longest element $w_0$ and $0 \le s \le t \le \ell$, we set
$$
D_{\utw_0}[s,t] \seteq \bc D(w_{\le t}\varpi_{i_t},w_{\le s-1}\varpi_{i_s}) &\text{ if } i_t=i_s \text{ and } s \ge 1, \\
D(w_{\le t}\varpi_{i_t},\varpi_{i_t}) &\text{ if } s =0, \\
\mathbf{1} &\text{ otherwise},
\ec
$$
by taking $\utw_0 \seteq (i_1,\ldots, i_\ell)$.
Then \eqref{eq: T-system} can be written as follows: For $1 \le a  <  b \le \ell$ with $i_a=i_b=i$,
\begin{align}\label{eq: T-system2}
D_{\utw_0}[a^+,b]D_{\utw_0}[a,b^-] = D_{\utw_0}[a,b]D_{\utw_0}[a^+,b^-] + \prod_{j;\; \sfc_{j,i}<0} D_{\utw_0}[a^+(j),b^-(j)]^{-\sfc_{j,i}}.
\end{align}
Here we ignore $\nu$-coefficients and understand $D[x,y]=\mathbf{1}$ for  $y < x$.

\subsection{Quantum torus for unipotent quantum coordinate algebra}
Now let us fix a reduced expression $\uw_0 =s_{i_1}\cdots s_{i_\ell}$ of the longest element $w_0 \in \sfW$.
Recall the notation $\be_k^{\uw_0}$ in~\eqref{eq: beta_k uw_0}
and $k^\pm$ in \eqref{eq: notations +,-} for $k \in  [1,\ell]$. For each $\al \in \Phi^+$, there exists a unique $k$ such that $\al =\be_k^{\uw_0}$. Then we define
\begin{align*}
\al^+ &\seteq \bc \be^{\uw_0}_{k^+}  & \text{ if } k^+ \le \ell, \\ 0 & \text{ if } k^+ =\ell+1,  \ec
& \al^- &\seteq \bc \be^{\uw_0}_{k^-}  & \text{ if } k^- \ge 1, \\ 0 & \text{ if } k^- = 0,   \ec  \\
 \la_\al &\seteq w_{\le k}\varpi_{i_k}.
\end{align*}
By the definition, we have
$$ \la_{\al^-} =\la_{\al}+\al $$
for $\al \in \Phi^+$, where we understand $\la_{\al^-}=\varpi_{i_k}$ if $\al^-=0$.

\smallskip
In the following proposition, we follow the convention that
$\be^+ \prec_{[\uw_0]} \al^+ $ is true if $\al^+=0$.

\begin{proposition}[{\cite[Proposition 10.1, Lemma 11.3]{GLS13}, \cite[Proposition 10.4]{GY17}}]  \label{prop: compatible pair2}
Set
$$  J \seteq \Phi^+, \quad J_f \seteq \{ \al \in \Phi^+ \mid \al^+ = 0 \} \quad \text{ and } \quad J_e \seteq J \setminus J_f. $$
Define the $J \times J_e$-integer matrix $\tsfB^{\uw_0}=(b_{\al\be})_{\al \in J, \be \in J_e}$ as
$$ b_{\al\be} =
\bc
1 & \text{ if } \be=\al^+,\\
\sfc_{i,j} & \text{ if } \al\prec_{[\uw_0]}\be\prec_{[\uw_0]}\al^+\prec_{[\uw_0]}\be^+ \text{ and } d(i,j)=1,\\
-1 & \text{ if } \be^+ =\al, \\
-\sfc_{i,j} & \text{ if } \be \prec_{[\uw_0]} \al \prec_{[\uw_0]} \be^+ \prec_{[\uw_0]} \al^+ \text{ and } d(i,j)=1,\\
0 & \text{ otherwise,}
\ec
$$
where $i=\res^{[\uw_0]}(\al)$ and $j=\res^{[\uw_0]}(\be)$.
Define the $J \times J$-skew symmetric matrix $\Lambda^{\uw_0}=(\Lambda^{\uw_0}_{\al\be})_{\al,\,\be \in J}$ by
$$
\Lambda^{\uw_0}_{\al\be} = (\varpi_{i}-\la_\al,\varpi_{j}+\la_{\be})
\quad \text{ for } \be \not\prec_{[\uw_0] } \al.
$$
Then, $\La^{\uw_0}_{\al\be}\in \Z$, and  $(\Lambda^{\utw_0},\tsfB^{\utw_0})$ is a compatible pair, that is
$$
\sum_{\ga \in J}b_{\ga\al}\Lambda_{\ga\be} = (\al,\al)\delta(\al=\be)
$$
for all $\al \in J_e$ and $\be \in J$.
\end{proposition}

\begin{remark} \label{rmk: simu reindexing}
It is easy to prove that
$\al^\pm$  
depends only on the commutation class of $\uw_0$.
Hence, the notations
$\Lambda^{[\uw_0]}$ and $\tsfB^{[\uw_0]}$ make sense.
\end{remark}

Now we define the following quantum torus by using the matrix $\La^{[\uw_0]}=(\La^{[\uw_0]}_{\al\be})_{\al,\be\in \Phi^+}$:

\begin{definition}
The quantum torus $\calT_{\nu,[\uw_0]}$ is the $\Z[\nu^{\pm \frac{1}{2}}]$-algebra given by the set of generators $\{ Y_{\al}^{\pm1} \mid \al \in \Phi^+\}$ and the following relations:
\begin{itemize}
\item $Y_\al Y_\al^{-1} = Y_\al^{-1}Y_\al =1$ for $\al \in \Phi^+$,
\item $Y_\al Y_\be = \nu^{\Lambda_{\al\be}}Y_\be Y_\al =1$ for $\al,\be \in \Phi^+$.
\end{itemize}
\end{definition}

Let $\calA_\nu(\La^{[\uw_0]}, \tsfB^{[\uw_0]})$ be the \emph{quantum cluster algebra} associated to the compatible pair $(\La^{[\uw_0]}, \tsfB^{[\uw_0]})$
in Proposition~\ref{prop: compatible pair2}. The algebra $\calA_\nu(\La^{[\uw_0]}, \tsfB^{[\uw_0]})$
is the $\Z[\nu^{\pm \frac{1}{2}}]$-subalgebra of the quantum torus  $\calT_{\nu,[\uw_0]}$  generated by the union of the elements called \emph{quantum cluster variables}, which are obtained by all possible sequences of \emph{mutations} (see \cite{FZ99,BZ05}
for more details).

For $\al,\be \in \Phi^+$ with $\res^{[\uw_0]}(\al)=\res^{[\uw_0]}(\be)=i$, we set
$$ \ttD^{[\uw_0]}(\al,\be) \seteq D(\la_\al,\la_\be), \qquad \ttD^{[\uw_0]}(\al,0) \seteq D(\la_\al,\varpi_i).$$
By the definition, $ \ttD^{[\uw_0]}(\al,\be)=0$ unless $\be \prec_{[\uw_0]} \al$. We set $ \ttD^{[\uw_0]}(0,0)\seteq 1$ by convention.

\begin{theorem}[{\cite[Theorem 12.3]{GLS13}, \cite[Theorem 10.1]{GY17}}] There exists a $\Q(\nu^{\frac{1}{2}})$-algebra isomorphism
\begin{align}\label{eq: CL}
 {\rm CL}\cl  \Q(\nu^{\frac{1}{2}}) \otimes_{\Z[\nu^{\pm \frac{1}{2}}]}\calA_\nu(\La^{[\uw_0]}, \tsfB^{[\uw_0]}) \isoto  \Q(\nu^{\frac{1}{2}}) \otimes_{\Z[\nu^{\pm \frac{1}{2}}]} A_\bbA(\n)
\end{align}
sending $ Y_\al$ to $\ttD^{[\uw_0]}(\al,0)$ for all $\al \in \Phi^+$.
\end{theorem}

\subsection{An isomorphism between the quantum tori} Recall the quantum torus $\calX_{\ttq,Q}$ for a Dynkin quiver $Q =(\Dynkin,\xi)$ in Definition~\ref{def: quantum tori XqQ}. Note that
the quantum torus has another presentation given by the set of generators $\{ \tm^{(i)}[p,\xi_i]^{\pm1} \mid (i,p) \in (\Gamma_Q)_0 \}$ and the following relations (see \eqref{def:m}):
\begin{itemize}
\item   $\tm^{(i)}[p,\xi_i]^{-1}$ is the inverse of $\tm^{(i)}[p,\xi_i] $ for $(i,p) \in (\Gamma_Q)_0$,
\item $\tm^{(i)}[p,\xi_i] \cdot \tm^{(j)}[s,\xi_j] = \nu^{\upkappa(i,p;j,s)} \tm^{(j)}[s,\xi_j] \cdot \tm^{(i)}[p,\xi_i]$  for $(i,p),(j,s) \in (\Gamma_Q)_0$,
\end{itemize}
where
\eqn
&&\tm^{(i)}[p,\xi_i] =\ttq^{m_{i,p}} \prod_{x:\;p\lec t\lec \xi_i}\tX_{i,t}, \allowdisplaybreaks\\
&&\upkappa(i,p;j,s) = \ucalN\bl\tm^{(i)}[p,\xi_i],\; \tm^{(j)}[s, \xi_j]\br
=\hs{-1ex}\sum_{\substack{x:\;p\lec x\lec \xi_i\\[.4ex]y:\,s\lec y\lec \xi_j}}\ucalN(i,x;j,y).
\eneqn
Here we choose $m_{i,p}\in \frac{1}{2}\Z$ such that
$\tm^{(i)}[p,\xi_i]$ is a bar-invariant $\calX_\ttq$-monomial.

  The following lemma follows from  Lemma~\ref{lem:admissible}.

\begin{lemma} \label{lem: iso tori}
  Assume that $\uw_0$ is adapted to a Dynkin quiver $Q$.
Then, for $\al \in \Phi^+$ with $\res^{[Q]}(\al)=i$, we have the followings:
\bnum
\item We have $\al = \tauQ \al^-$ provided that $\al^- \ne 0$. If $\al^-=0$, then $\al=\ga_i^Q= \varpi_i - \tauQ \varpi_i$.
\item \label{it: la al}We have $\la_\al  = \tauQ^{(\xi_i-p)/2+1}\varpi_i$, where $\phiQ(i,p)=(\al,0)$.
\ee
\end{lemma}

\begin{theorem} \label{thm: tori iso}
For each Dynkin quiver $Q=(\Dynkin,\xi)$,
there is   an   algebra isomorphism
$$
\Uppsi_Q \cl \calT_{\nu,[Q]} \to\calX_{\ttq,Q}
$$
given by
$$  \nu^{\pm \frac{1}{2}} \mapsto \ttq^{\pm \frac{1}{2}}, \qquad Y_\al \mapsto \tm^{(i)}[p,\xi_i] $$
for all $\al\in \Phi^+$, where $\phiQ(i,p)=(\al,0)$.
\end{theorem}

\begin{proof}
It suffices to show
$$   \upkappa(i,p;j,s) = \Lambda_{\al\be}   $$
for $\al,\be\in \Phi^+$ such that $\phiQ(i,p)=(\al,0)$ and $\phiQ(j,s)=(\be,0)$.
 Since $\Lambda^{[Q]}$ is skew-symmetric, we may assume $p \le s$. Then we have
$\al \not\preceq_{[Q]} \be$, and the conditions $p-s \le d(i,j)$ and
$  \xi_j-\xi_i  \le d(i,j)$ hold  obviously.
We shall apply Proposition~\ref{prop: NnKR} with $p'=\xi_i$ and $s'=\xi_j$.
Then, together with Lemma~\ref{lem: iso tori}~\eqref{it: la al}, we obtain
\eqn
\ucalN( \tm^{(i)}[p,\xi_i] , \tm^{(j)}[s,\xi_j]) &&=
\Bigl(\cox{Q}^{(\xi_{i} -p)/2 + 1}\varpi_{i} + \cox{Q}^{(\xi_{i} - p')/2}\varpi_{i},\;
\cox{Q}^{(\xi_{j}-s)/2 + 1}\varpi_{j} - \cox{Q}^{(\xi_j -s')/2}\varpi_j\Bigr)
\\
&&=
(\la_\al+\varpi_i,\la_\be-\varpi_j) = -\Lambda_{\be\al}.     \qquad\qquad \qquad \qquad\qquad \qedhere
\eneqn
\end{proof}

\begin{definition} Let $Q=(\Dynkin,\xi)$ be a Dynkin quiver.
We call the image  of $\calA_\nu(\La^{[Q]}, \tsfB^{[Q]}) $ under $\Uppsi_Q $ the \emph{quantum virtual Grothendieck ring associated to $Q$}
and denote it by $\frakK_{\ttq,Q}$.
\end{definition}

\begin{remark}
For a Dynkin quiver $Q$ of simply-laced type, $\frakK_{\ttq,Q}$ and $\frakK_{\ttq,Q}\vert_{\ttq=1}$ are known as the \emph{quantum} Grothendieck ring
and  the Grothendieck ring of a certain subcategory  $\scrC_Q$ of modules over the  quantum affine algebra associated to $Q$.
We expect  that $\frakK_{\ttq,Q}$ (resp.\ $\frakK_{\ttq,Q}\vert_{\ttq=1}$)   is contained in $\frakK_{\ttq}$ (resp.\ $\frakK_{\ttq}\vert_{\ttq=1}$) when $Q$ is a Dynkin quiver of type $BCFG$.
\end{remark}

\section{Compatible pairs} \label{Sec: bi dec quivers and conjecture} \label{sec: conj}
 In this section, we first give a generalization of Proposition
~\ref{prop: compatible pair2} as a conjecture.  Then we prove the conjecture under the condition~\eqref{cond:din add} below.

\smallskip
Let $\tw=\sseq{i_1,\ldots,i_r}$ be a sequence of elements of $\Dynkin_0$.
We set $J \seteq [1,r]$.
We define the
matrix $\tsfB^\tw$ and $\La^\tw$
as in Proposition~\ref{prop: compatible pair2}.
Namely, we set
\eqn
j^+&&\seteq\max\st{s\in J\mid j<s,\;i_s=i_j}\cup\st{r+1}\qt{for any $j\in J$,}\\
J_f &&\seteq \{j \in J \mid j^+ = r+1 \},\quad  J_e \seteq J\setminus J_f,\\
w_{\le t}&&\seteq s_{i_1}\cdots s_{i_t}\qt{for any $t\in J$,}
\eneqn
and
the $J \times J_e$-matrix
$\tsfB^{\tw}=(b_{s,t})_{s \in J,\; t\in J_e}$ is defined by
\begin{align*}
b_{s,t} = \bc
1 &\text{ if } t=s^+, \\
-1 &\text{ if } s=t^+, \\
\sfc_{i_s,i_t} & \text{if $s<t<s^+<t^+$,} \\
-\sfc_{i_s,i_t} &\text{if $t<s<t^+<s^+$,} \\
0 &\text{otherwise,}
\ec
\end{align*}
and the  skew-symmetric $J \times J$-matrix $\Lambda^{\tw} =
(\Lambda_{s,t})_{s,t \in J}$ is defined by
\begin{align}\label{eq: Lambdast}
\Lambda_{s,t} = (\varpi_{i_s}-w_{\le
s}\varpi_{i_s},\varpi_{i_t}+w_{\le t}\varpi_{i_t})\quad \text{ for }
s<t.
\end{align}
Note that $\tsfB^{\tw}$ is skew-symmetrizable by $\diag(d_{i_k})_{k\in J_e}$.

Consider the condition on $\tw$:
\eq
&&\parbox[t]{70ex}{$s_{i_a}\cdots s_{i_b}$ has length $b-a+1$
for any $a,b\in \Z$ such that $1\le a\le b\le r$ and $1+b-a\le\ell(w_0)$.}\label{cond:din}
\eneq

\begin{conjecture} \label{conj: compatible}
For any $\tw$ satisfying \eqref{cond:din},
the pair $(\Lambda^\tw,\tsfB^\tw)$ is compatible, i.e.,
$$\Lambda^\tw\tsfB^\tw=(-2d_{i_s}\delta(s=t) )_{s\in J,\;t\in J_e}.$$
\end{conjecture}

\subsection{$Q$-adapted case}
  In this subsection, we prove Conjecture~\ref{conj: compatible} when
$\tw$ is a $Q$-adapted sequence for some Dynkin quiver $Q$
but \eqref{cond:din} is not assumed.

For a sequence $\tw$ in~\eqref{cond:din}, we assume
the following condition:
\eq
&&\parbox[t]{70ex}{there exists a Dynkin quiver $Q$
such that $\sseq{i_k}_{1\le k\le r}$ is $Q$-adapted, i.e., $i_k$ is a source of
$s_{i_{k-1}}\cdots s_{i_1}Q$ for any $k$ with $1\le k\le r$.
}\label{cond:din add}
\eneq

\emph{Throughout this subsection, we fix $\tw$ satisfying
~\eqref{cond:din add}. }

For each $j \in\Dynkin_0$, we set
$$   n_j \seteq   \big\vert \st{k\in\Z\mid \text{$1 \le k \le r$ and $i_k =j$}} \big\vert  \in \Z_{\ge 0} \sqcup \{ \infty \}.$$

We set
$$
\Gamma_0^{[\tw]}\seteq \{ (j,p) \in \hDynkin_0 \mid \xi_j -2 n_j < p \le \xi_j \}.
$$
For $k\in[1,r]$, set
$$p_k=\xi_{i_k}-2\times \big\vert\st{s\in\Z\mid 1\le s<k, i_s=i_k}\big\vert.$$
Then, by Lemma~\ref{lem:admissible},
$\sseq{(i_k,p_k)}_{1\le k\le r}$ is a compatible reading of
$\Gamma^{[\tw]}_0$ \eqref{eq: comreading}. Hence
$k\leftrightarrow(i_k,p_k)$ gives a one-to-one correspondence
$$ \eta\cl \Gamma^{[\tw]}_0\isoto J= \{ 1,2,\ldots,r\}.$$
We set
$$(\Gamma_0^{[\tw]})_e\seteq\eta^{-1}(J_e)
=\st{(j,p)\in\hDynkin_0\mid \xi_j-2n_j-2<p\le\xi_j}.$$
Then the exchange matrix $\tsfB^{\tw}=(b_{s,t})_{s \in J,\; t\in J_e}$
and the skew-symmetric matrix $J \times J$-matrix $\Lambda^{\tw}$
can be re-written as follows:
\eq
&&b'_{(i,p),(j,s)}\seteq b_{\eta(i,p),\eta(j,s)}= \bc
(-1)^{\delta(s>p)}\sfc_{i,j} & \text{ if } |p-s|=1 \text{ and } d(i,j)=1, \\
(-1)^{\delta(s>p)}  & \text{ if } |p-s|=2 \text{ and } i = j, \\
0 &\text{otherwise,}
\ec
\label{eq:Bd}
\eneq
for $((i,p), (j,s) ) \in    \Gamma^{[\tw]}_0\times \bl\Gamma^{[\tw]}_0\br_e$,
and
\eq
&& \Lambda'_{(i,p),(j,s)}\seteq\Lambda_{\eta(i,p),\eta(j,s)} =
\ucalN\bl\tm^{(i)}[p,\xi_i],\; \tm^{(j)}[s, \xi_j]\br
 =\hs{-1ex}\sum_{\substack{x:\;p\lec x\lec \xi_i\\[.2ex] y:\;s\lec y\lec \xi_j}}
\hs{-1ex}\ucalN(i,x;j,y)\label{eq:Lamd}
\eneq
for $(i,p), (j,s) \in \Gamma^{[\tw]}_0$.
Here $u\lec v$ means that $u\le v$ and $u\equiv v\bmod 2$.

Indeed, \eqref{eq:Bd} is obvious, and
\eqref{eq:Lamd} is obtained from
Lemma~\ref{lem:admissible}\,\eqref{item:svarp} and Proposition~\ref{prop: NnKR}
as follows: assuming $p\le s$, we have
\eqn
\Lambda'_{(i,p),(j,s)}=-\Lambda'_{(j,s),(i,p)}&&=
-\bl \varpi_{j}-\tauQ^{(\xi_j-s)/2}\varpi_j,
\varpi_{i}+\tauQ^{(\xi_i-p)/2}\varpi_i\br\\
&&=\ucalN\bl\tm^{(i)}[p,\xi_i],\; \tm^{(j)}[s, \xi_j]\br.
\eneqn

\begin{theorem} \label{thm: compatible conj}
Assume that $\tw$ satisfies \eqref{cond:din add}.
Then, we have
$$\Lambda^{\tw}\tsfB^{\tw}= (-2d_{i_s}\delta(s=t))_{(s,t)\in J\times J_e}.$$
\end{theorem}

\begin{proof}
Let $(i,p)\in\Gamma^{[\tw]}_0$ and $(j,s)\in(\Gamma^{[\tw]}_0)_e$.
By \eqref{eq:Bd}, we have
\eqn
\bl \Lambda^{\tw}\tsfB^{\tw} \br_{ \eta(i,p),\eta(j,s) }&& =
\delta(s\not=\xi_j)\Lambda'_{(i,p),(j,s+2)} - \Lambda'_{(i,p),(j,s-2)}\\
&&\hs{10ex} + \sum_{k:\; d(k,j)=1} \lan h_k,\al_j \ran \left(\delta(s\not=\xi_k+1) \Lambda'_{(i,p),(k,s+1)} -
\Lambda'_{(i,p),(k,s-1)} \right).
\eneqn
Note that
\eqn
&&\delta(s\not=\xi_j)\Lambda'_{(i,p),(j,s+2)} - \Lambda'_{(i,p),(j,s-2)} \allowdisplaybreaks\\
&&\hs{5ex}=
\sum_{t:\;s+2\lec t\lec\xi_i}\ucalN(\tm^{(i)}[p,\xi_i],\tX_{j,t})-
\sum_{t:\;s-2\lec t\lec\xi_i}\ucalN(\tm^{(i)}[p,\xi_i],\tX_{j,t}) \allowdisplaybreaks\\
&&\hs{5ex}=-\ucalN(\tm^{(i)}[p,\xi_i],\tX_{j,s-2})
-\ucalN(\tm^{(i)}[p,\xi_i],\tX_{j,s})
\eneqn
On the other hand, we have
\eqn
\ucalN(\tm^{(i)}[p,\xi_i],\tX_{j,s})&&=\hs{-3ex}\sum_{t:\;p\lec t\lec\xi_i}
\ucalN(i,t;j,s)
=\sum_{t:\;p\lec t\lec\xi_i}\bl\teta_{i,j}(t-s-1)-\teta_{i,j}(t-s+1)\br\\
&&=\sum_{t:\;p-1\lec t\lec\;\xi_i-1}\teta_{i,j}(t-s)-
\sum_{t:\;p+1\lec t\lec\;\xi_i+1}\teta_{i,j}(t-s)\\
&&=\teta_{i,j}(p-s-1)-\teta_{i,j}(\xi_i-s+1).
\eneqn
Hence we obtain
\eqn
&&\delta(s\not=\xi_j)\Lambda'_{(i,p),(j,s+2)} - \Lambda'_{(i,p),(j,s-2)}\\
&&\hs{5ex}=-\Bigl(\teta_{i,j}(p-s+1)-\teta_{i,j}(\xi_i-s+3)+\teta_{i,j}(p-s-1)-\teta_{i,j}(\xi_i-s+1)\Bigr)\\
&&\hs{5ex}=-\teta_{i,j}(p-s-1)-\teta_{i,j}(p-s+1)+\teta_{i,j}(\xi_i-s+1)
+\teta_{i,j}(\xi_i-s+3).
\eneqn
Similarly
\eqn
&&\delta(s\not=\xi_k+1) \Lambda_{(i,p),(k,s+1)} - \Lambda_{(i,p),(k,s-1)} \\
&&\hs{5ex}=\sum_{t:\;s+1\lec t\lec\xi_k}\ucalN(\tm^{(i)}[p,\xi_i],\tX_{k,t})
-\sum_{t:\;s-1\lec t\lec\xi_k}\ucalN(\tm^{(i)}[p,\xi_i],\tX_{k,t}) \allowdisplaybreaks\\
&&\hs{5ex}=-\ucalN(\tm^{(i)}[p,\xi_i],\tX_{k,s-1}) \allowdisplaybreaks\\
&&\hs{5ex} =-\teta_{i,k}(p-s)+\teta_{i,k}(\xi_i-s+2).
\eneqn

Hence we obtain

\eqn
\bl \Lambda^{\tw}\tsfB^{\tw} \br_{ \eta(i,p),\eta(j,s) }
&&=-\teta_{i,j}(p-s-1)-\teta_{i,j}(p-s+1)+\teta_{i,j}(\xi_i-s+1)
+\teta_{i,j}(\xi_i-s+3) \allowdisplaybreaks\\
&&\hs{5ex}+ \sum_{k; d(k,j)=1} \ang{h_k,\al_j}\bl- \teta_{i,k}(p-s) +\teta_{i,k}(\xi_i-s+2)\br \allowdisplaybreaks\\
&&\eqs2d_i\; \delta(i=j) \bl-\delta(p-s=0)+\delta(\xi_i-s+2=0)\br \allowdisplaybreaks\\
&&=-2d_i\;\delta\bl(i,p)=(j,s)\br,
\eneqn
where $\eqs$ follows from Corollary~\ref{cor: additive application}.
\end{proof}

\begin{remark}
When (a) $\Dynkin$ is simply-laced, (b) $Q=(\Dynkin,\xi)$ has a
sink-source orientation, that is  $\xi_i \in \{0,1\}$ for all $i \in
\Dynkin_0$, and (c) $r = \infty$, the above theorem is proved
in~\cite[Proposition 5.1.1]{B20} (see also \cite[Proposition 5.26 in
Arxiv version (arXiv:2007.03159v1)]{FO21}).

\end{remark}

\appendix

\section{$\tde_{i,j}(t)$ for $E_7$ and $E_8$} \label{appeA: tde}

\subsection{$E_7$}   Here is the list of $\tde_{i,j}(t)$ for $E_7$.
\begin{align*}
&\tde_{1,1}(t) = t^1+t^{7}+t^{11}+t^{17},  \hspace{17.4ex}\tde_{1,2}(t) = t^4+t^{8}+t^{10}+t^{14}, \allowdisplaybreaks \\
&\tde_{1,3}(t)=t^2+t^{6}+t^{8}+t^{10}+t^{12}+t^{16}, \ \hspace{7ex} \tde_{1,4}(t)=t^3+t^{5}+t^{7}+2t^{9}+t^{11}+t^{13}+t^{15},\allowdisplaybreaks \\
&\tde_{1,6}(t)=t^5+t^{7}+t^{11}+t^{13}, \hspace{17.8ex}  \tde_{1,7}(t)=t^6+t^{12}, \allowdisplaybreaks\\
&\tde_{2,2}(t)=t^1+t^{5}+t^{7}+t^{9}+t^{11}+t^{13}+t^{17}, \hspace{3.2ex}  \tde_{2,3}(t)=\tde_{1,4}(t)  \allowdisplaybreaks\\
&\tde_{2,4}(t)=t^2+t^{4}+2t^{6}+2t^{8}+2t^{10}+2t^{12}+t^{14}+t^{16}, \allowdisplaybreaks\\
&\tde_{2,5}(t)=t^3+t^{5}+2t^{7}+t^{9}+2t^{11}+t^{13}+t^{15}, \  \tde_{2,7}(t)=t^{5}+t^{9}+t^{13}, \allowdisplaybreaks\\
&\tde_{2,6}(t)=t^{4}+t^{6}+t^{8}+t^{10}+t^{12}+t^{14}, \allowdisplaybreaks\\
&\tde_{3,3}(t) = t^1+t^{3}+t^{5}+2t^{7}+2t^{9}+2t^{11}+t^{13}+t^{15}+t^{17}, \allowdisplaybreaks\\
&\tde_{3,4}(t) = t^2+2t^{4}+2t^{6}+3t^{8}+3t^{10}+2t^{12}+2t^{14}+t^{16}, \allowdisplaybreaks\\
&\tde_{3,5}(t) = t^{3}+2t^{5}+2t^{7}+2t^{9}+2t^{11}+2t^{13}+t^{15}, \allowdisplaybreaks\\
&\tde_{3,6}(t) = t^{4}+2t^{6}+t^{8}+t^{10}+2t^{12}+t^{14}, \hspace{6ex}   \tde_{3,7}(t) = t^{5}+t^{7}+t^{11}+t^{13}, \allowdisplaybreaks\\
&\tde_{4,4}(t) = t^1+2t^{3}+3t^{5}+4t^{7}+4t^{9}+4t^{11}+3t^{13}+2t^{15}+t^{17}, \allowdisplaybreaks\\
&\tde_{4,5}(t) = t^2+2t^{4}+3t^{6}+3t^{8}+3t^{10}+3t^{12}+2t^{14}+t^{16}, \allowdisplaybreaks\\
&\tde_{4,6}(t) =\tde_{3,5}(t),  \hspace{29ex}   \tde_{4,7}(t) = t^{4}+t^{6}+t^{8}+t^{10}+ t^{12}+ t^{14}, \allowdisplaybreaks\\
&\tde_{5,5}(t) = t^1+t^{3}+2t^{5}+2t^{7}+3t^{9}+2t^{11}+2t^{13}+t^{15}+t^{17}, \allowdisplaybreaks\\
&\tde_{5,6}(t) = t^2+t^{4}+t^{6}+2t^{8}+2t^{10}+t^{12}+t^{14}+t^{16}, \quad \tde_{5,7}(t) = t^{3}+t^{7}+t^{11}+t^{15}, \allowdisplaybreaks\\
&\tde_{6,6}(t) = t^1+t^{3}+t^{7}+2t^{9}+t^{11}+t^{15}+t^{17},  \hspace{4ex}   \tde_{6,7}(t) = t^2+t^{8}+t^{10}+t^{16}, \allowdisplaybreaks\\
&\tde_{7,7}(t) = t^1+t^{9}+t^{17} \hspace{10.65ex}  \text{ and }  \hspace{10.45ex}  \tde_{i,j}(t) =   \tde_{j,i}(t).
\end{align*}

\subsection{$E_8$}   Here is the list of $\tde_{i,j}(t)$ for $E_8$.
\begin{align*}
\tde_{1,1}(t) &= t^1+t^{7}+t^{11}+t^{13}+t^{17}+t^{19}+t^{23}+t^{29}, \allowdisplaybreaks\\
\tde_{1,2}(t) &= t^4+t^{8}+t^{10}+t^{12}+t^{14}+t^{16}+t^{18}+t^{20}+t^{22}+t^{26}, \allowdisplaybreaks\\
\tde_{1,3}(t) &= t^2+t^{6}+t^{8}+t^{10}+2t^{12}+t^{14}+t^{16}+2t^{18}+t^{20}+t^{22}+t^{24}+t^{28}, \allowdisplaybreaks\\
\tde_{1,4}(t) &= t^3+t^{5}+t^{7}+2t^{9}+2t^{11}+2t^{13}+2t^{15}+2t^{17}+2t^{19}+2t^{21}+t^{23}+t^{25}+t^{27}, \allowdisplaybreaks\\
\tde_{1,5}(t) &= t^4+t^{6}+t^{8}+2t^{10}+t^{12}+2t^{14}+2t^{16}+t^{18}+2t^{20}+t^{22}+t^{24}+t^{26}, \allowdisplaybreaks\\
\tde_{1,6}(t) &= t^{5}+t^{7}+t^{9}+t^{11}+t^{13}+2t^{15}+t^{17}+t^{19}+t^{21}+t^{23}+t^{25}, \allowdisplaybreaks\\
\tde_{1,7}(t) &= t^{6}+t^{8}+t^{12}+t^{14}+t^{16}+t^{18}+t^{22}+t^{24},  \allowdisplaybreaks\\
\tde_{1,8}(t) &= t^{7}+t^{13}+t^{17}+t^{23}, \allowdisplaybreaks\\
\tde_{2,2}(t) &= t^1+t^{5}+t^{7}+t^{9}+2t^{11}+t^{13}+2t^{15}+t^{17}+2t^{19}+t^{21}+t^{23}+t^{25}+t^{29}, \allowdisplaybreaks\\
\tde_{2,3}(t) &= t^3+t^{5}+t^{7}+2(t^{9}+t^{11}+t^{13}+t^{15}+t^{17}+t^{19}+t^{21})+t^{23}+t^{25}+t^{27}, \allowdisplaybreaks\\
\tde_{2,4}(t) &= t^2+t^{4}+2t^{6}+2t^{8}+3(t^{10}+t^{12}+t^{14}+t^{16}+t^{18}+t^{20})+2t^{22}+2t^{24}+t^{26}+t^{28}, \allowdisplaybreaks\\
\tde_{2,5}(t) &= t^3+t^{5}+2t^{7}+2t^{9}+2t^{11}+3t^{13}+3t^{15}+3t^{17}+2t^{19}+2t^{21}+2t^{23}+t^{25}+t^{27}, \allowdisplaybreaks\\
\tde_{2,6}(t) &= t^4+t^{6}+2t^{8}+t^{10}+2t^{12}+2t^{14}+2t^{16}+2t^{18}+t^{20}+2t^{22}+t^{24}+t^{26}, \allowdisplaybreaks\\
\tde_{2,7}(t) &=\tde_{1,6}(t), \qquad  \qquad  \qquad \tde_{2,8}(t) = t^{6}+t^{10}+t^{14}+t^{16}+t^{20}+t^{24},  \allowdisplaybreaks\\
\tde_{3,3}(t) &= t^1+t^{3}+t^{5}+2t^{7}+2t^{9}+3t^{11}+3t^{13}+2t^{15}+3t^{17}+3t^{19}+2t^{21}+2t^{23} \allowdisplaybreaks\\
                         & \hspace{2ex}+t^{25}+t^{27}+t^{29}, \allowdisplaybreaks\\
\tde_{3,4}(t) &= t^2+2(t^{4}+t^{6})+3t^{8}+4(t^{10}+t^{12}+t^{14}+t^{16}+t^{18}+t^{20})+3t^{22}+2(t^{24}+t^{26})+t^{28}, \allowdisplaybreaks\\
\tde_{3,5}(t) &= t^3+2t^{5}+2t^{7}+3t^{9}+3t^{11}+3t^{13}+4t^{15}+3t^{17}+3t^{19}+3t^{21}+2t^{23}+2t^{25}+t^{27}, \allowdisplaybreaks\\
\tde_{3,6}(t) &= t^4+2t^{6}+2t^{8}+2t^{10}+2t^{12}+3t^{14}+3t^{16}+2t^{18}+2t^{20}+2t^{22}+2t^{24}+t^{26}, \allowdisplaybreaks\\
\tde_{3,7}(t) &= t^{5}+2t^{7}+t^{9}+t^{11}+2t^{13}+2t^{15}+2t^{17}+t^{19}+t^{21}+2t^{23}+t^{25}, \allowdisplaybreaks\\
\tde_{3,8}(t) &= t^{6}+t^{8}+t^{12}+t^{14}+t^{16}+t^{18}+t^{22}+t^{24}, \allowdisplaybreaks\\
\tde_{4,4}(t) &= t^1+2t^{3}+3t^{5}+4t^{7}+5t^{9}+6t^{11}+6t^{13}+6t^{15}+6t^{17}+6t^{19}+5t^{21}\allowdisplaybreaks\\
                         & \hspace{2ex}+4t^{23} +3t^{25}+2t^{27}+t^{29}, \allowdisplaybreaks\\
\tde_{4,5}(t) &= t^2+2t^{4}+3t^{6}+4t^{8}+4t^{10}+5t^{12}+5t^{14}+5t^{16}+5t^{18}+4t^{20}+4t^{22}\allowdisplaybreaks\\
                         & \hspace{2ex}+3t^{24}+2t^{26}+t^{28}, \allowdisplaybreaks\\
\tde_{4,6}(t) &= t^3+2t^{5}+3t^{7}+3t^{9}+3t^{11}+4t^{13}+4t^{15}+4t^{17}+3t^{19}+3t^{21}+3t^{23}+2t^{25}+t^{27}, \allowdisplaybreaks\\
\tde_{4,7}(t) &= \tde_{3,6}(t), \ \
\tde_{4,8}(t) = t^{5}+t^{7}+t^{9}+t^{11}+t^{13}+2t^{15}+t^{17}+t^{19}+t^{21}+t^{23}+t^{25}, \allowdisplaybreaks\\
\tde_{5,5}(t) &= t^1+t^{3}+2t^{5}+3t^{7}+3t^{9}+4t^{11}+4t^{13}+4t^{15}+4t^{17}+4t^{19}+3t^{21}\allowdisplaybreaks\\
                         & \hspace{2ex}+3t^{23}+2t^{25}+t^{27}+t^{29}, \allowdisplaybreaks\\
\tde_{5,6}(t) &= t^2+t^{4}+2t^{6}+2t^{8}+3t^{10}+3t^{12}+3t^{14}+3t^{16}+3t^{18}+3t^{20}+2t^{22}\allowdisplaybreaks\\
                         & \hspace{2ex}+2t^{24}+t^{26}+t^{28}, \allowdisplaybreaks\\
\tde_{5,7}(t) &= t^3+t^{5}+t^{7}+2t^{9}+2t^{11}+2t^{13}+2t^{15}+2t^{17}+2t^{19}+2t^{21}+t^{23}+t^{25}+t^{27}, \allowdisplaybreaks\\
\tde_{5,8}(t) &= t^{4}+t^{8}+t^{10}+t^{12}+t^{14}+t^{16}+t^{18}+t^{20}+t^{22}+t^{26}, \allowdisplaybreaks\\
\tde_{6,6}(t) &= t^1+t^{3}+t^{5}+t^{7}+2t^{9}+3t^{11}+2t^{13}+2t^{15}+2t^{17}+3t^{19}+2t^{21}\allowdisplaybreaks\\
                         & \hspace{2ex}+t^{23}+t^{25}+t^{27}+t^{29}, \allowdisplaybreaks\\
\tde_{6,7}(t) &= t^2+t^{4}+t^{8}+2t^{10}+2t^{12}+t^{14}+t^{16}+2t^{18}+2t^{20}+t^{22}+t^{26}+t^{28}, \allowdisplaybreaks\\
\tde_{6,8}(t) &= t^3+t^{9}+t^{11}+t^{13}+t^{17}+t^{19}+t^{21}+t^{27}, \allowdisplaybreaks\\
\tde_{7,7}(t) &= t^1+t^{3}+t^{9}+2t^{11}+t^{13}+t^{17}+2t^{19}+t^{21}+t^{27}+t^{29}, \allowdisplaybreaks\\
\tde_{7,8}(t) &= t^2+t^{10}+t^{12}+t^{18}+t^{20}+t^{28}, \allowdisplaybreaks\\
\tde_{8,8}(t) &= t^1+t^{11}+t^{19}+t^{29} \qquad\qquad
\qt{and} \qquad\qquad\quad
 \tde_{i,j}(t) =   \tde_{j,i}(t).
\end{align*}

\end{document}